\documentclass[12pt,reqno,twoside]{amsart}
\usepackage[asymmetric,top=3.5cm,bottom=4.3cm,left=3.1cm,right=3.1cm]{geometry}
\usepackage{xcolor}
\colorlet{BLUE}{blue}
\usepackage{amsmath,amsfonts,amsthm,mathrsfs,amssymb,bm,cite}
\usepackage{graphicx,epstopdf,booktabs,array}
\usepackage{paralist,verbatim,caption,subcaption}
\usepackage{cases}
\usepackage[textsize=tiny,textwidth=1in,shadow,backgroundcolor=yellow!40]{todonotes}

\newtheorem{thm}{Theorem}[section]

\theoremstyle{definition}
\newtheorem{defn}{Definition}[section]
\theoremstyle{remark}

\numberwithin{equation}{section}

\newcommand{\RR}{\mathbb{R}}
\newcommand{\CC}{\mathbb{C}}

\title[Geometric Structure of Transmission Eigenfunctions]{\bf  On vanishing and localizing of transmission eigenfunctions near singular points: a numerical study}

\author{Eemeli Bl{\aa}sten}
\address{HKUST Jockey Club Institute for Advanced Study, Hong Kong University of Science and Technology, Hong Kong SAR.}
\email{iaseemeli@ust.hk}

\author{Xiaofei Li}
\address{Department of Mathematics, Southern University of Science and Technology, Shenzhen, China.}
\email{xiaofeilee@hotmail.com}

\author{Hongyu Liu}
\address{Department of Mathematics, Hong Kong Baptist University,Kowloon Tong, Hong Kong SAR.\vspace*{-4mm}}
\address{\vspace*{-4mm}and}
\address{HKBU Institute of Research and Continuing Education, Virtual University Park, Shenzhen, P. R. China.}
\email{hongyu.liuip@gmail.com}

\author{Yuliang Wang}
\address{Department of Mathematics, Hong Kong Baptist University,Kowloon Tong, Hong Kong SAR.\vspace*{-4mm}}
\address{\vspace*{-4mm}and}
\address{HKBU Institute of Research and Continuing Education, Virtual University Park, Shenzhen, P. R. China.}
\email{yuliang@hkbu.edu.hk}

\begin{document}
\maketitle

\begin{abstract}
This paper is concerned with the intrinsic geometric structure of interior transmission eigenfunctions arising in wave scattering theory. We numerically show that the aforementioned geometric structure can be very delicate and intriguing. The major findings can be roughly summarized as follows. We say that a point on the boundary of the perturbation is singular, we mean that the surface tangent is discontinuous there. The interior transmission eigenfunction then vanishes near singular points where the  interior angle is less than $\pi$, whereas the interior transmission eigenfunction localizes near the singular point if its interior angle is bigger than $\pi$. Furthermore, we show that the vanishing and blowup orders are inversely proportional to the interior angle of the singular point: the sharper the corner, the higher the convergence order. Our results are first of its type in the spectral theory for transmission eigenvalue problems, and the existing studies in the literature concentrate more on the intrinsic properties of the transmission eigenvalues instead of the transmission eigenfunctions. Due to the finiteness of computing resources, our study is by no means exclusive and complete. We consider our study only in a certain geometric setup including corner, curved corner and edge singularities. Nevertheless, we believe that similar results hold for more general singularities and rigorous theoretical justifications are much desirable. Our study enriches the spectral theory for transmission eigenvalue problems. We also discuss its implication to inverse scattering theory.

\medskip\medskip


\end{abstract}
\medskip

\section{Introduction}

Let $D$ be a bounded domain in $\RR^d,d=2,3$ and $n \in L^\infty(D)$ be a measurable function such that $n\equiv\hspace*{-4mm}\backslash \ 1$ in $D$. Consider the following interior transmission problem
\begin{equation}\label{eq:1}
\begin{cases}
  \ \  \Delta u + k^2 n u = 0 & \ {\rm in}\ \ D,  \\
  \ \  \Delta u_0 + k^2 u_0 = 0 & \ {\rm in}\ \ D,  \\
  \ \ u-u_0\in H_0^2(D),
    \end{cases}
\end{equation}
where the Sobolev space $H_0^2(D)$ is defined as the completion of test functions in $D$ by the standard $H^2(D)$-norm. For Lipschitz domains
\begin{align*}
  H_0^2(D) = \left\{v \in H^2(D): v =0, ~ \partial_\nu v = 0 \ {\rm on}\ \partial D \right\}.
\end{align*}
with $\nu$ signifying the unit normal vector directed into the exterior of $D$.

\begin{defn}\label{def:1}
  A value $k \in \CC, k \neq 0$ for which the transmission problem \eqref{eq:1} has nontrivial solutions $(u,u_0) \in L^2(D) \times L^2(D)$ such that $u-u_0 \in H_0^2(D)$ is called an interior transmission eigenvalue associated with $(D; n)$. The nontrivial solutions $(u,u_0)$ are called the corresponding interior transmission eigenfunctions.
\end{defn}
In what follows, for simplicity, we call $k$ and $(u, u_0)$ in Definition~\ref{def:1}, respectively, the transmission eigenvalue and transmission eigenfunctions. It is our objective of this paper to demonstrate and investigate certain intrinsic geometric properties of the transmission eigenfunctions through numerical experiments.

The interior transmission eigenvalue problem arises in inverse scattering theory and has become a very important area of research. In fact, both linear sampling method \cite{Colton--Kirsch} and the factorization method \cite{Kirsch--Grinberg}, two important reconstruction methods for inverse scattering problems, succeed at wavenumbers that are not transmission eigenvalues.  The study of the transmission eigenvalue problem is mathematically interesting and challenging since it is a type of non-elliptic and non self-adjoint problem. Recently, the transmission eigenvalue problems were also connected to invisibility cloaking \cite{JL,LWZ}. In the literature, the existing studies concentrate more on the intrinsic properties of the transmission eigenvalues, including the existence, discreteness and infiniteness, see for example \cite{BYY,BY,CGH10,CPS,CKP,S}. However, there are few results concerning the intrinsic properties of the transmission eigenfunctions. Only very recently there is some started appearing. An example is \cite{YXB} where the authors have applied the theory of nodal sets to the transmission eigenfunctions of a radially symmetric refractive index to determine the latter from the eigenvalues.

In a recent article \cite{EL} by the first and third author of this paper, it is shown that if the perturbation $1-n$ possesses a corner on its support, and the interior angle of the corner is less than $\pi$, then the associated transmission eigenfunctions must vanish near the corner under a certain generic condition. In this article, with the help of numerical methods, we shall show that the geometric structure of the transmission eigenfunction can be very delicate and intriguing, and the study in \cite{EL} together with the current one opens up a new direction of research on the spectral theory of transmission eigenvalue problems.

We first recall two important topics in the classical spectral theory for Dirichlet/Neumann Laplacian: nodal sets and eigenfunction localization. The former is the set of points in the domain where the eigenfunction vanishes. For the latter, an eigenfunction is said to be localized if most of its $L^2$-energy concentrates in a subdomain which is a fraction of the total domain. The vanishing and localizing of Dirichlet/Neumann Laplacian eigenfunctions have many important applications, both in pure and applied areas, and they still remain active in the mathematical investigation; see \cite{GN} for a more recent survey. The aim of our study is to show that the transmission eigenfunctions also possess the intrinsic vanishing and localizing behaviours.

The major numerical findings of the present paper can be roughly summarized as follows. If there is a singular point on the support of the underlying perturbation to the refractive index, then the transmission eigenfunction vanishes near the singular point if its interior angle is less than $\pi$, whereas the transmission eigenfunction localizes near the sungular point if its interior angle is bigger than $\pi$. By a singular point we mean a point where the tangent at the boundary is discontinuous.  Furthermore, we show that the vanishing and blowup orders are proportional to the interior angle of the singularity by inversion: the sharper the angle, the higher the convergence order.
Due to the limitedness of the computing resources, our study is by no means exclusive and complete. We consider our study only in a certain geometric setup including corner, curved corner and edge singularities. Nevertheless, we believe that similar results hold for more general singularities and rigorous theoretical justifications are very desirable. Our study enriches the spectral theory for transmission eigenvalue problems. More relevant discussions of our study shall be given in Section~\ref{sect:con}.

The rest of the paper is organized as follows. In the next section, we briefly review some existing results on the transmission eigenvalues and eigenfunctions. In Section \ref{sec:numerical-method}, we present the numerical method to solve the transmission eigenvalue problem. Sections \ref{sec:numerical-experiment-2D} and \ref{sec:numerical-experiment-3D} are devoted to the main results on numerically showing the vanishing and localizing properties, as well as the convergence orders in two and three dimensions. We conclude our study in Section~\ref{sect:con} with some discussions.

\section{Preliminaries on transmission eigenvalue problem}

Throughout the rest of the paper, we assume that $n(x)$, $x\in D$ is real-valued. It is remarked that we allow the presence of complex eigenvalues in our study.

We first collect some theoretical results for the interior transmission eigenvalue problem \eqref{eq:1}. Denote by
\begin{align*}
  n_* = \inf_{x \in D} n(x), \quad n^* = \sup_{x \in D} n(x)
\end{align*}
the essential infimum and supremum of the refractive index. The following theorem shows the existence of interior transmission eigenvalues \cite{CGH10}.
\begin{thm}
 Let $n\in L^{\infty}(D)$ satisfies either one of the following assumptions
  \begin{enumerate}
  \item $1+\alpha \leq n_*\leq n \leq n^* < \infty$ for some constants $\alpha>0$,
  \item $0 < n_* \leq n\leq n^* \leq 1-\beta$ for some constants $\beta>0$,
  \end{enumerate}
then there exists an infinite set of interior transmission eigenvalues with $+\infty$ as the only accumulation point.
\end{thm}

In the numerical computation of the transmission eigenvalue problem, it is desirable to have an estimate of the lower bound of the first positive transmission eigenvalue. Recall from \cite{CGH10} that
\begin{thm} \label{thm:bounds}
  Let $R$ be the radius of the smallest ball containing $D$. Denote by $\lambda_*$ and $\lambda^*$ the first positive transmission eigenvalue corresponding to the ball of radius $1$ with refractive index $n \equiv n_*$ and $n \equiv n^*$, respectively. Let $\lambda_1$ be the first Dirichlet eigenvalue for $-\Delta$ in $D$. Denote by $k_1$ the first positive transmission eigenvalue corresponding to $D$ and refractive index $n$.
  \begin{enumerate}
  \item If $1+\alpha \leq n_*\leq n \leq n^* < \infty$ for some constant $\alpha>0$, then
    \begin{align*}
      k_1 \geq \max \left( \frac{\lambda^*}{R}, \sqrt{\frac{\lambda_1}{n^*}} \right).
    \end{align*}
    \item If $0 < n_*\leq n \leq n^* < 1-\beta$ for some constant $\beta>0$, then
      \begin{align*}
        k_1 \geq \max \left( \frac{\lambda_*}{R}, \sqrt{\lambda_1} \right).
      \end{align*}
  \end{enumerate}
\end{thm}

Our numerical experiments indicate the existence of complex transmission eigenvalues, but this has not been established theoretically in general. The following theorem shows the non-existence of purely imaginary transmission eigenvalues \cite{CMS10}.
\begin{thm}
  If $n>1$ or $n<1$ almost everywhere in $D$, then there exists no purely imaginary transmission eigenvalues.
\end{thm}

Next, we give a more definite description of the vanishing and localizing of transmission eigenfunctions.

\begin{defn}\label{def:2}
Assume that $k\in\mathbb{C}$ is a transmission eigenvalue, then there exist $u_0,u\in L^2(D)$ such that
\begin{align*}
   & \Delta u + k^2 n u = 0  \quad\, {\rm in}\ D,  \\
   & \Delta u_0 + k^2 u_0 = 0  \quad {\rm in}\ D,  \\
  & u-u_0\in H^2_0(D),~\|u_0\|_{L^2(D)}=1.
\end{align*}
Let $P\in\partial D$ be a point and $B_r(P)$ be a ball of radius $r\in\mathbb{R}_+$ centred at $P$. Set $D_r(P):=B_r(P)\cap D$. Assume that
\begin{equation}\label{eq:ii1}
\|1-n\|_{L^\infty(D_r(p))}\geq \epsilon_0, \quad\epsilon_0\in\mathbb{R}_+.
\end{equation}
Then we say that vanishing occurs near $P$ if
$$\lim_{r\rightarrow +0}\frac{1}{\sqrt{|D_r(P)|}}\|u_0(x)\|_{L^2(D_r(P))}=0,$$
whereas we say that localizing occurs near $P$ if
$$\lim_{r\rightarrow +0}\frac{1}{\sqrt{|D_r(P)|}}\|u_0(x)\|_{L^2(D_r(P))}=+\infty,$$
where $|D_r(P)|$ signifies the area or volume of the region $D_r(P)$, respectively, in two and three dimensions.
\end{defn}

If $P$ is the vertex of a corner with an interior angle less than $\pi$, the vanishing of the transmission eigenfunctions has been rigorously verified in \cite{EL}. Our observation of the vanishing and localizing near a corner point comes from the corner scattering study in \cite{EL1,BPS,PSV,HSV,EH1,EH2}. An important consequence of the study in \cite{EL1,BPS,PSV,EH1} is the fact that the transmission eigenfunction $u_0$ cannot be analytically extended across a corner point to form an entire solution to the Hemholtz equation, $\Delta u_0+k^2 u_0=0$ in $\mathbb{R}^d$. However, we note that due to the interior regularity, the transmission eigenfunction $u_0$ is always analytic away from the corner point. Hence, heuristically, the failure of the analytic extension may indicate that $u_0$ either vanishes or blows up when approaching the corner point. Clearly, the failure of the analytical extension should also hold across any irregular point on the support of $1-n$, and hence we conjecture that the vanishing or localizing behaviours of the transmission eigenfunctions would occur near any singular point on the support of the underlying. Theoretical proof of such a conjecture is fraught with significant difficulties. Indeed, the proof in \cite{EL} of the vanishing of the transmission eigenfunction in the special case near a corner with an interior angle less than $\pi$ already involves much technical analysis and advanced tools. In the next two sections, we conduct extensive numerical experiments to verify the aforementioned conjecture, and hopefully, the numerical results can also inspire the rigorous theoretical proof in the future.

\section{Finite element method for the transmission problem} \label{sec:numerical-method}
If $D$ is smooth and $n$ is a constant, then one may solve the transmission eigenvalue problem \eqref{eq:1} by using the integral equation formulation \cite{CH13,Kleefeld13,ZSX16}. If $D$ is nonsmooth or $n$ is nonconstant, then the finite element method is more appropriate \cite{CMS10,HSSZ16,MS12,Nedelec80}. We shall use the continuous finite element approximation as proposed in \cite{CMS10}, which we describe briefly in the sequel. Multiplying the second equation in  \eqref{eq:1} by a test function $\phi \in H_0^1(D)$ and integrating by parts, we obtain
\begin{align}
  \label{eq:6}
  (\nabla u_0, \nabla \phi) - k^2 (u_0, \phi) = 0 \quad \forall\, \phi \in H_0^1(D).
\end{align}
Multiplying the first equation in \eqref{eq:1} by a test function $\psi \in H_0^1(D)$ and integrating by parts, we obtain
\begin{align}
  \label{eq:7}
  (\nabla u, \nabla \psi) - k^2 (nu, \psi) = 0 \quad \forall\, \psi \in H_0^1(D).
\end{align}
To enforce the boundary conditions in \eqref{eq:1} weakly, we multiply it by a test function $\psi \in H^1(D)$ and integrate by parts to obtain
\begin{align}
  \label{eq:8}
  (\nabla u, \nabla \psi) - k^2 (nu, \psi) = (\nabla u_0, \nabla \psi) - k^2 (u_0, \psi) \quad \forall\, \psi \in H^1(D).
\end{align}
Note that \eqref{eq:7} is already implied by \eqref{eq:6} and \eqref{eq:8}. Hence a variational formulation for \eqref{eq:1} is:
\begin{center}
  \begin{minipage}{0.8\linewidth}
    to find $(u,u_0) \in H^1(D) \times H^1(D)$ satisfying \eqref{eq:6} and
    \eqref{eq:8}, together with the essential boundary condition $u=u_0$ on
    $\partial D$.
  \end{minipage}
\end{center}

Let $\mathcal{T}$ be a regular triangular or tetrahedral mesh of $D$. Denote by $V$ the finite element subspace of $H^1(D)$ consisting of piecewise linear functions on each element of $\mathcal{T}$ and $V_0=V \cap H_0^1(D)$. Denote by $\{\eta_j: j=1,\cdots,I\}$ the set of linear Lagrange basis functions for $V_0$ and $\{\eta_j: j=1,\cdots,I\} \cup \{\zeta_j: j=1,\cdots,J\}$ the set of linear Lagrange basis functions for $V$, respectively. Let
\begin{align*}
  u_0^h = \sum_{j=1}^I u_0^{(j)} \eta_j + \sum_{j=1}^J w^{(j)} \zeta_j, \quad
  u^h = \sum_{j=1}^I u^{(j)} \eta_j + \sum_{j=1}^J w^{(j)} \zeta_j
\end{align*}
represent the finite element approximations of $u_0$ and $u$, respectively. Note that the essential boundary condition $u=u_0$ on $\partial D$ is already enforced in the above representation. Let
\begin{align*}
  X=\left[ u_0^{(1)}, \cdots, u_0^{(I)}, u^{(1)}, \cdots, u^{(I)}, w^{(1)}, \cdots, w^{(J)} \right] \in \CC^{2I+J}
\end{align*}
denote the vector of unknown coefficients, then the finite element discretization of \eqref{eq:6} and \eqref{eq:8} can be written as
\begin{align}
  \label{eq:5}
  AX=k^2BX,
\end{align}
where
\begin{align*}
  A=
  \begin{bmatrix}
    S_1^{II} & 0 & S_1^{IB} \\[5pt]
    0 & S_1^{II} & S_1^{IB} \\[5pt]
    S_1^{BI} & -S_1^{BI} & 0
  \end{bmatrix},
\quad
B=
                           \begin{bmatrix}
                             M_2^{II} & 0 & M_2^{IB} \\[5pt]
                             0 & M_1^{II} & M_1^{IB} \\[5pt]
                             M_2^{BI} & -M_1^{BI} & M_2^{BB}-M_1^{BB}
                           \end{bmatrix}
\end{align*}
with the block stiffness and mass matrices given by
\begin{align*}
  S_1^{II}(i,j) &= (\nabla \eta_j, \nabla \eta_i), &\quad 1 \leq i \leq I, \ 1 \leq j \leq I, \\
  S_1^{IB}(i,j) &= (\nabla \zeta_j, \nabla \eta_i), &\quad 1 \leq i \leq I, \ 1 \leq j \leq J, \\
  S_1^{BI}(i,j) &= (\nabla \eta_j, \nabla \zeta_i), &1 \leq i \leq J, \ 1 \leq j \leq I, \\
  M_1^{II}(i,j) &= (\eta_j, \eta_i), &\quad 1 \leq i \leq I, \ 1 \leq j \leq I, \\
  M_1^{IB}(i,j) &= (\zeta_j, \eta_i), &\quad 1 \leq i \leq I, \ 1 \leq j \leq J, \\
  M_1^{BI}(i,j) &= (\eta_j, \zeta_i), &1 \leq i \leq J, \ 1 \leq j \leq I, \\
  M_1^{BB}(i,j) &= (\zeta_j, \zeta_i),& \quad 1 \leq i \leq J, \ 1 \leq j \leq J,\\
  M_2^{II}(i,j) &= (n\eta_j, \eta_i), &\quad 1 \leq i \leq I, \ 1 \leq j \leq I, \\
  M_2^{IB}(i,j) &= (n\zeta_j, \eta_i),& \quad 1 \leq i \leq I, \ 1 \leq j \leq J, \\
  M_2^{BI}(i,j) &= (n\eta_j, \zeta_i), &1 \leq i \leq J, \ 1 \leq j \leq I, \\
  M_2^{BB}(i,j) &= (n\zeta_j, \zeta_i),& \quad 1 \leq i \leq J, \ 1 \leq j \leq J.
\end{align*}
The assembly of the above matrices is implemented with the FEM package given by \cite{AV15}. The non-Hermitian generalized eigenvalue problem \eqref{eq:5} is then solved by the {\tt sptarn} function in MATLAB, which is based on the Arnoldi algorithm with spectral transformation. The lower bound for the search interval of the eigenvalues may be chosen according to Theorem \ref{thm:bounds}.

\section{Numerical experiments: two dimension}\label{sec:numerical-experiment-2D}
\subsection{Example: Equilateral triangle}
\begin{figure}[t]
  \centering
  \begin{subfigure}[b]{0.32\textwidth}
    \includegraphics[width=\textwidth]{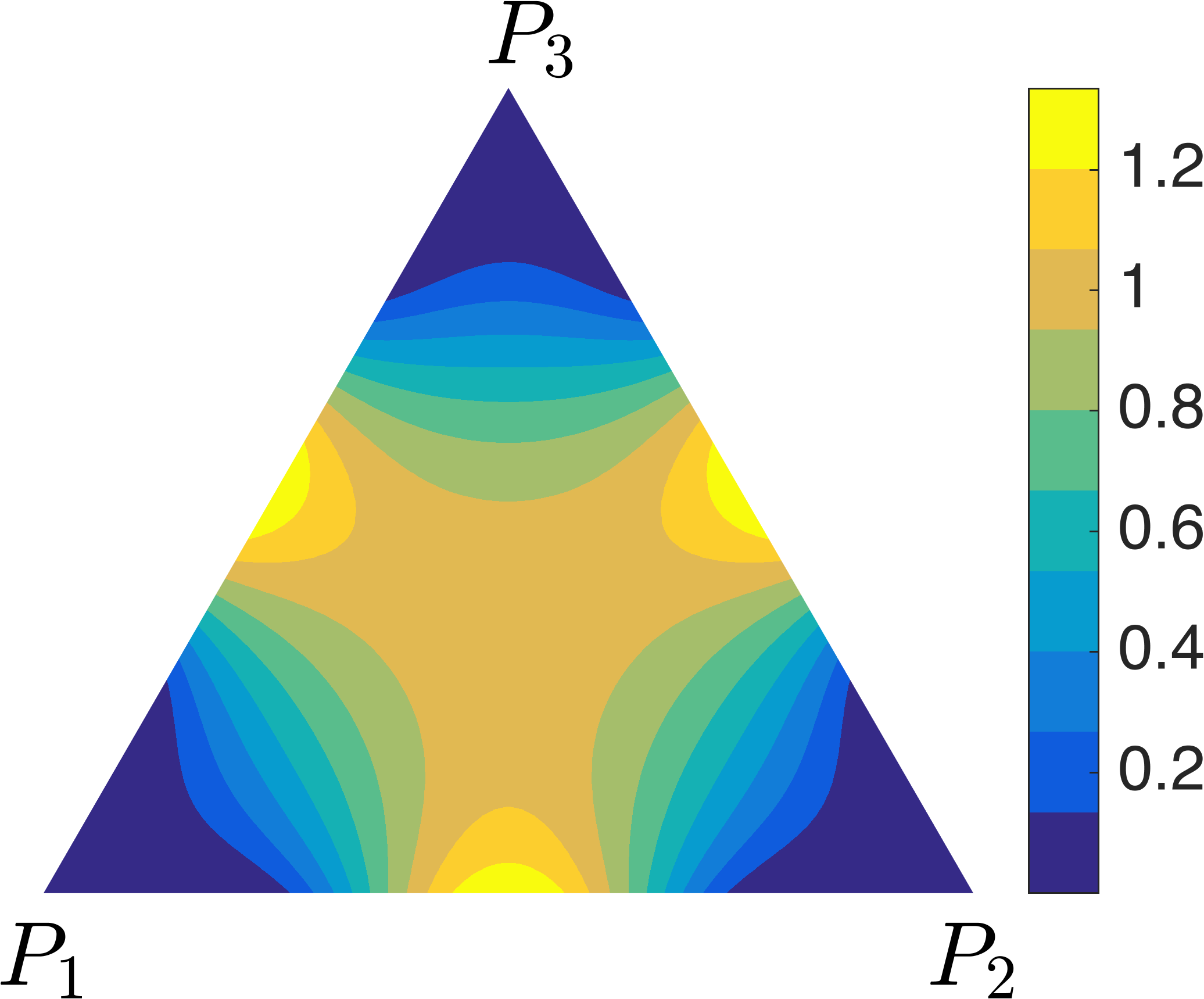}
    \caption{}
  \end{subfigure}
  \hfill
  \begin{subfigure}[b]{0.32\textwidth}
    \includegraphics[width=\textwidth]{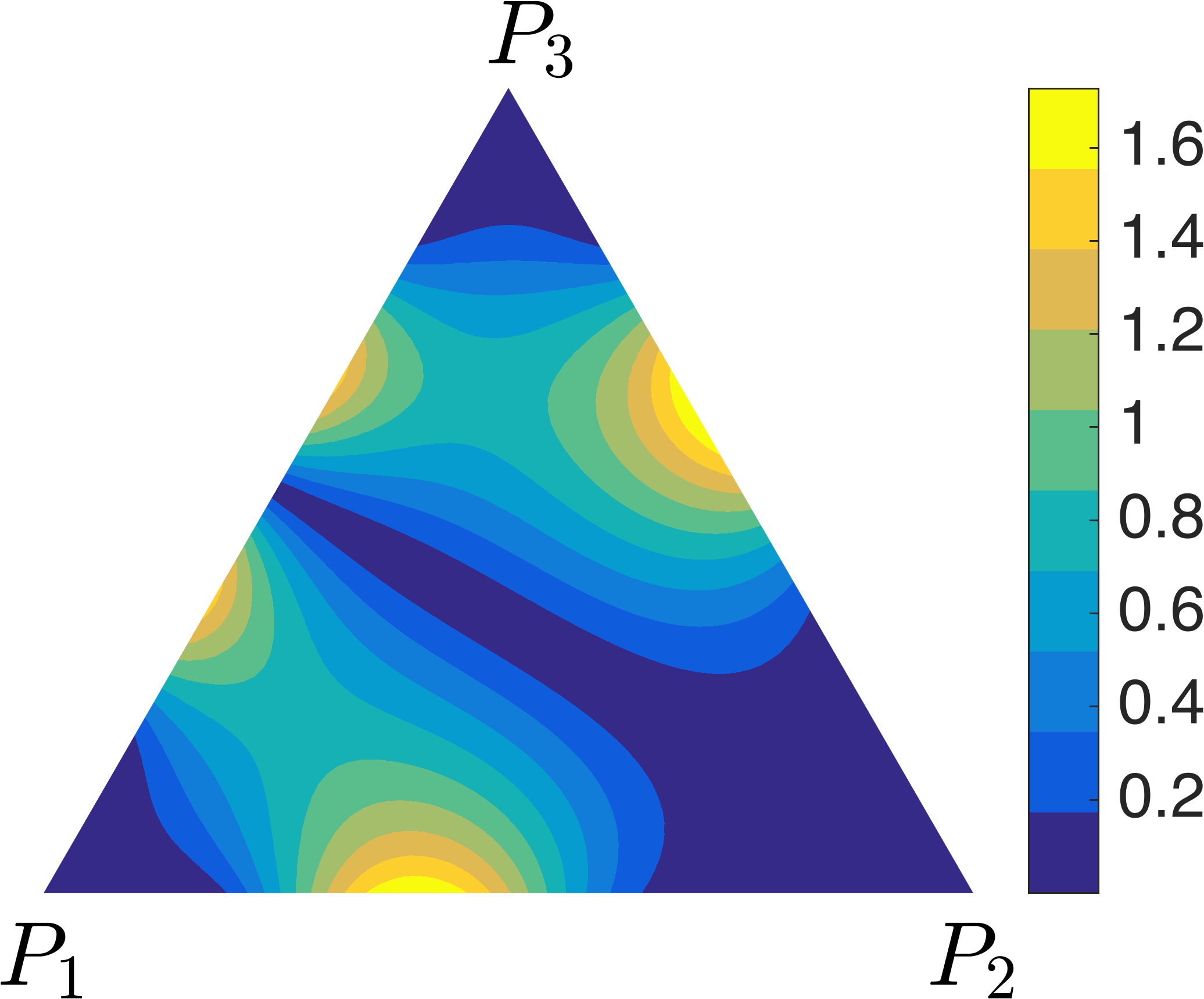}
    \caption{}
  \end{subfigure}
  \hfill
  \begin{subfigure}[b]{0.32\textwidth}
    \includegraphics[width=\textwidth]{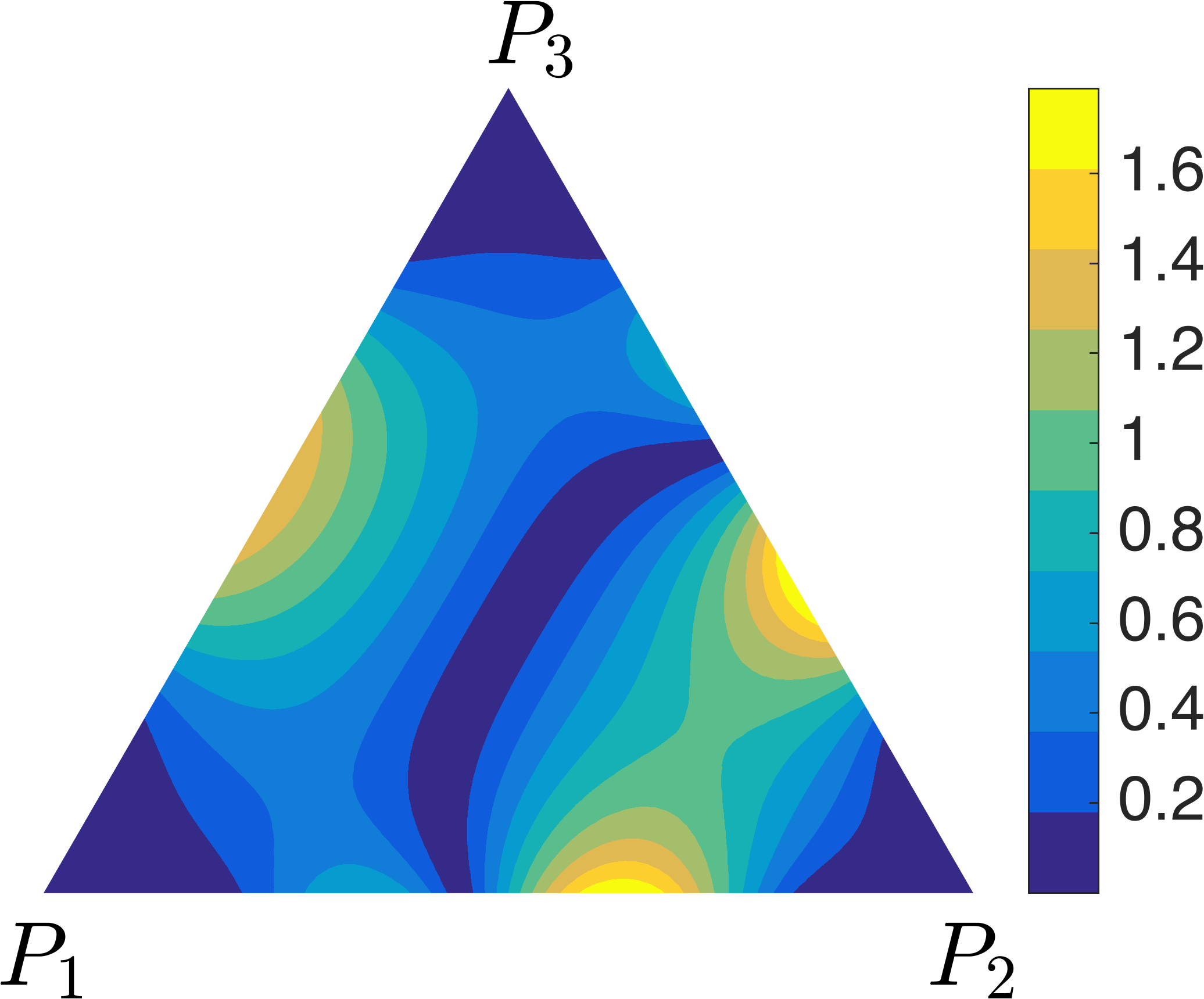}
    \caption{}
  \end{subfigure}
  \hfill
  \begin{subfigure}[b]{0.32\textwidth}
    \includegraphics[width=\textwidth]{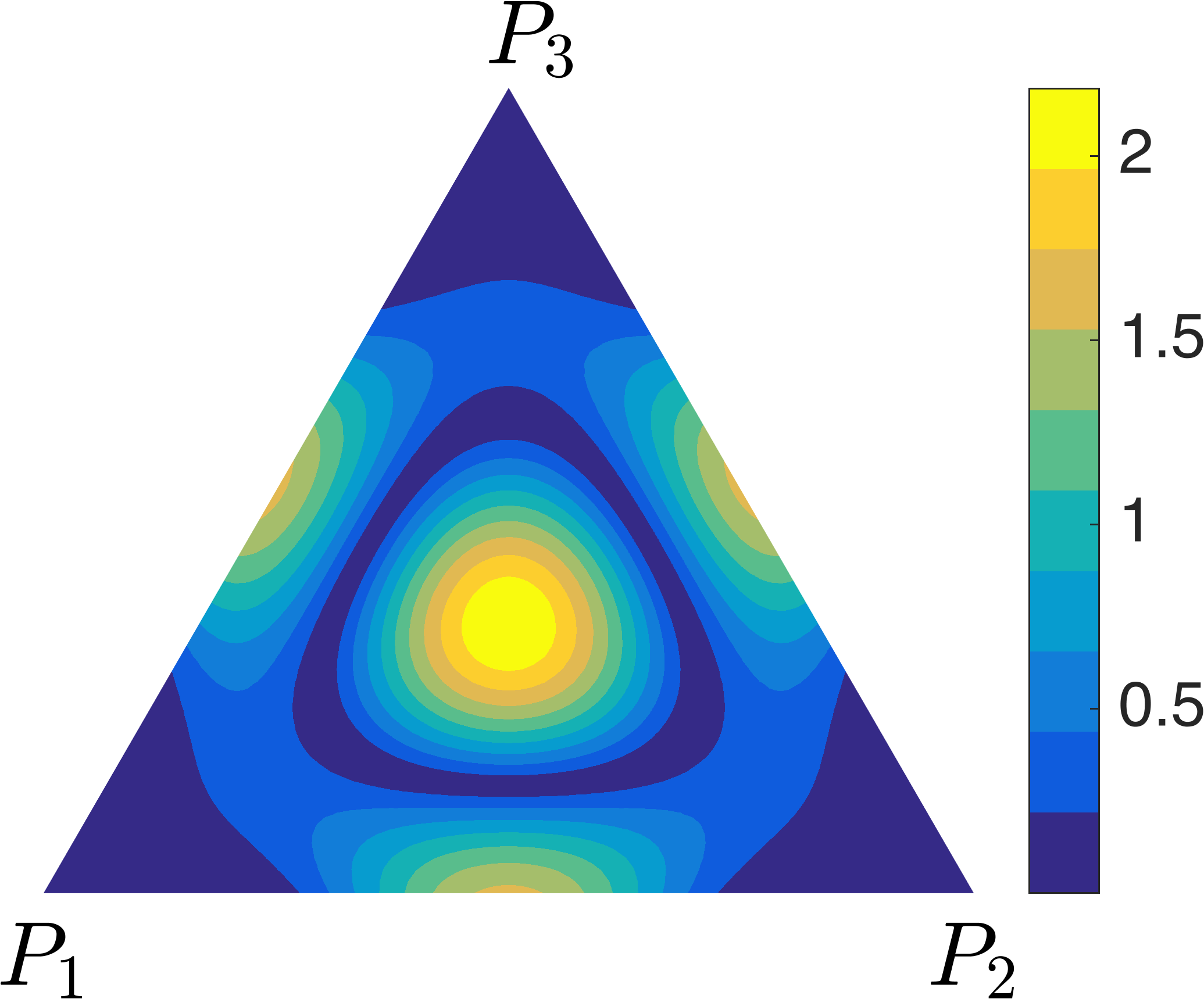}
    \caption{}
  \end{subfigure}
  \hfill
  \begin{subfigure}[b]{0.32\textwidth}
    \includegraphics[width=\textwidth]{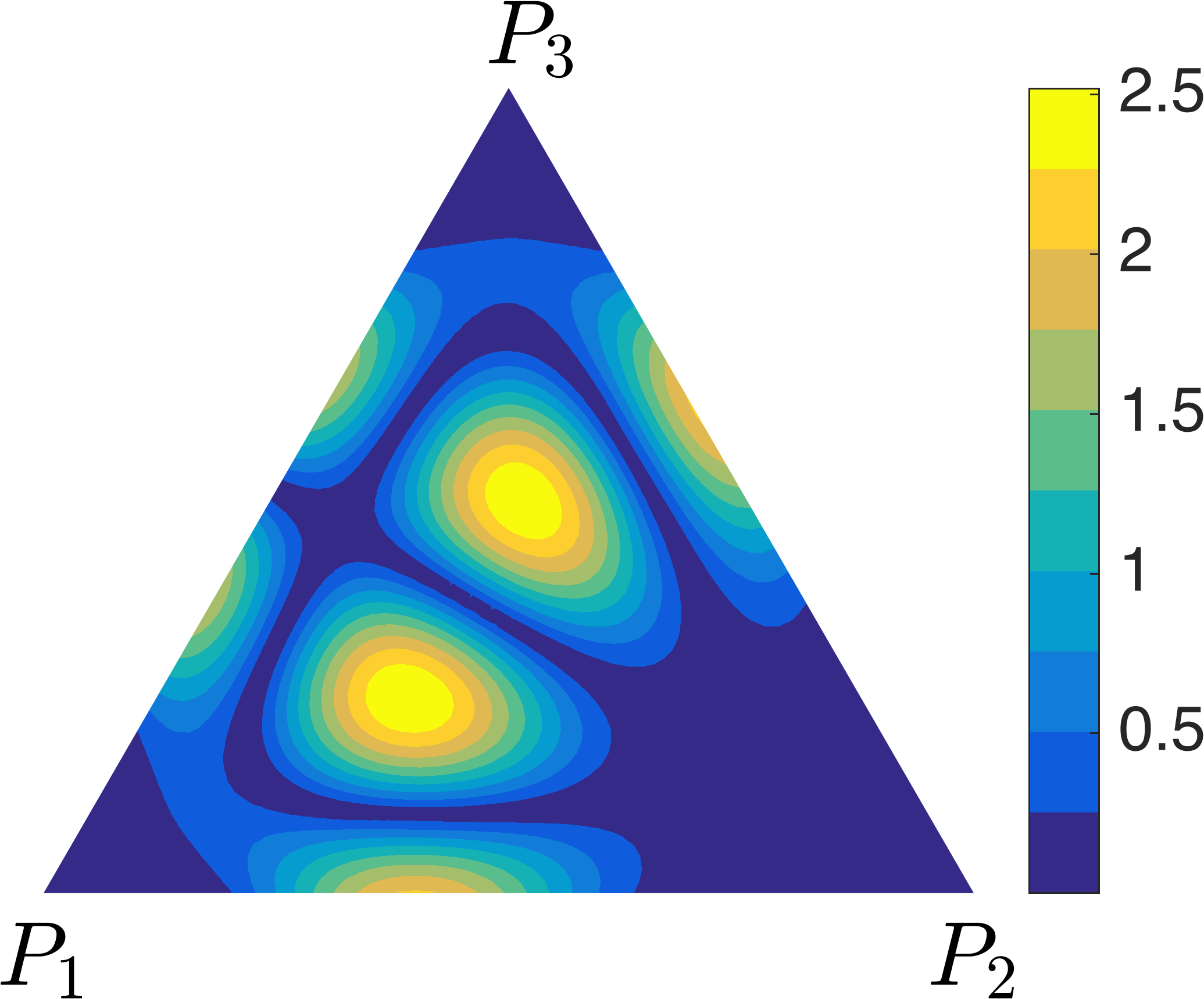}
    \caption{}
  \end{subfigure}
  \hfill
  \begin{subfigure}[b]{0.32\textwidth}
    \includegraphics[width=\textwidth]{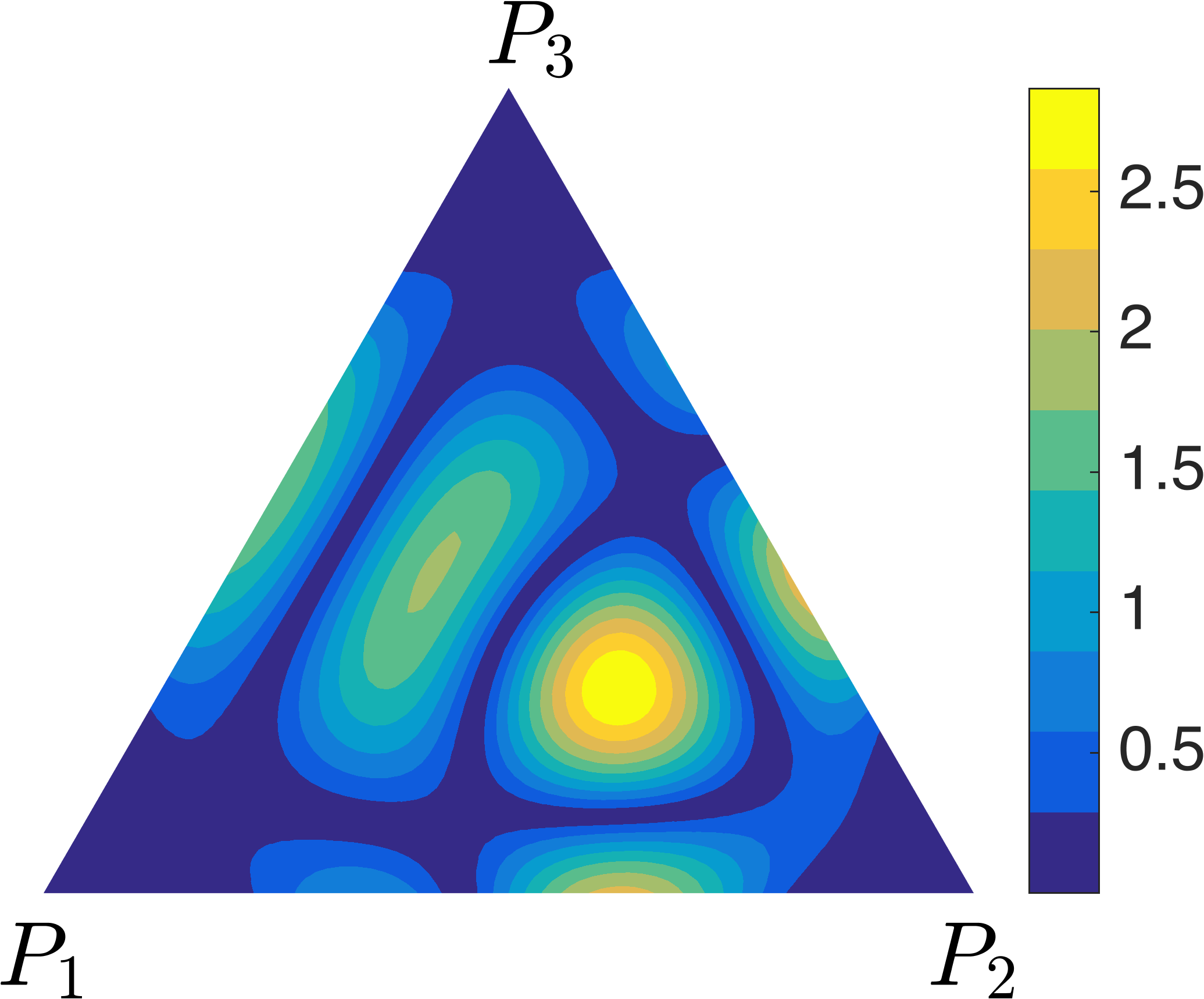}
    \caption{}
  \end{subfigure}
  \caption{The magnitude of transmission eigenfunctions for the equilateral triangle with refractive index $n=16$. (A) $|u_0^{(1)}|$; (B) $|u_0^{(2)}|$; (C) $|u_0^{(3)}|$; (D) $|u^{(1)}|$; (E) $|u^{(2)}|$; (F) $|u^{(3)}|$.}
  \label{fig:equilateral_triangle}
\end{figure}

Let $D$ be the equilateral triangle with vertices at $P_1=(-1,0), P_2=(1,0)$ and $P_3=(0,\sqrt{3})$. Let the refractive index $n=16$. The first three real transmission eigenvalues are found to be $1.5748, 1.9807$ and $1.9807$ with multiplicities. Denote by $u_0^{(i)}$ and $u^{(i)}$ the $i$-th transmission eigenfunctions corresponding to the $i$-th transmission eigenvalue. The magnitude of the transmission eigenfunctions for the first three eigenvalues are shown in Figure \ref{fig:equilateral_triangle}.
We observe that each transmission eigenfunction tends to zero towards every corner point of $D$.

\begin{figure}[t]
  \centering
  \begin{subfigure}[b]{0.32\textwidth}
    \includegraphics[width=\textwidth]{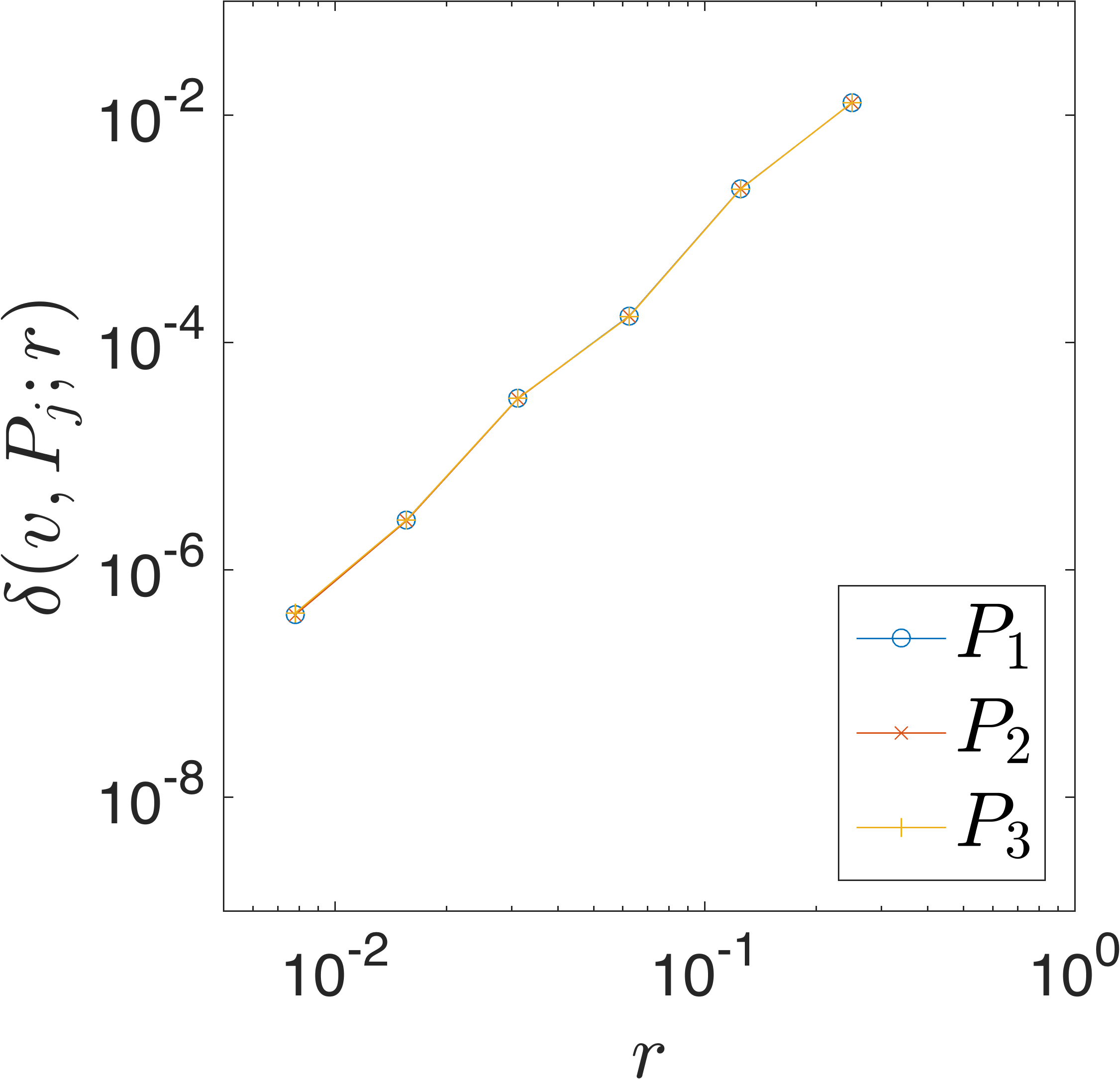}
    \caption{}
  \end{subfigure}
  \hfill
  \begin{subfigure}[b]{0.32\textwidth}
    \includegraphics[width=\textwidth]{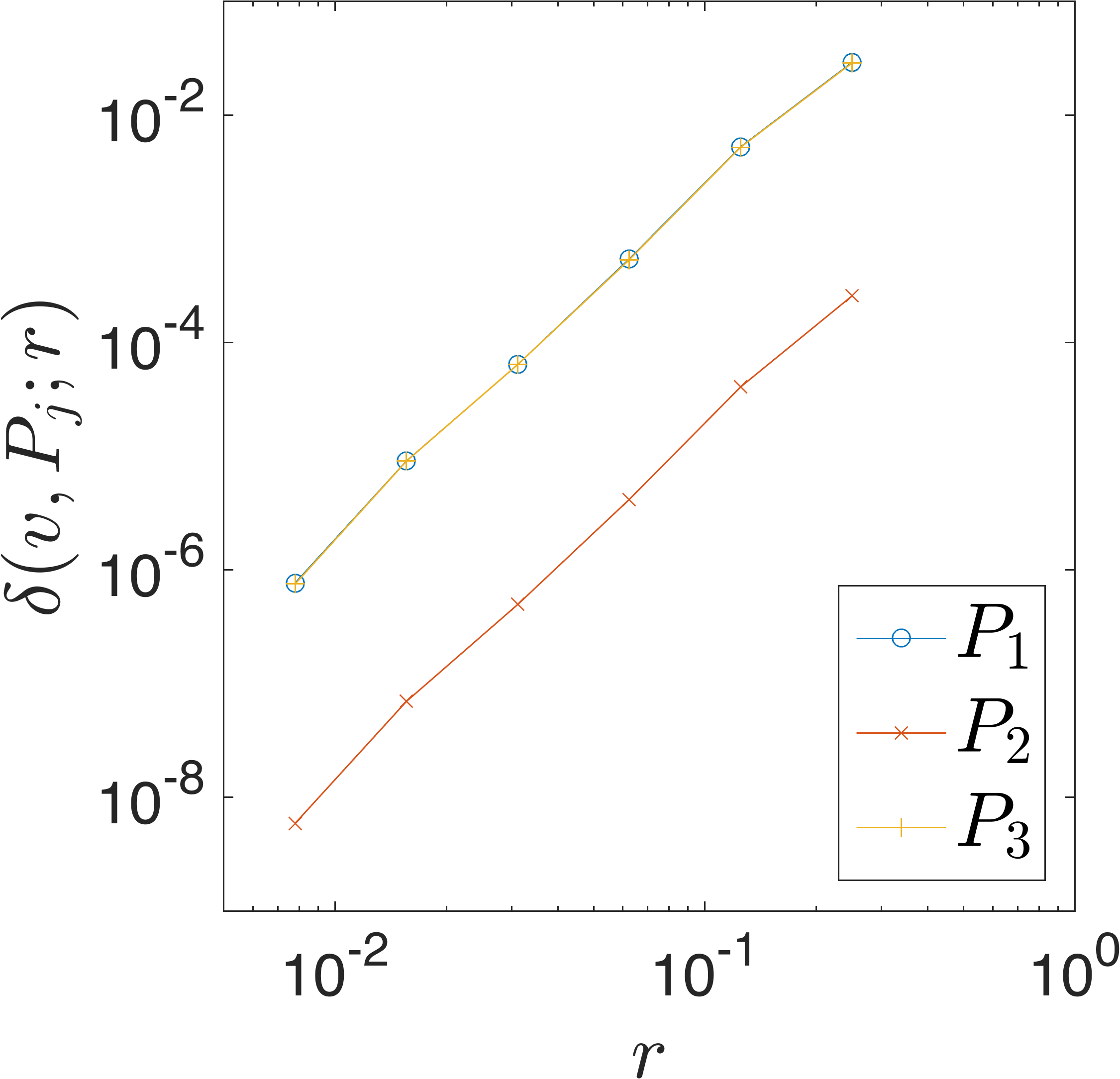}
    \caption{}
  \end{subfigure}
  \hfill
  \begin{subfigure}[b]{0.32\textwidth}
    \includegraphics[width=\textwidth]{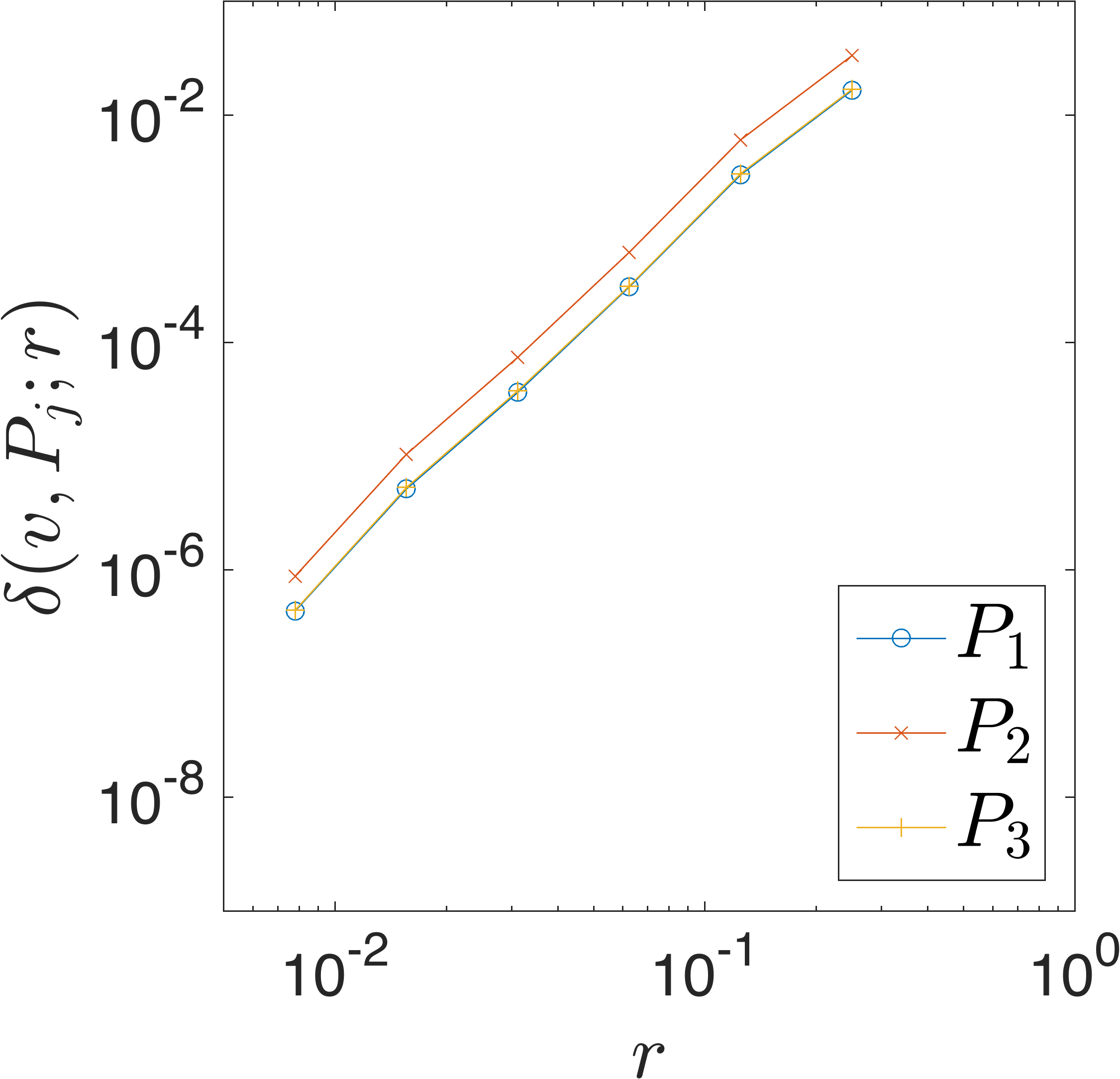}
    \caption{}
  \end{subfigure}
  \hfill
  \begin{subfigure}[b]{0.32\textwidth}
    \includegraphics[width=\textwidth]{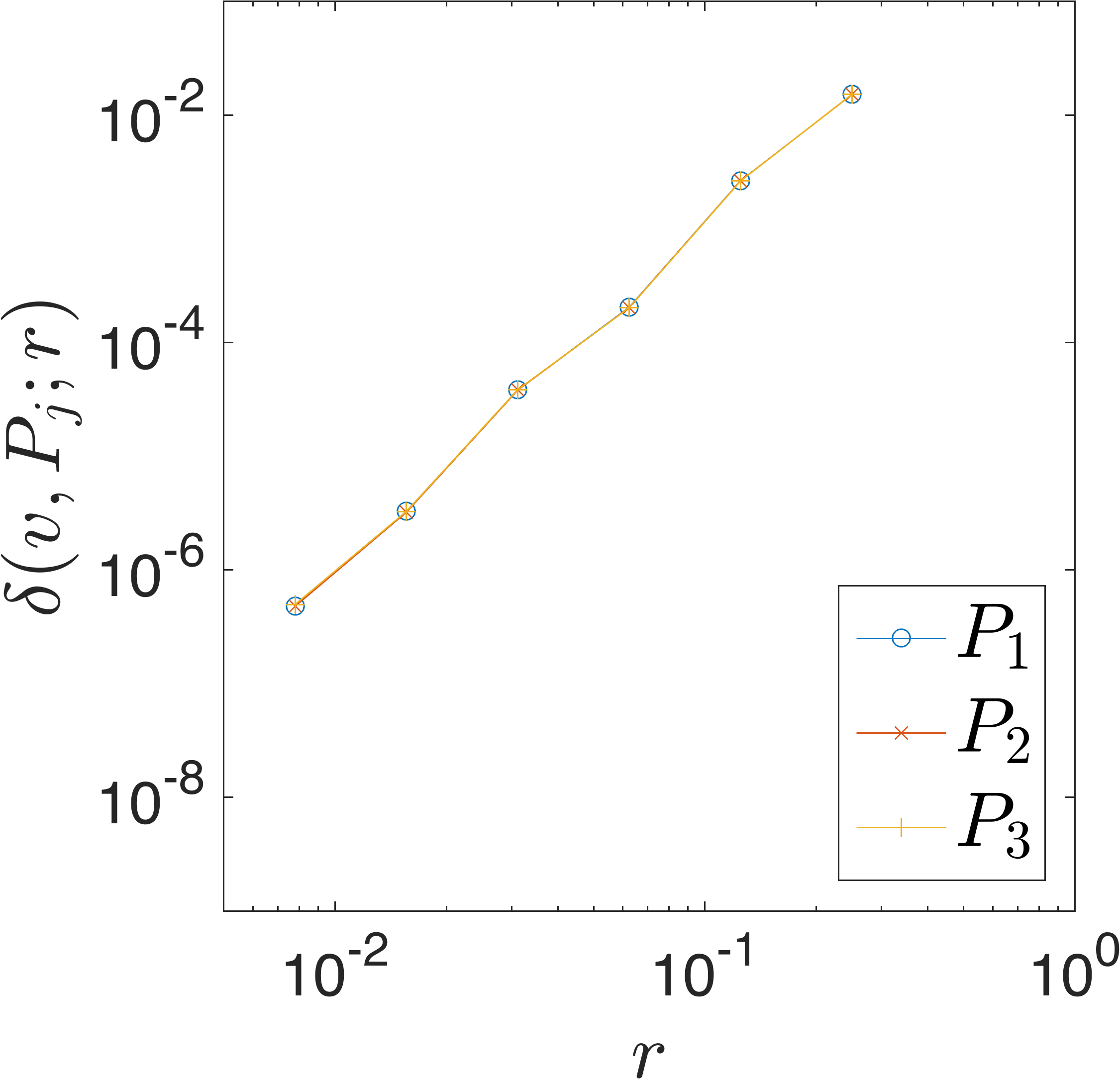}
    \caption{}
  \end{subfigure}
  \hfill
  \begin{subfigure}[b]{0.32\textwidth}
    \includegraphics[width=\textwidth]{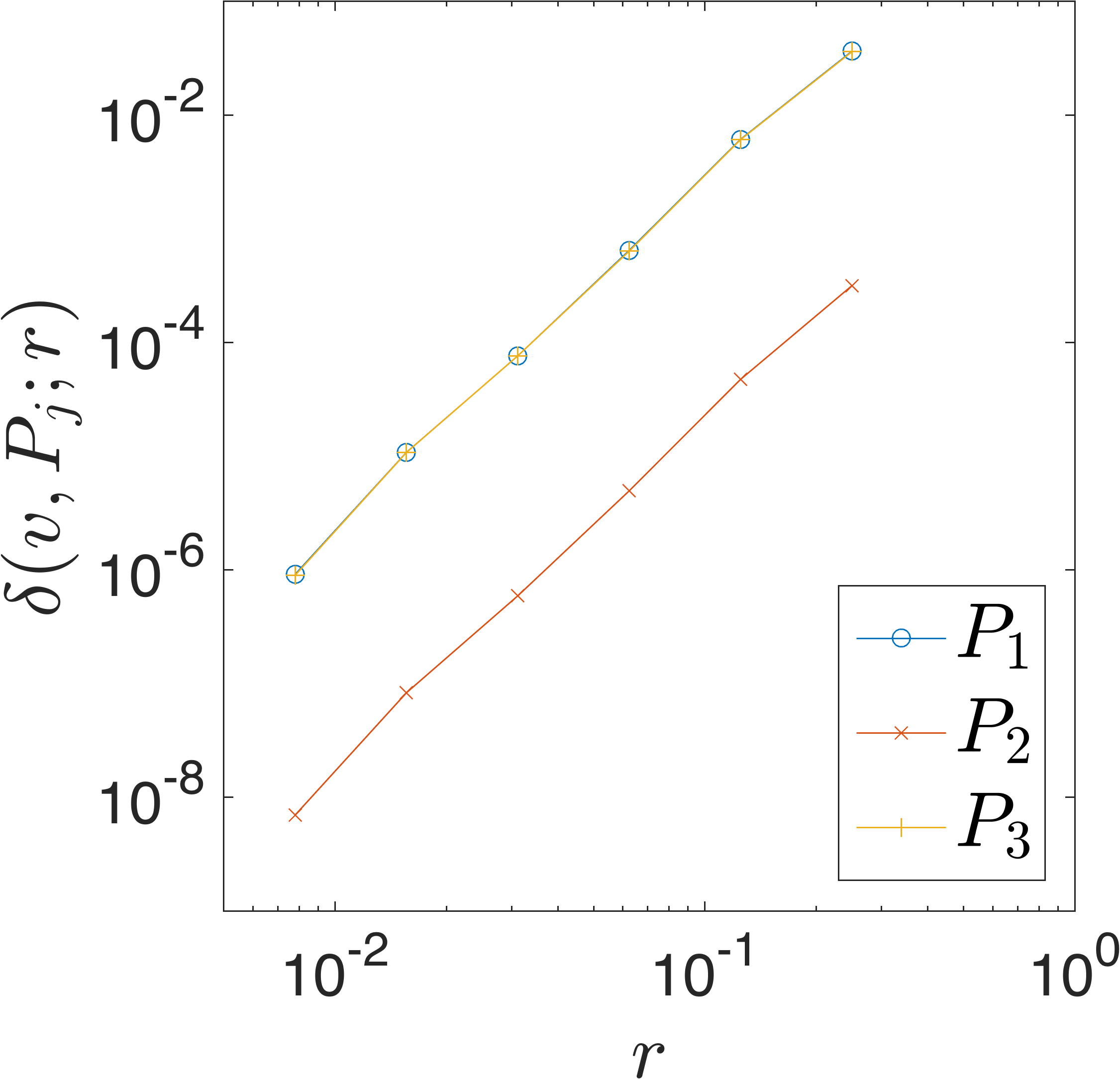}
    \caption{}
  \end{subfigure}
  \hfill
  \begin{subfigure}[b]{0.32\textwidth}
    \includegraphics[width=\textwidth]{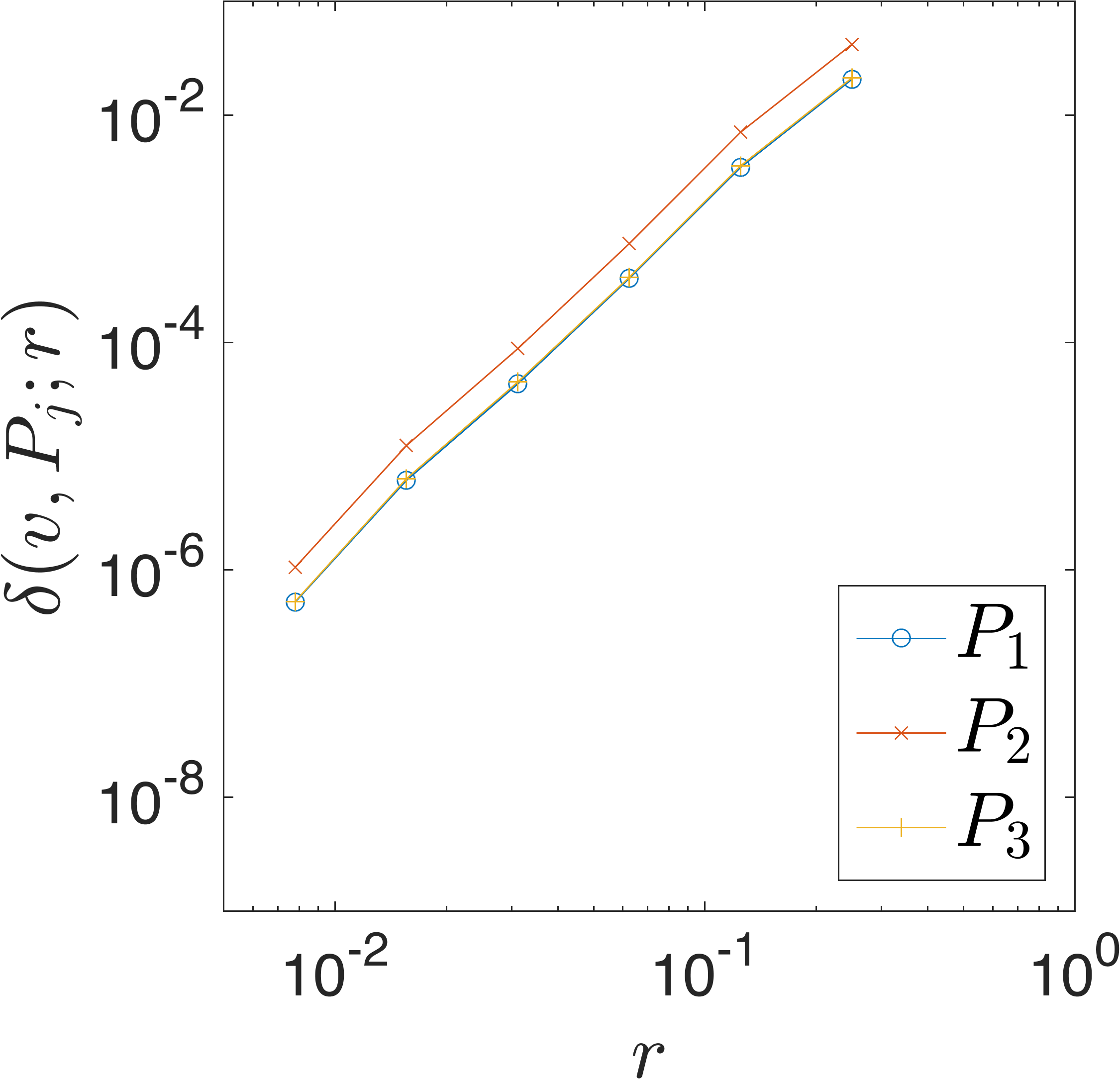}
    \caption{}
  \end{subfigure}
  \caption{The vanishing of the transmission eigenfunctions towards the corner points of the equilateral triangle with refractive index $n=16$. (A) $v=u_0^{(1)}$; (B) $v=u_0^{(2)}$; (C) $v=u_0^{(3)}$; (D) $v=u^{(1)}$; (E) $v=u^{(2)}$; (F) $v=u^{(3)}$.}
  \label{fig:equilateral_triangle_delta}
\end{figure}

We refer corners as the singular points. Let $P$ be a corner point of $D$ and $B_r(P)$ and $D_r(P)$ be defined as before. Define the average $L^2$ norm of a transmission eigenfunction $v$ in $D_r(P)$ as
\begin{align*}
  \delta(v,P;r) = \frac{1}{\sqrt{|D_r(P)|}} \left\|v\right\|_{L^2 \left(D_r(P)\right)}.
\end{align*}
We numerically demonstrate that $\delta(v,P;r) \to 0$ as $r \to 0$ and estimate the rate of convergence. In Figure
\ref{fig:equilateral_triangle_delta} we plot $\delta(u_0^{(i)},
P_j;r)$ and $\delta(u^{(i)}, P_j; r)$ versus $r$ for $i=1,2,3$, $j=1,2,3$.
From Figure \ref{fig:equilateral_triangle_delta} we further confirm that each transmission eigenfunction vanishes towards every corner point of $D$ in the sense that $\delta(v,P;r) \to 0$ as $r \to 0$. Furthermore, define the convergence rate $\alpha$ through the asymptotic behavior of $\delta(v,P;r)$
  \begin{align}
    \label{eq:9}
    \delta(v,P;r) \sim c(v,P) r^{-\alpha(v,P)} \quad {\rm as}~r \to 0
  \end{align}
  for each $v=u_0^{(i)}$ or $v=u^{(i)}$ and each $P=P_j$, where $c>0$ and $\alpha>0$ are constants independent of $r$ but may dependent on $v$ and $P$. The constants $c$ and $\alpha$ may also depend on the domain $D$ and the refractive index $n$.

We are more concerned with the order of convergence $\alpha$. Fitting the data in Figure \ref{fig:equilateral_triangle_delta} by linear polynomials we obtain the estimates of the convergence order shown in Table \ref{tab:equilateral_triangle_order}. From this result we conjecture that the order of convergence $\alpha(v,P)$ is independent of the transmission eigenfunction $v$.
\begin{table}[t]
  \centering
  \begin{subtable}[c]{0.4\textwidth}
    \centering
    \begin{tabular}{cccc}
      \toprule
      & $P_1$ & $P_2$ & $P_3$ \\
      \midrule
      $u_0^{(1)}$ & 3.03 & 3.04 & 3.03 \\[5pt]
      $u_0^{(2)}$ & 3.05 & 3.08 & 3.05 \\[5pt]
      $u_0^{(3)}$ & 3.05 & 3.05 & 3.05 \\
      \bottomrule
    \end{tabular}
    \caption{}
  \end{subtable}
  \begin{subtable}[c]{0.4\textwidth}
    \centering
    \begin{tabular}{cccc}
      \toprule
      & $P_1$ & $P_2$ & $P_3$ \\
      \midrule
      $u^{(1)}$ & 3.03 & 3.03 & 3.02 \\[5pt]
      $u^{(2)}$ & 3.05 & 3.08 & 3.06 \\[5pt]
      $u^{(3)}$ & 3.05 & 3.05 & 3.06 \\
      \bottomrule
    \end{tabular}
    \caption{}
  \end{subtable}
  \caption{Order of convergence $\alpha(v,P)$ of the transmission eigenfunctions for the equilateral triangle with refractive index $n=16$; (A) $v=u_0^{(j)},j=1,2,3$; (B) $v=u^{(j)},j=1,2,3$.}
  \label{tab:equilateral_triangle_order}
\end{table}

\subsection{Example: Right triangle}
\begin{figure}[t]
  \centering
  \begin{subfigure}[b]{0.32\textwidth}
    \includegraphics[width=\textwidth]{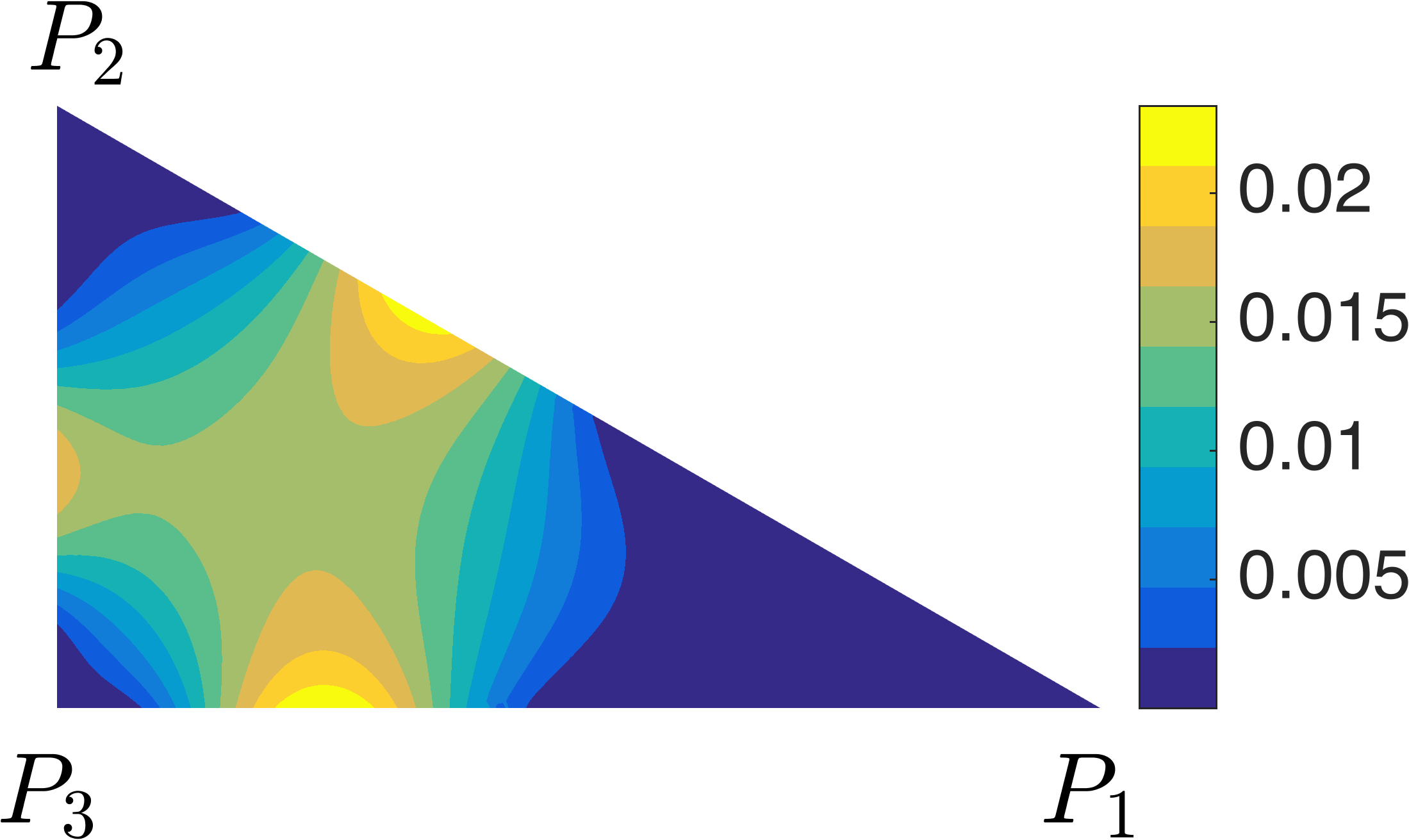}
    \caption{}
  \end{subfigure}
  \hfill
  \begin{subfigure}[b]{0.32\textwidth}
    \includegraphics[width=\textwidth]{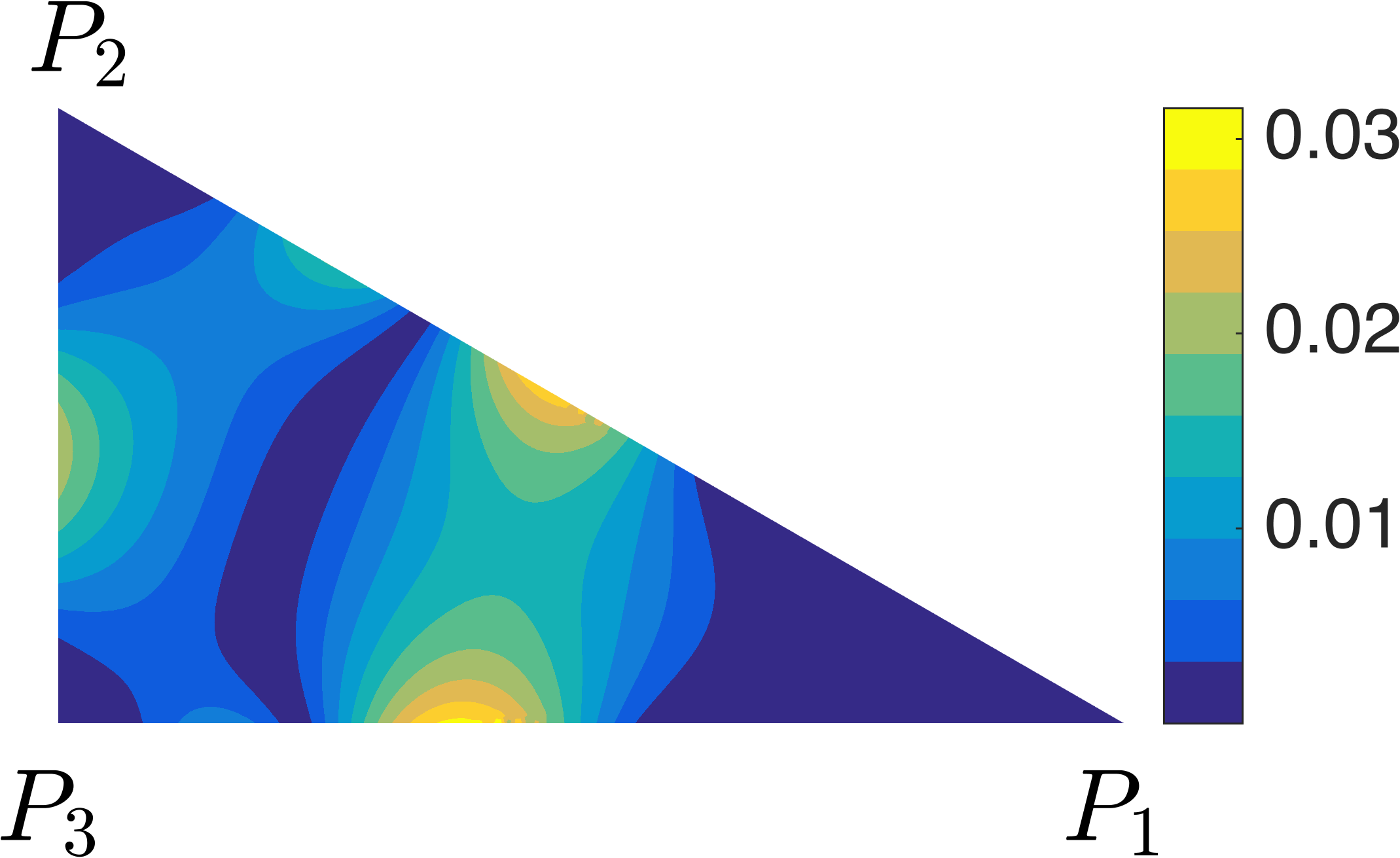}
    \caption{}
  \end{subfigure}
  \hfill
  \begin{subfigure}[b]{0.32\textwidth}
    \includegraphics[width=\textwidth]{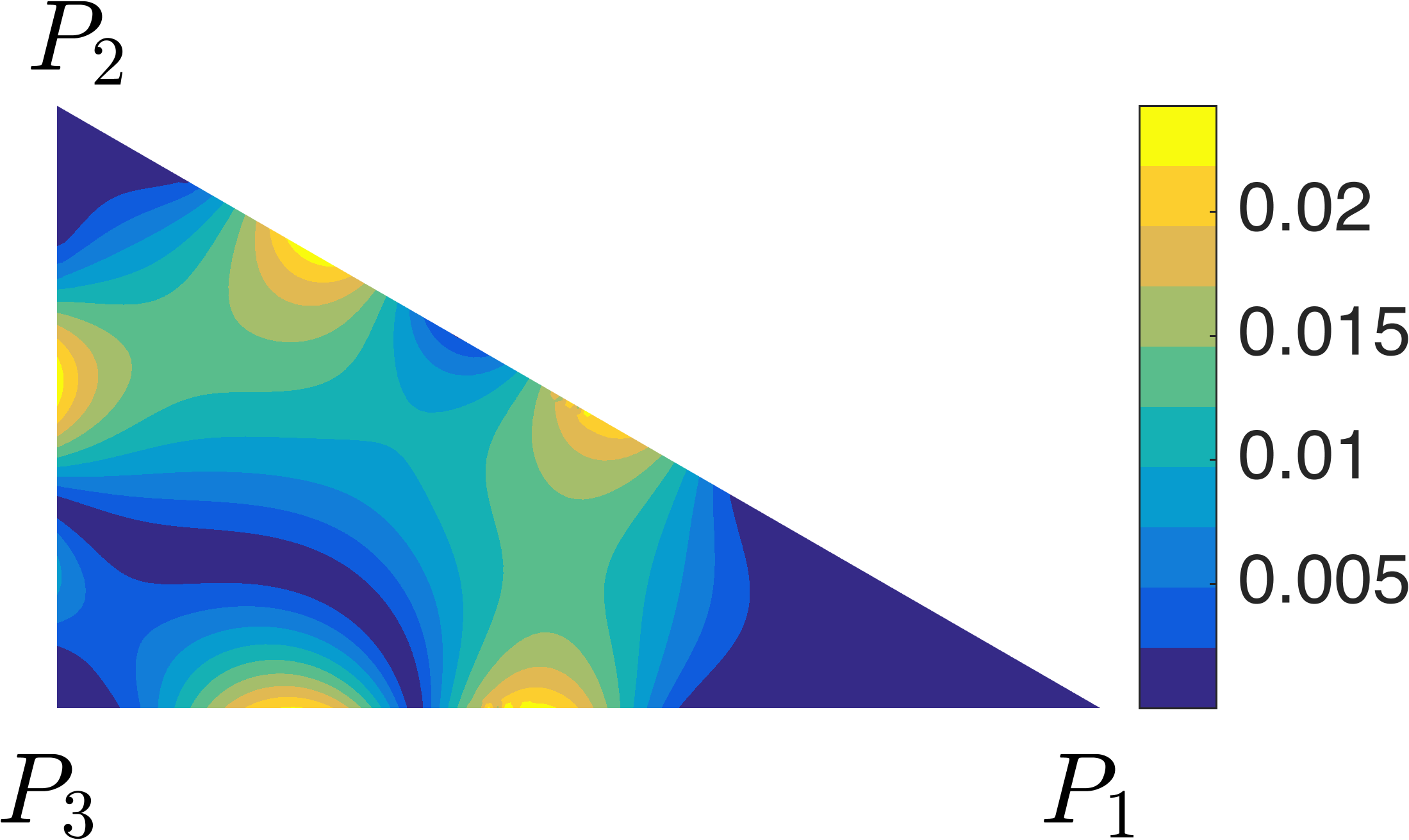}
    \caption{}
  \end{subfigure}
  \caption{The magnitude of transmission eigenfunctions for the right triangle with refractive index $n=16+8\sin(4xy)$. (A) $|u_0^{(1)}|$; (B) $|u_0^{(2)}|$; (C) $|u_0^{(3)}|$.}
  \label{fig:right_triangle_u0}
\end{figure}

\begin{figure}[t]
  \centering
  \begin{subfigure}[b]{0.32\textwidth}
    \includegraphics[width=\textwidth]{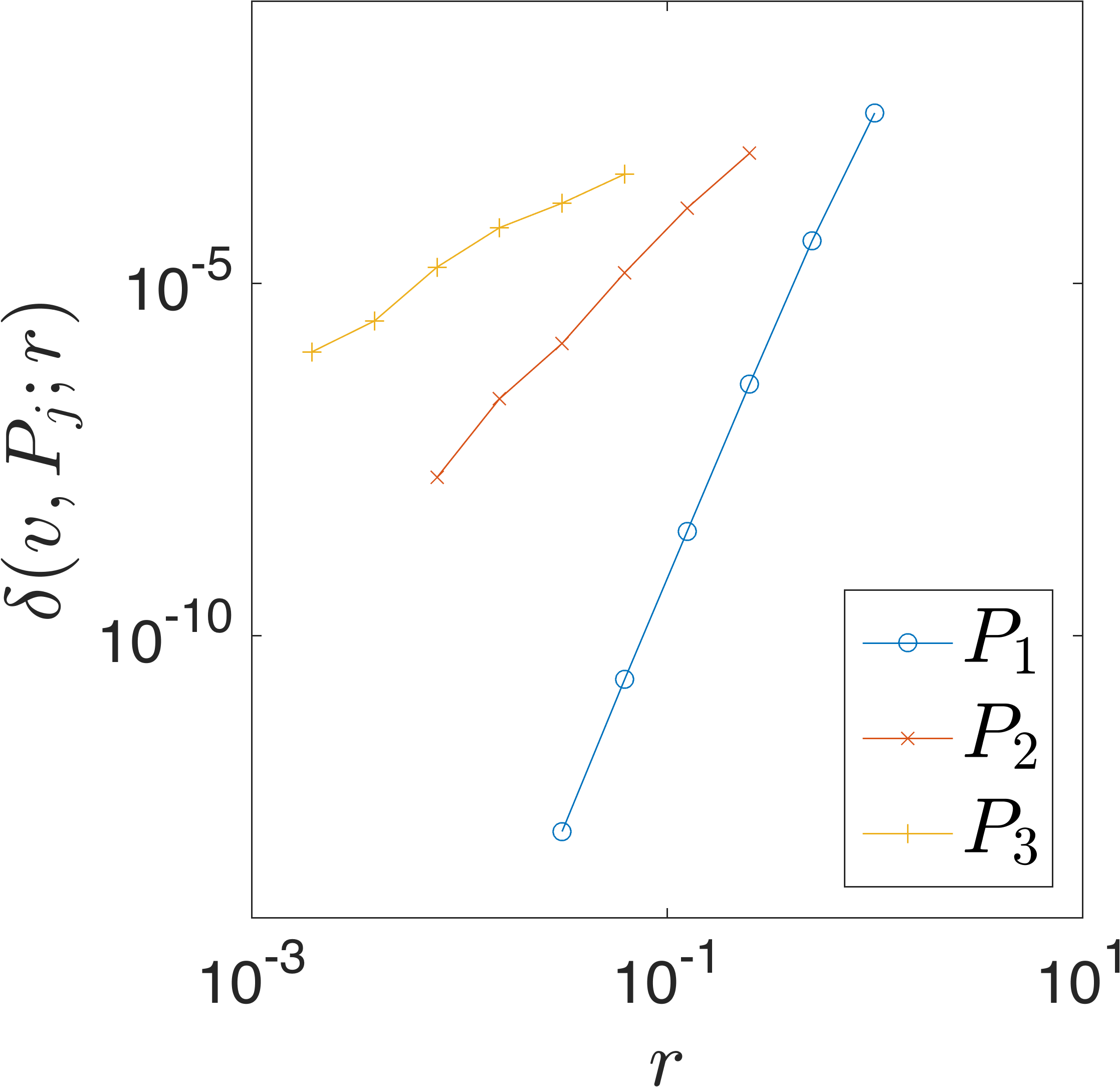}
    \caption{}
  \end{subfigure}
  \hfill
  \begin{subfigure}[b]{0.32\textwidth}
    \includegraphics[width=\textwidth]{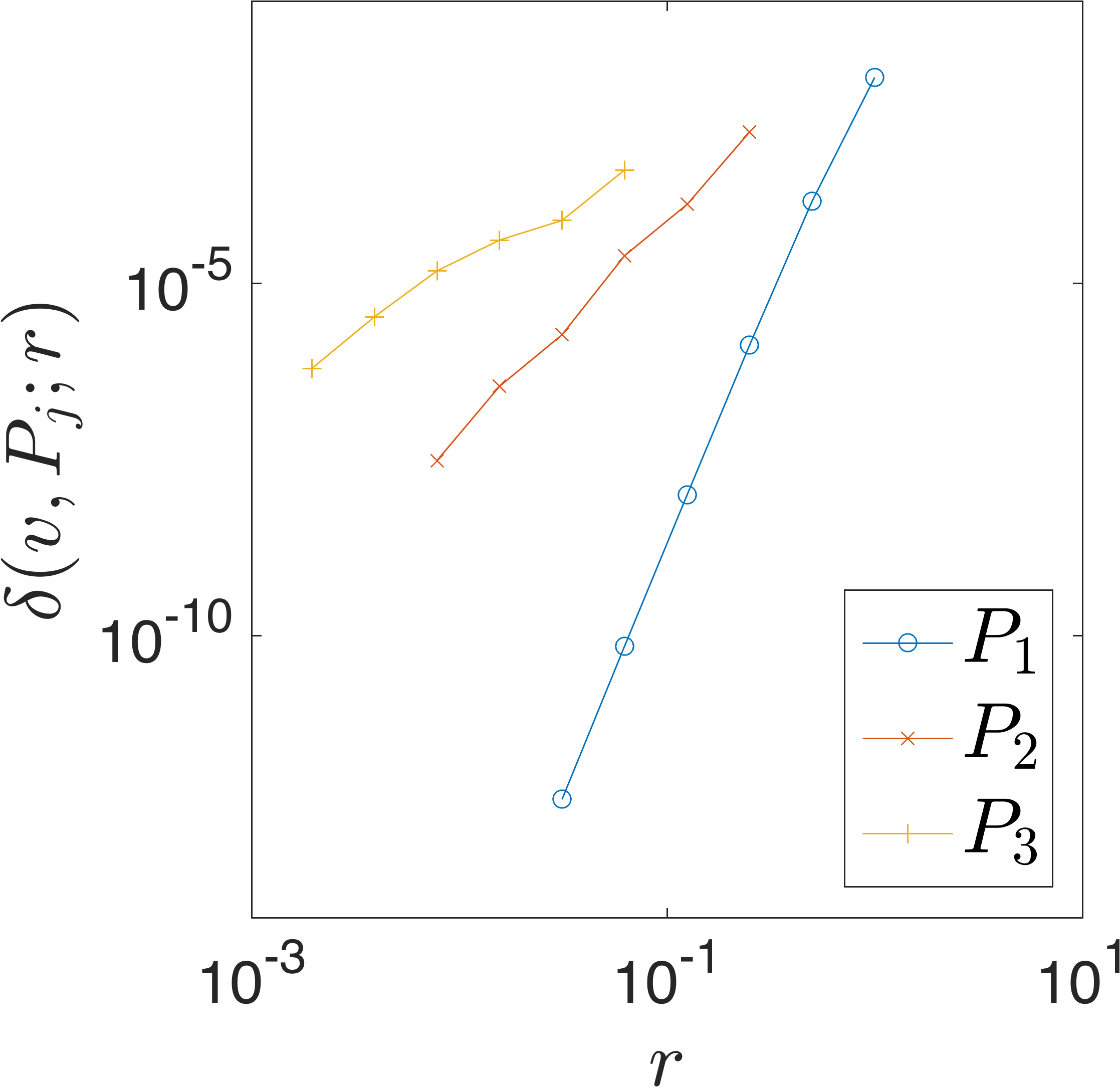}
    \caption{}
  \end{subfigure}
  \hfill
  \begin{subfigure}[b]{0.32\textwidth}
    \includegraphics[width=\textwidth]{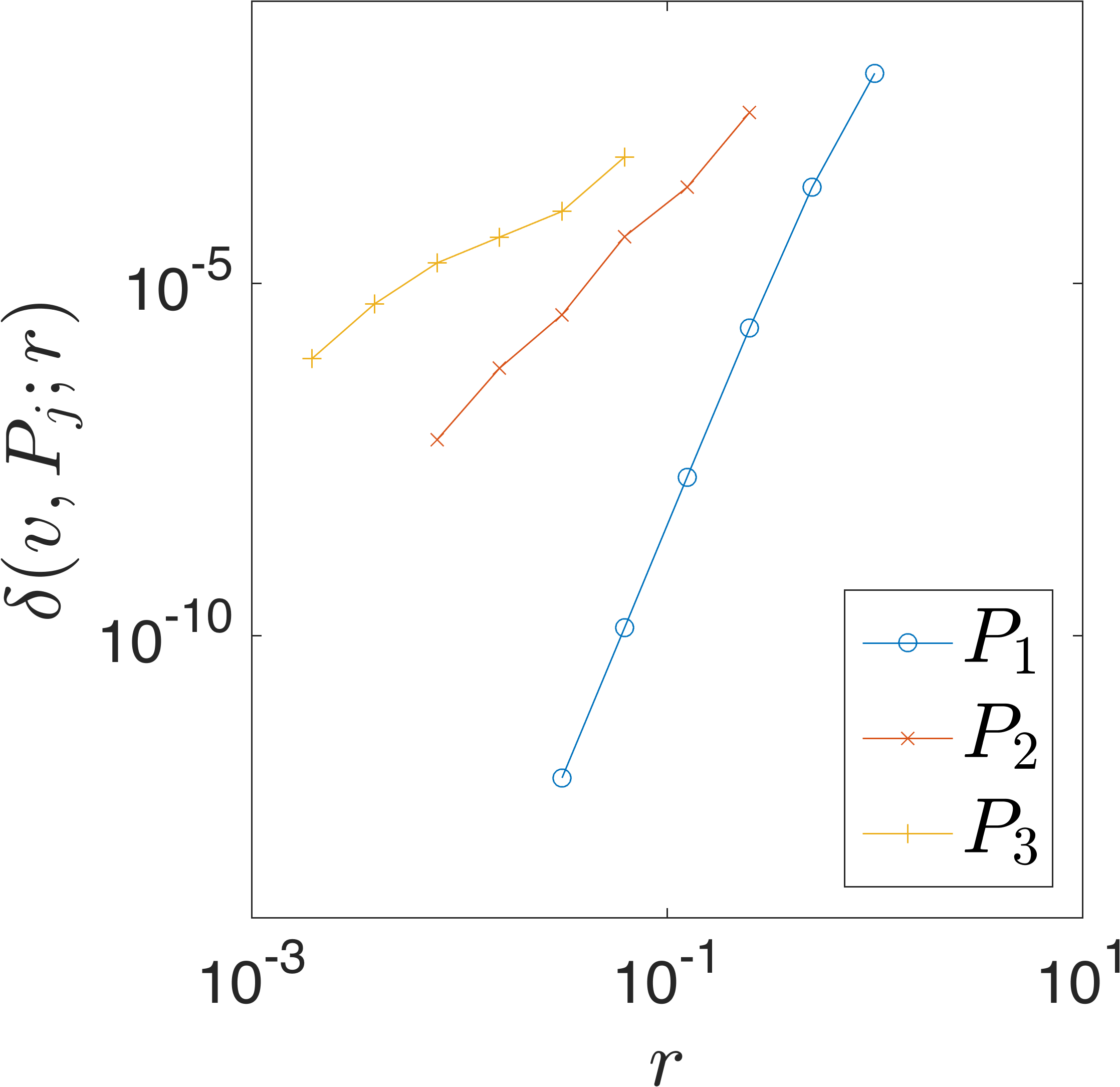}
    \caption{}
  \end{subfigure}
  \caption{The vanishing of the transmission eigenfunctions towards the corner points of the right triangle with refractive index $n=16+8\sin(4xy)$. (A) $v=u_0^{(1)}$; (B) $v=u_0^{(2)}$; (C) $v=u_0^{(3)}$. }
  \label{fig:right_triangle_delta}
\end{figure}

In this example, we investigate the dependence of convergence rate $\alpha$ on the domain $D$, the corner point $P_j$ and the refractive index $n$. Let $D$ be the right triangle with vertices at $P_1=(\sqrt{3},0), P_2=(0,1)$ and $P_3=(0,0)$. Let the refractive index be $n=16+8\sin(4xy)$. The first three real transmission eigenvalues are computed to be $2.1651, 2.4993$ and $2.8814$ with corresponding eigenfunctions shown in Figure \ref{fig:right_triangle_u0}. The convergence of the average $L^2$ norm of the eigenfunctions towards the corner points is shown in Figure \ref{fig:right_triangle_delta} and the order of convergence is listed in Table \ref{tab:right_triangle_order_u0}.
\begin{table}[t]
  \centering
  \begin{tabular}{cccc}
    \toprule
   & $P_1$ & $P_2$ & $P_3$ \\
    \midrule
    $u_0^{(1)}$ & 6.81 & 3.05 & 1.72 \\[5pt]
    $u_0^{(2)}$ & 6.86 & 3.05 & 1.77 \\[5pt]
    $u_0^{(3)}$ & 6.74 & 3.04 & 1.77 \\
    \bottomrule
  \end{tabular}
  \caption{Order of convergence $\alpha(u_0^{(i)},P_j)$ of the transmission eigenfunctions for the right triangle with refractive index $n=16+8\sin(4xy)$.}
  \label{tab:right_triangle_order_u0}
\end{table}

Comparing the results in Table \ref{tab:equilateral_triangle_order} and Table \ref{tab:right_triangle_order_u0} we conjecture that the order of convergence $\alpha(v,P)$ in \eqref{eq:9} is also independent of the domain $D$ and the refractive index $n$. Furthermore, $\alpha(v,P)$ increases as the sharpness of the corner $P$ increases.

In the previous examples all the corners of the domain have angles less than $\pi$. In the next example we consider corners with angle greater than $\pi$.

\subsection{Example: Arrow shaped polygon}
\begin{figure}[t]
  \centering
  \begin{subfigure}[b]{0.32\textwidth}
    \includegraphics[width=\textwidth]{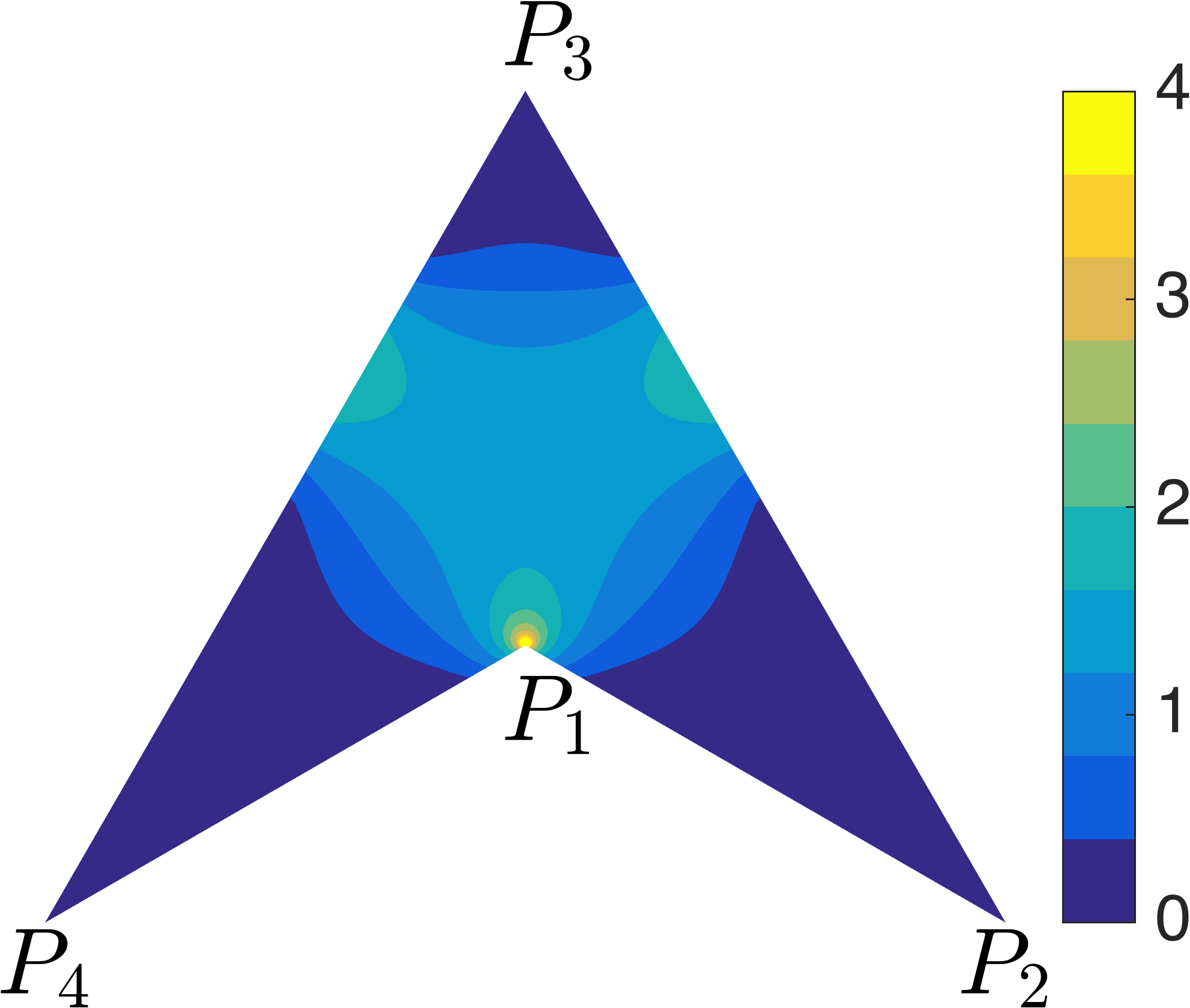}
    \caption{}
  \end{subfigure}
  \hfill
  \begin{subfigure}[b]{0.32\textwidth}
    \includegraphics[width=\textwidth]{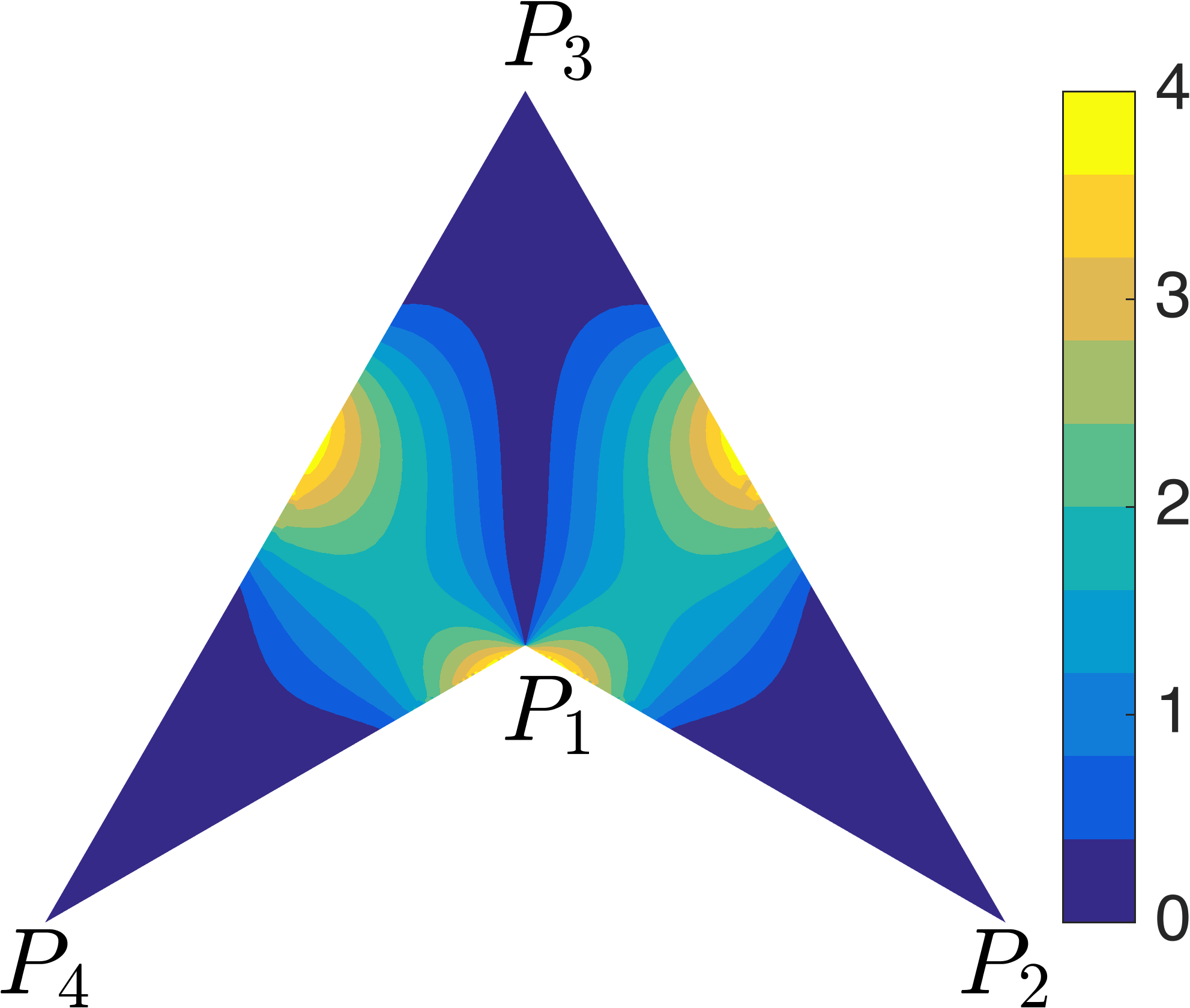}
    \caption{}
  \end{subfigure}
  \hfill
  \begin{subfigure}[b]{0.32\textwidth}
    \includegraphics[width=\textwidth]{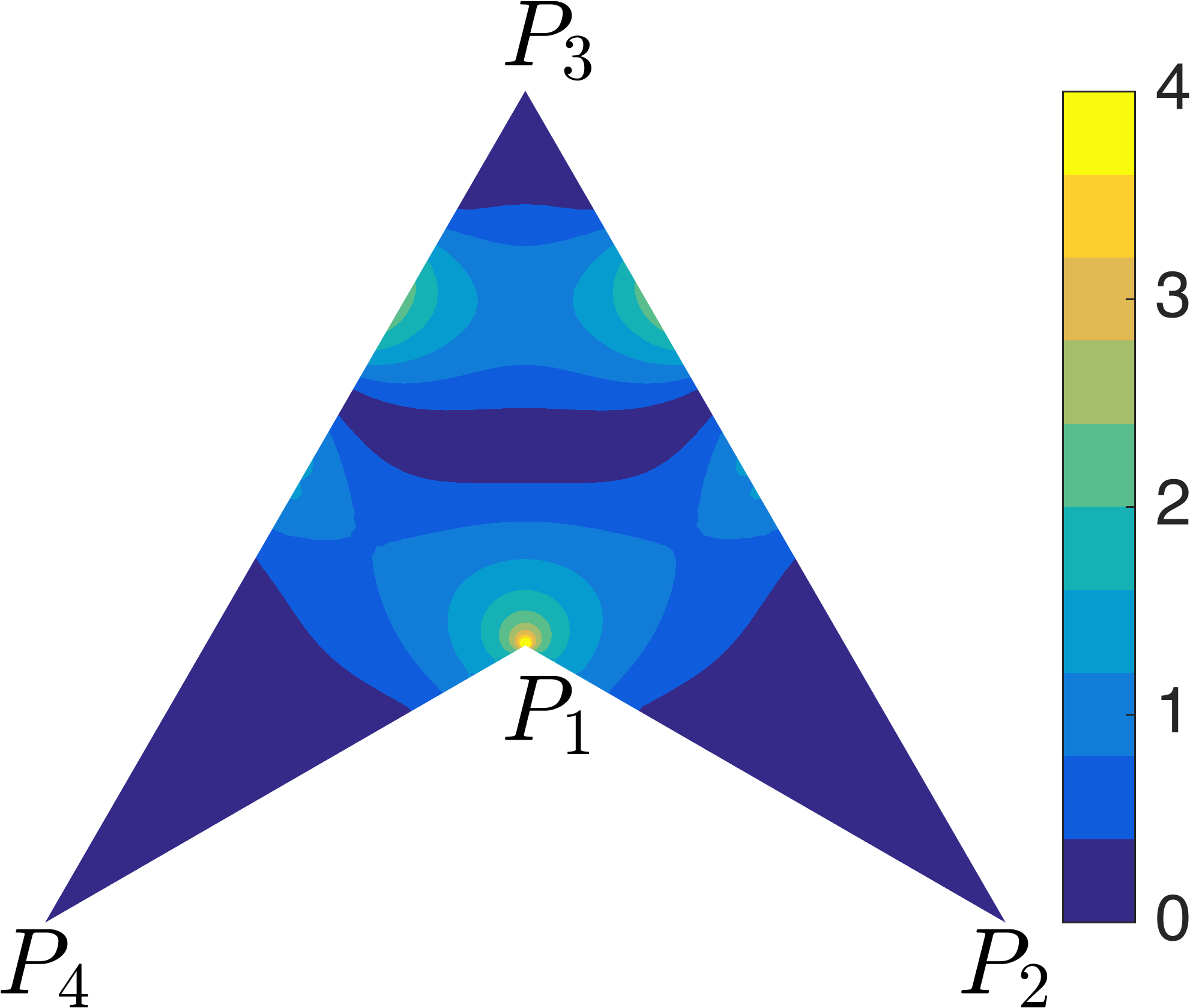}
    \caption{}
  \end{subfigure}
  \hfill
  \centering
  \begin{subfigure}[b]{0.32\textwidth}
    \includegraphics[width=\textwidth]{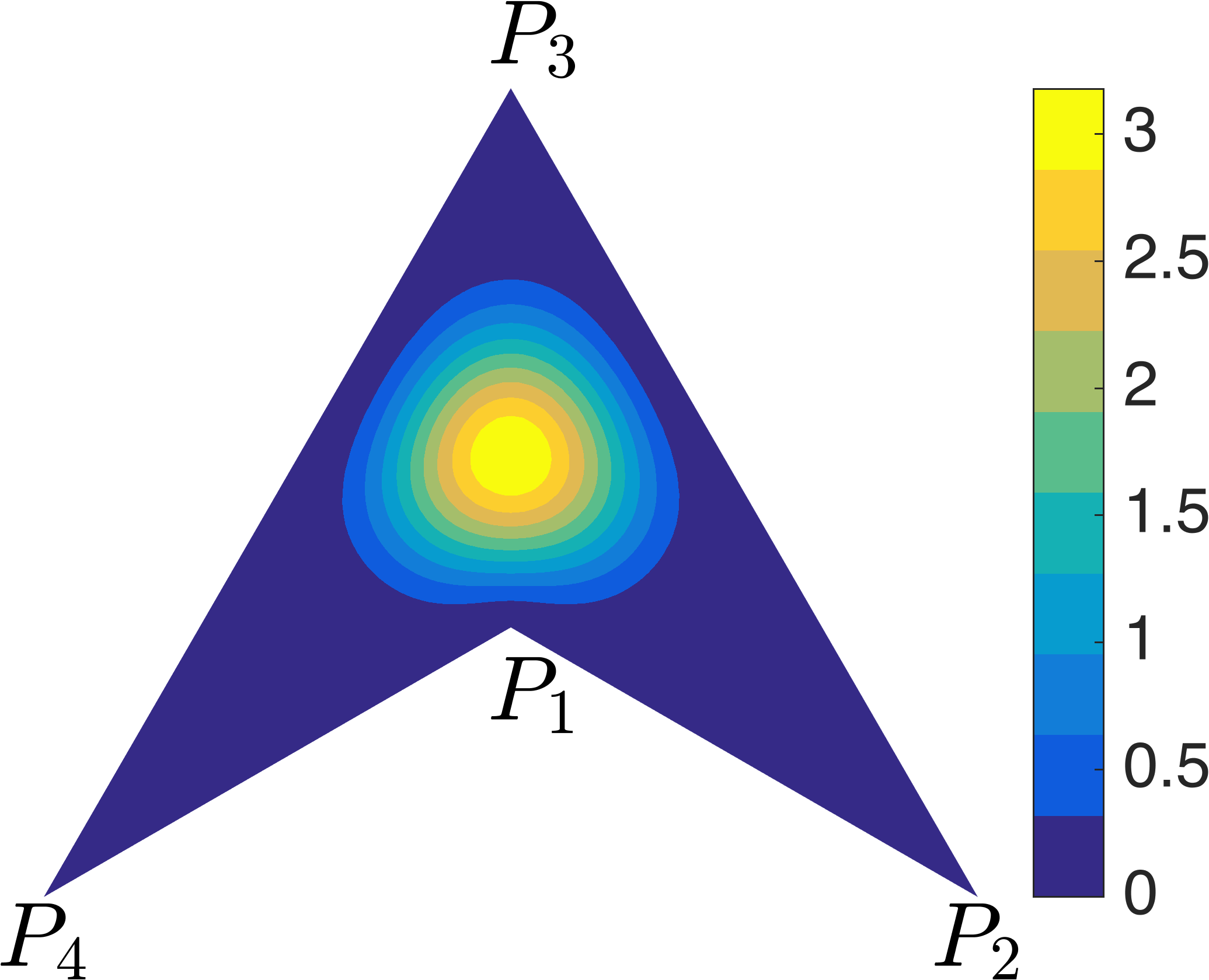}
    \caption{}
  \end{subfigure}
  \hfill
  \begin{subfigure}[b]{0.32\textwidth}
    \includegraphics[width=\textwidth]{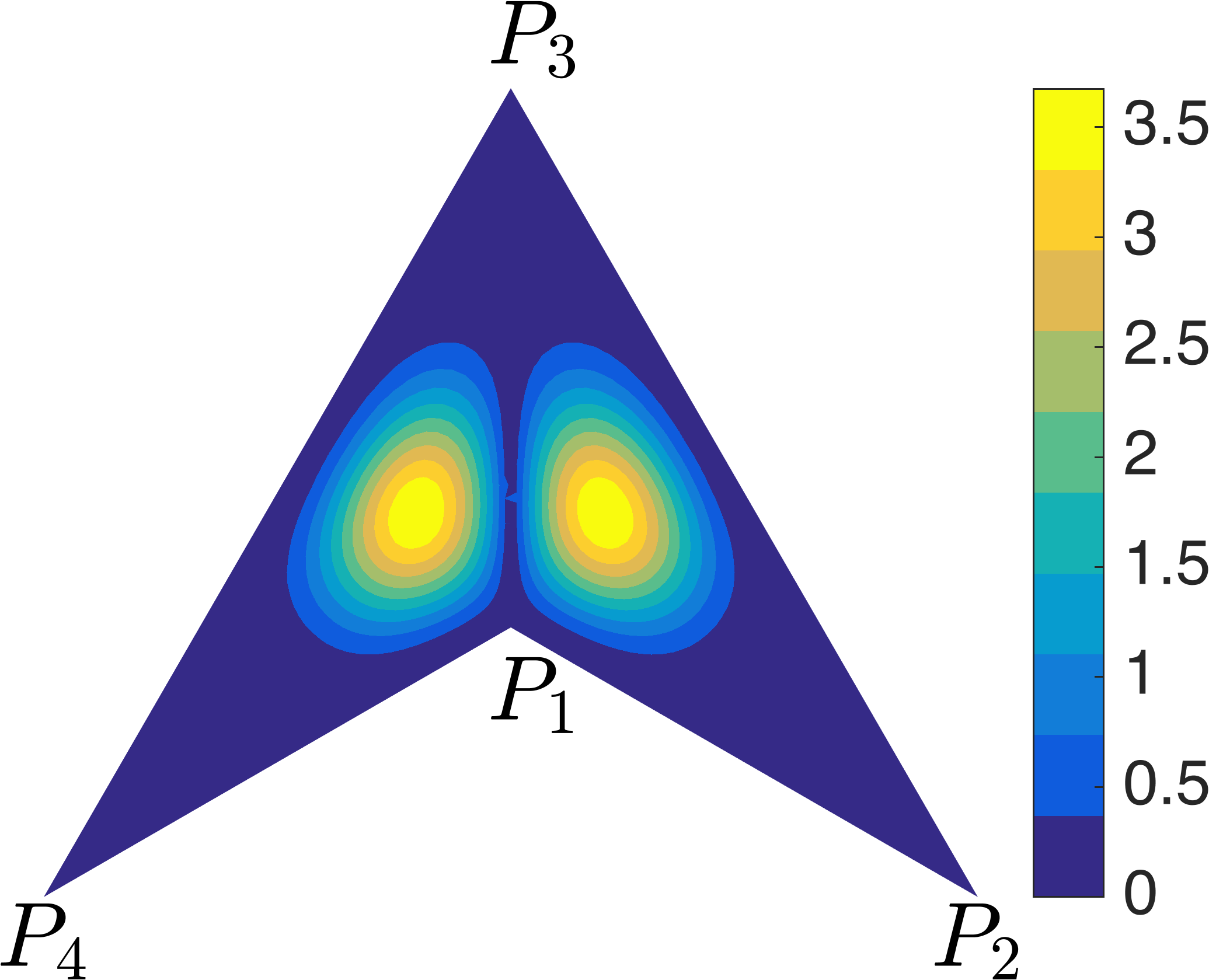}
    \caption{}
  \end{subfigure}
  \hfill
  \begin{subfigure}[b]{0.32\textwidth}
    \includegraphics[width=\textwidth]{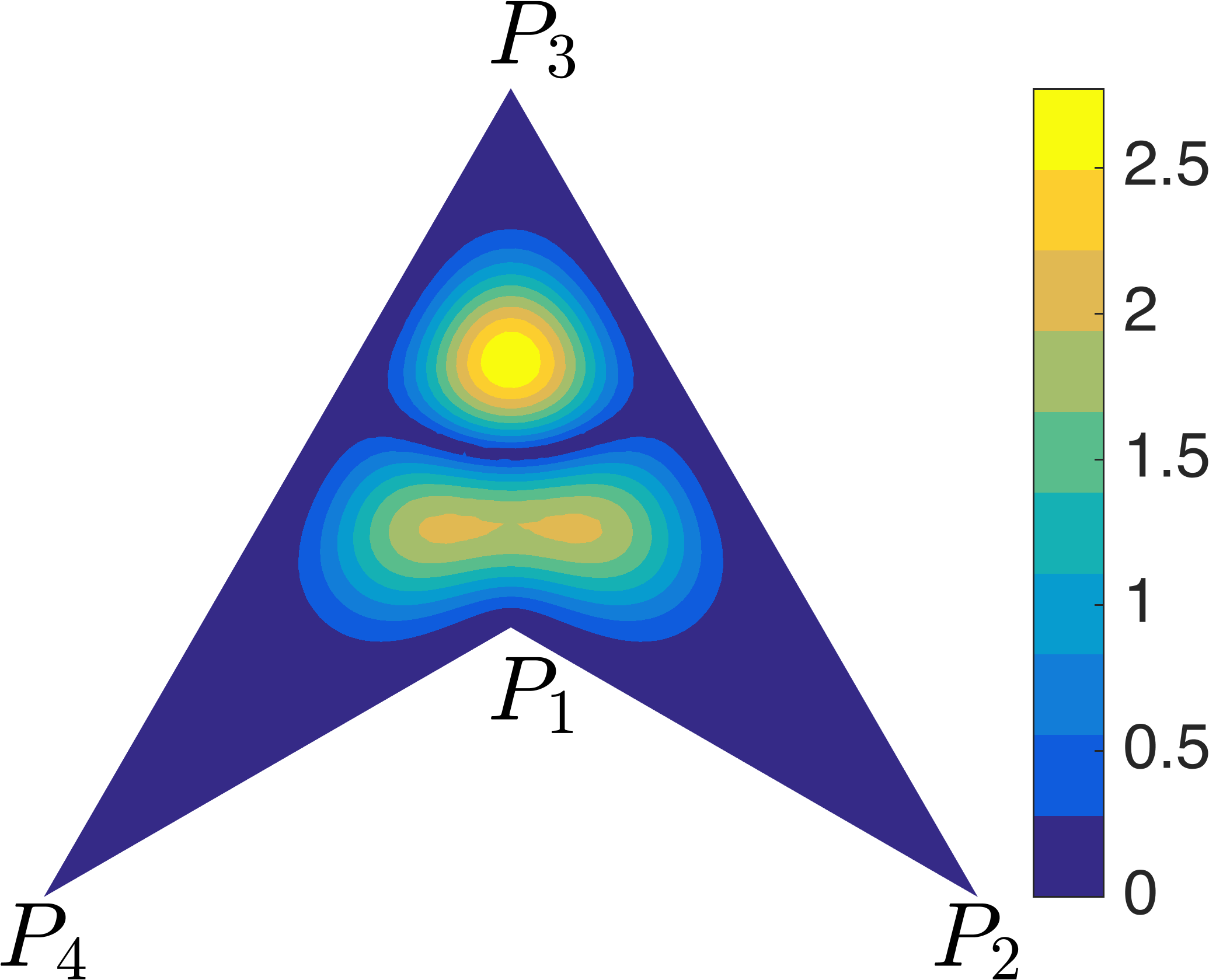}
    \caption{}
  \end{subfigure}
  \caption{The magnitude of transmission eigenfunctions for the arrow shaped domain with refractive index $n=16$. (A) $|u_0^{(1)}|$; (B) $|u_0^{(2)}|$; (C) $|u_0^{(3)}|$; (D) $|u^{(1)}-u_0^{(1)}|$; (E) $|u^{(2)}-u_0^{(2)}|$; (F) $|u^{(3)}-u_0^{(3)}|$.}
  \label{fig:arrow}
\end{figure}

\begin{figure}[t]
  \centering
  \begin{subfigure}[b]{0.32\textwidth}
    \includegraphics[width=\textwidth]{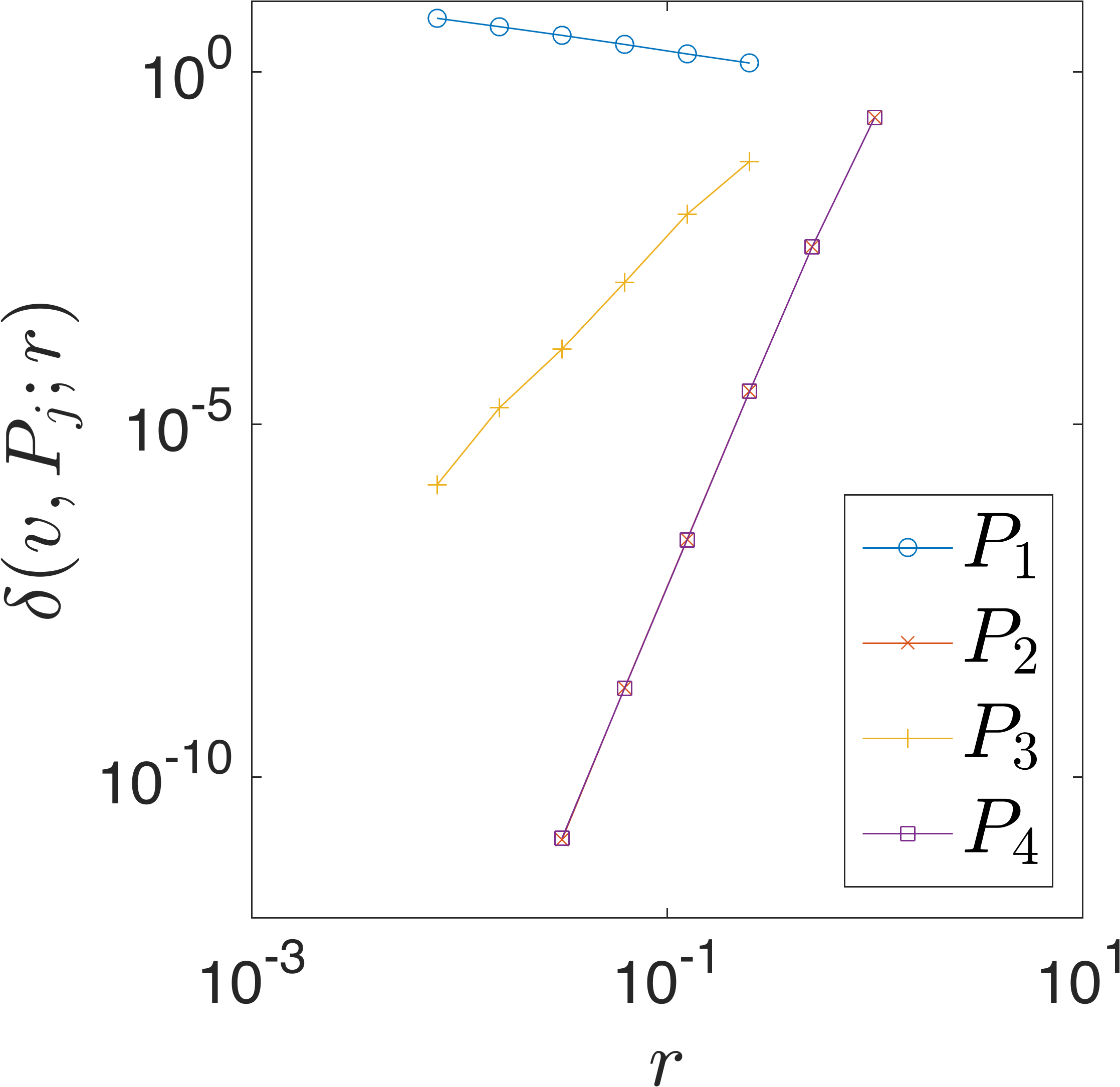}
    \caption{}
  \end{subfigure}
  \hfill
  \begin{subfigure}[b]{0.32\textwidth}
    \includegraphics[width=\textwidth]{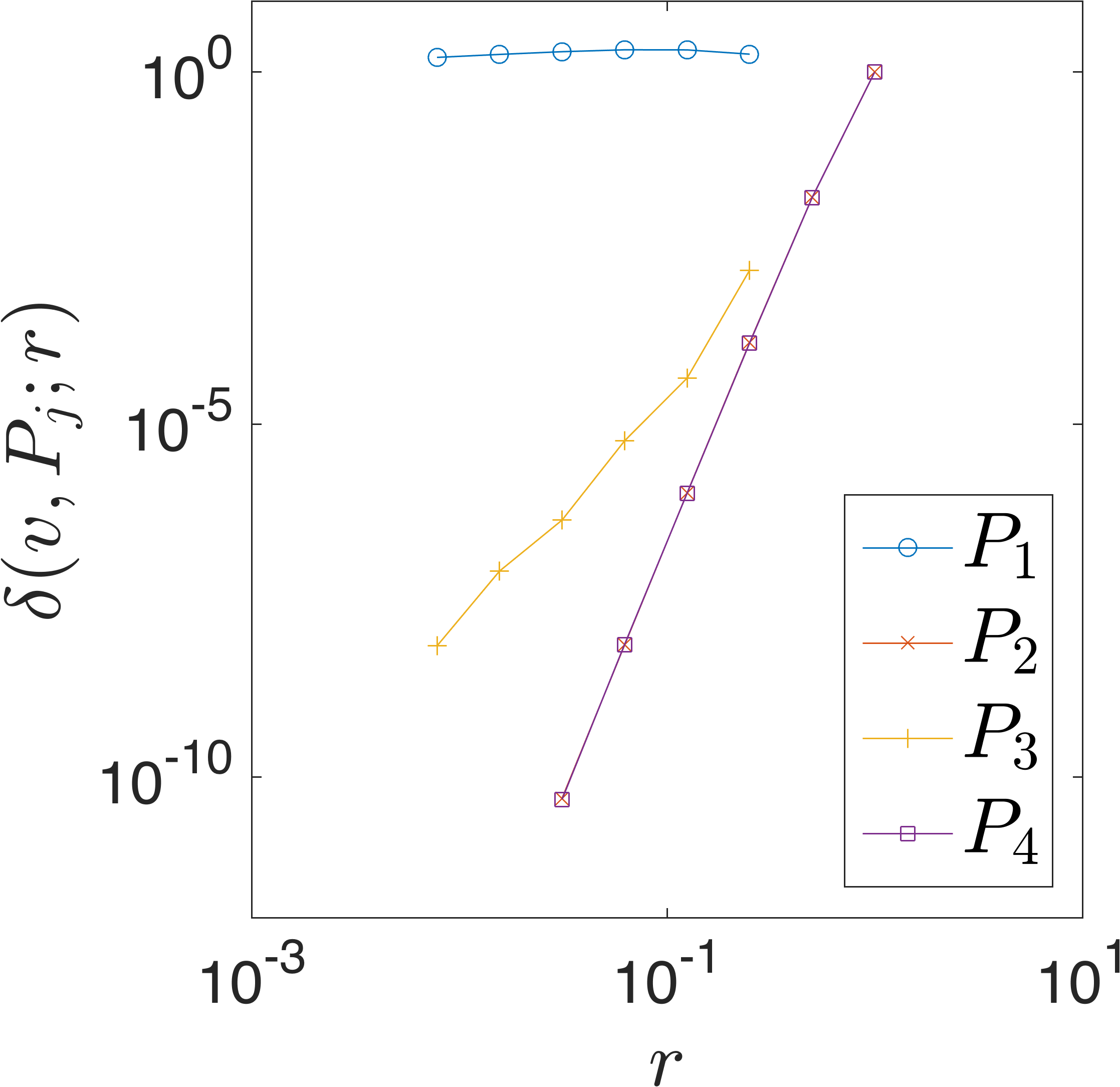}
    \caption{}
  \end{subfigure}
  \hfill
  \begin{subfigure}[b]{0.32\textwidth}
    \includegraphics[width=\textwidth]{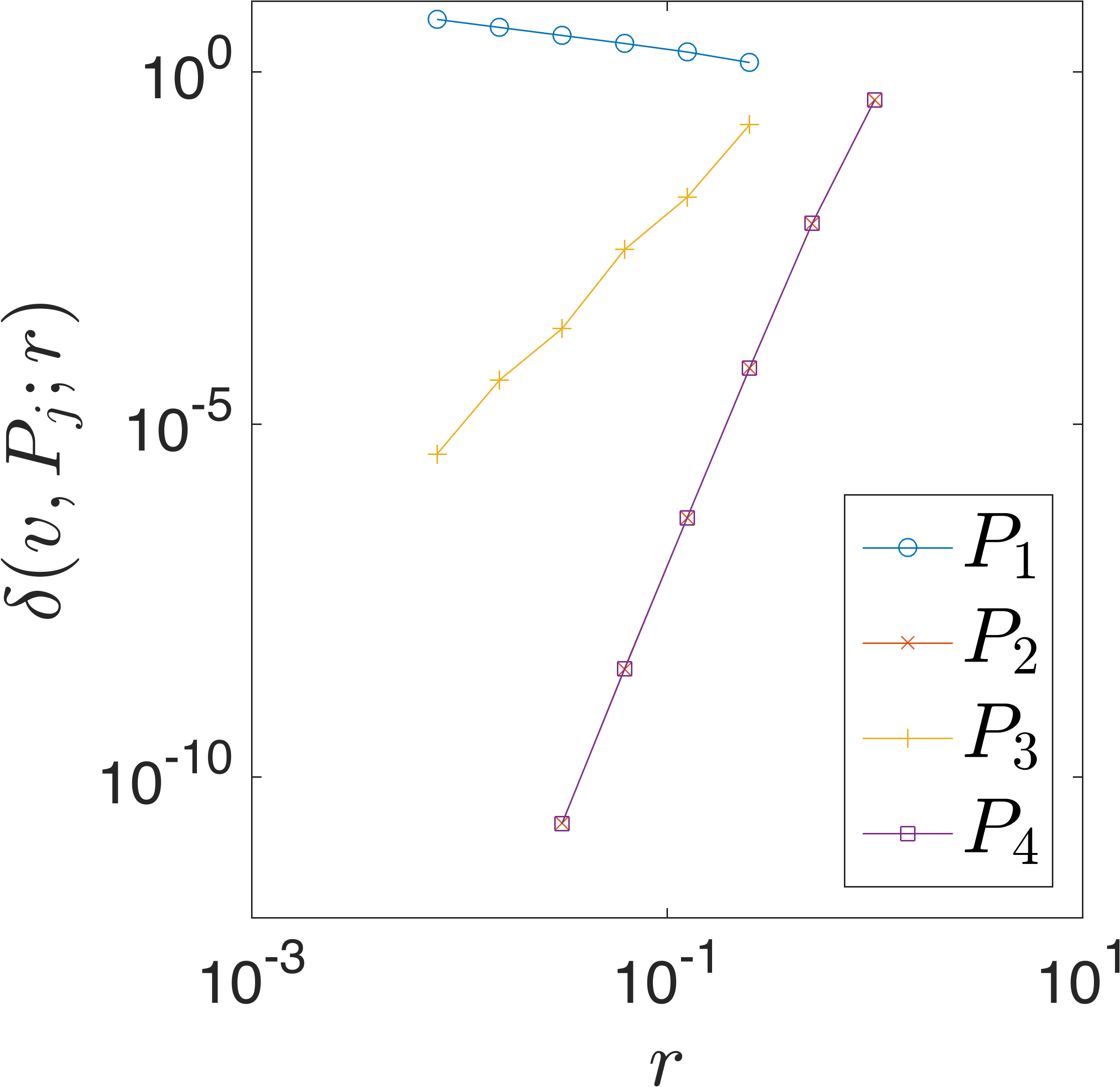}
    \caption{}
  \end{subfigure}
  \hfill
  \begin{subfigure}[b]{0.32\textwidth}
    \includegraphics[width=\textwidth]{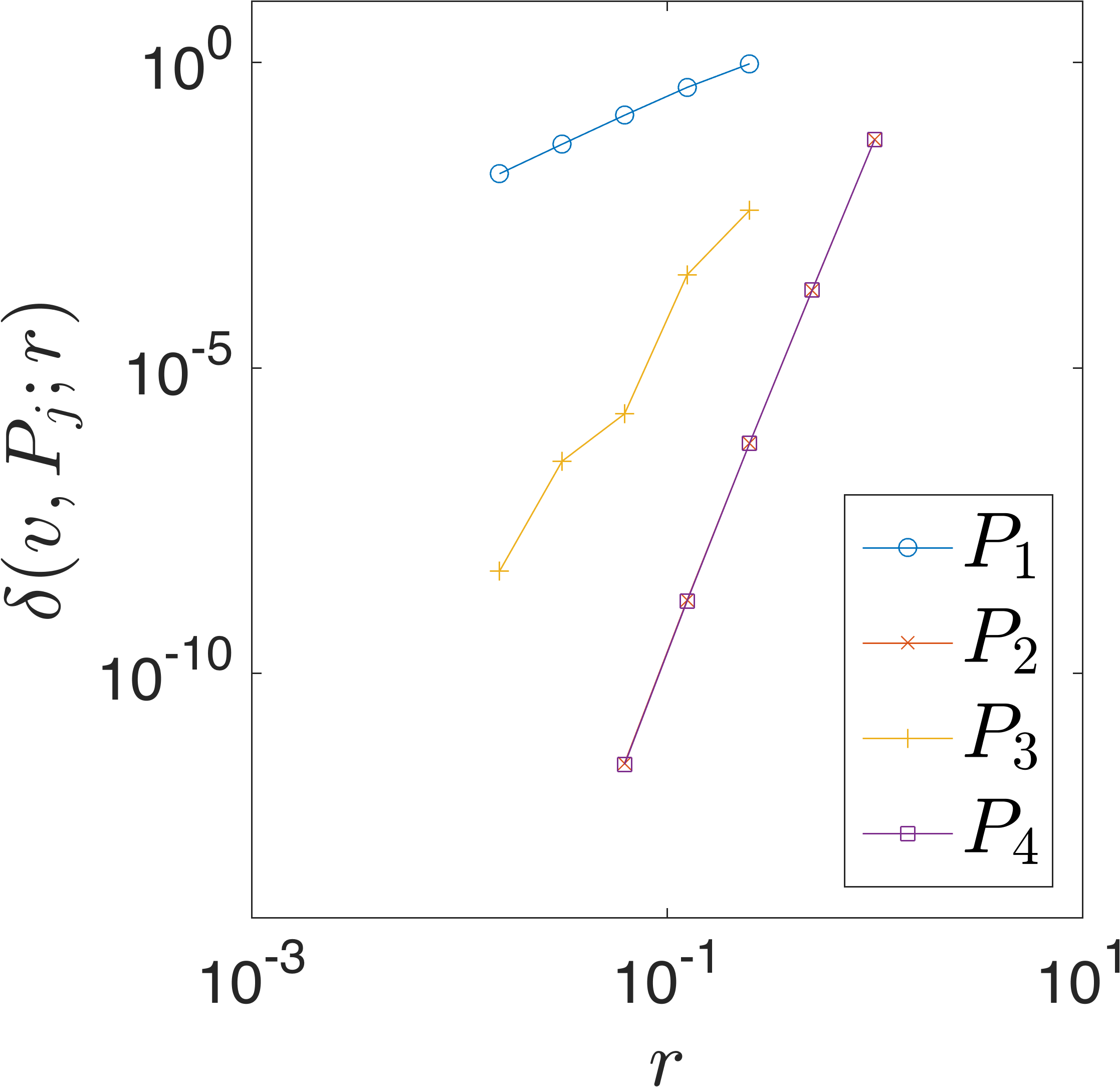}
    \caption{}
  \end{subfigure}
  \hfill
  \begin{subfigure}[b]{0.32\textwidth}
    \includegraphics[width=\textwidth]{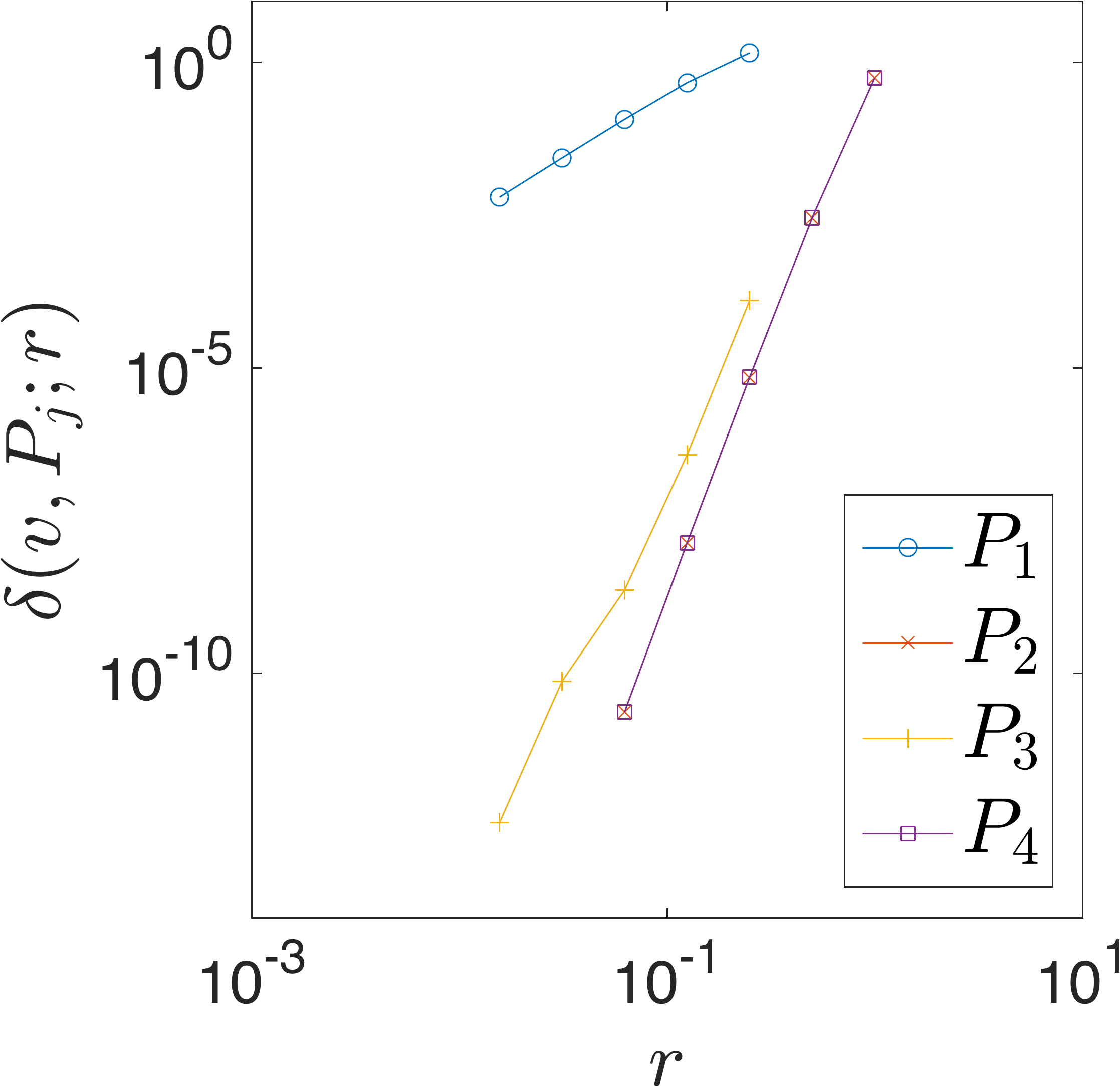}
    \caption{}
  \end{subfigure}
  \hfill
  \begin{subfigure}[b]{0.32\textwidth}
    \includegraphics[width=\textwidth]{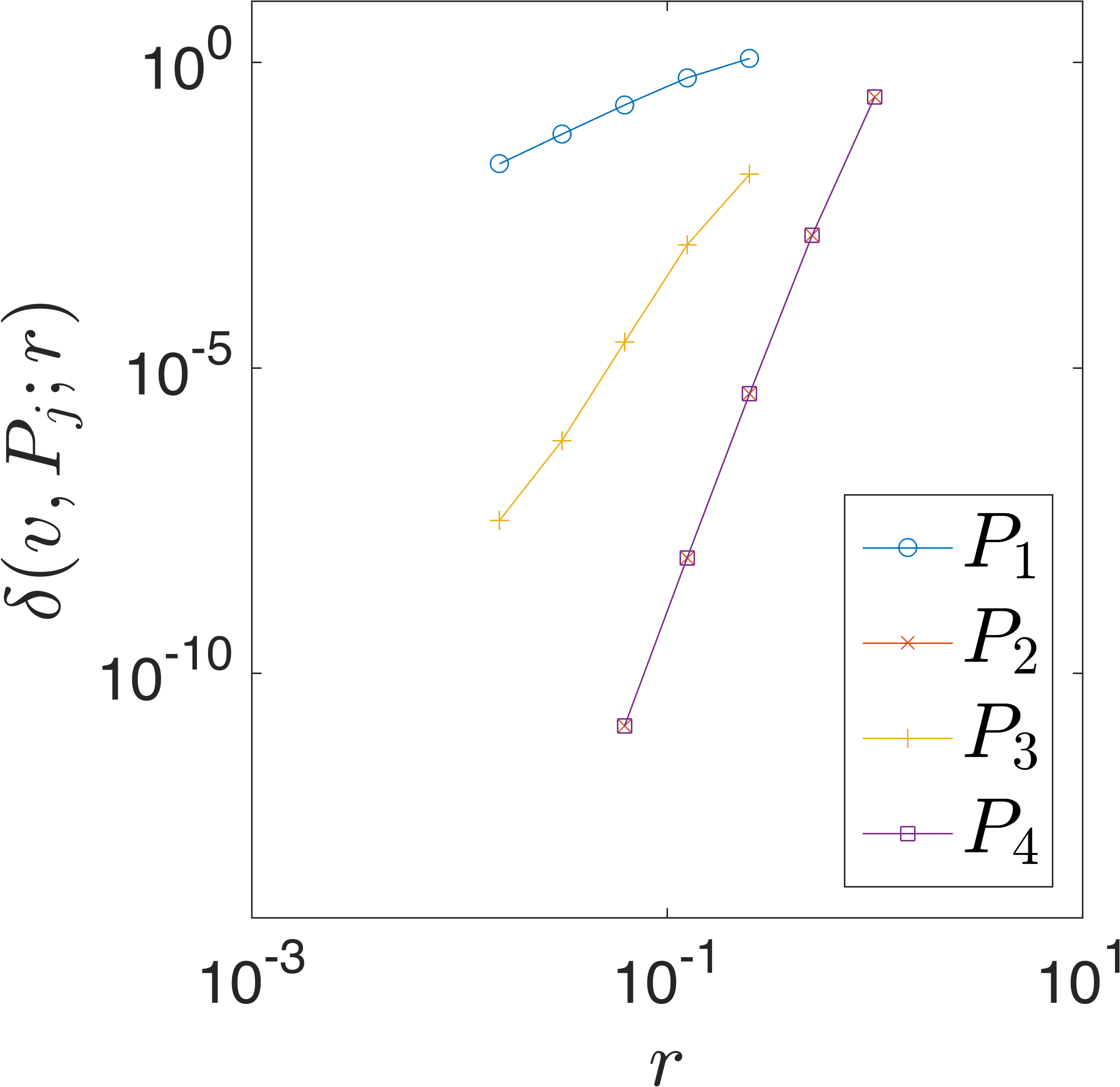}
    \caption{}
  \end{subfigure}
  \caption{The vanishing of the transmission eigenfunctions towards the corner points of the arrow shaped domain with refractive index $n=16$. (A) $u_0^{(1)}$; (B) $u_0^{(2)}$; (C) $u_0^{(3)}$; (D) $u^{(1)}-u_0^{(1)}$; (E) $u^{(2)}-u_0^{(2)}$; (F) $u^{(3)}-u_0^{(3)}$. }
  \label{fig:arrow_delta}
\end{figure}

Let $D$ be an arrow shaped polygon with vertices at $P_1=(0,\sqrt{3}/3), P_2=(1,0), P_3=(0, \sqrt{3})$ and $P_4=(-1,0)$. Let the refractive index be $n=16$. The first three real transmission eigenvalues are found to be $2.2954, 2.7635$ and $2.8532$ with corresponding eigenfunctions shown in Figure \ref{fig:arrow} (A)--(C). The convergence of the average $L^2$ norm of the eigenfunctions towards the corner points is shown in Figure \ref{fig:arrow_delta} (A)--(C) and the order of convergence for the corner points $P_2,P_3,P_4$ is listed in Table \ref{tab:arrow_order} (A).
\begin{table}[t]
  \centering
  \begin{subtable}[c]{0.4\textwidth}
    \centering
    \begin{tabular}{cccc}
      \toprule
      & $P_2$ & $P_3$ & $P_4$ \\
      \midrule
      $u_0^{(1)}$ & 6.84 & 3.05 & 6.83 \\[5pt]
      $u_0^{(2)}$ & 6.90 & 3.42 & 6.90 \\[5pt]
      $u_0^{(3)}$ & 6.87 & 3.06 & 6.87 \\
      \bottomrule
    \end{tabular}
    \caption{}
  \end{subtable}
  \begin{subtable}[c]{0.4\textwidth}
    \centering
    \begin{tabular}{ccccc}
      \toprule
      & $P_1$ & $P_2$ & $P_3$ & $P_4$ \\
      \midrule
      $v^{(1)}$ & 1.50 & 8.47 & 4.94 & 8.48 \\[5pt]
      $v^{(2)}$ & 1.98 & 8.66 & 6.91 & 8.66 \\[5pt]
      $v^{(3)}$ & 1.98 & 8.66 & 6.91 & 8.66 \\
      \bottomrule
    \end{tabular}
    \caption{}
  \end{subtable}
  \caption{Order of convergence of the transmission eigenfunctions for the arrow shaped domain with refractive index $n=16$.}
  \label{tab:arrow_order}
\end{table}
From Figure \ref{fig:arrow} and \ref{fig:arrow_delta} we conjecture
that the vanishing properties of the eigenfunctions $u_0^{(i)}$ and
$u^{(i)}$ are only valid if the angle of the corner is less than
$\pi$. However, if consider the difference
$v^{(j)}=u^{(i)}-u_0^{(i)}$, then $\delta(v^{(j)},P,r) \to 0$ as $r
\to 0$ even for corner points with angle greater than $\pi$. But the order of convergence depends on the eigenfunction. This can be seen from the Figure \ref{fig:arrow} (D)--(F), Figure \ref{fig:arrow_delta} (D)--(F) and Table \ref{tab:arrow_order} (B).

\subsection{Example: Moon shaped domain}
We consider domain with curved angles. Denote $\mathcal{K}_\omega=\{(r,\theta): r>0, 0 < \theta < \omega\}$ a dihedral domain in $\RR^2$ with opening angle $\omega \in (0,2\pi)$ at the origin $O$.
\begin{defn}
  \cite{MNP00} Let $D$ be a bounded open set in $\RR^2$. A point $P \in \partial D$ is called a corner point if there exists a neighborhood $V$ of $P$, a diffeomorphism $\Psi$ of class $C^2$ and an opening angle $\omega = \omega(P) \in (0,\pi) \cup (\pi,2\pi)$ such that
  \begin{align*}
    \Psi(P) = O, \quad \Psi(D \cap V) = \mathcal{K}_\omega \cap B_O(1), \quad \nabla \Psi(P) = I_{2 \times 2}.
  \end{align*}
\end{defn}

\begin{figure}[t]
  \centering
  \begin{subfigure}[b]{0.27\textwidth}
    \includegraphics[width=\textwidth]{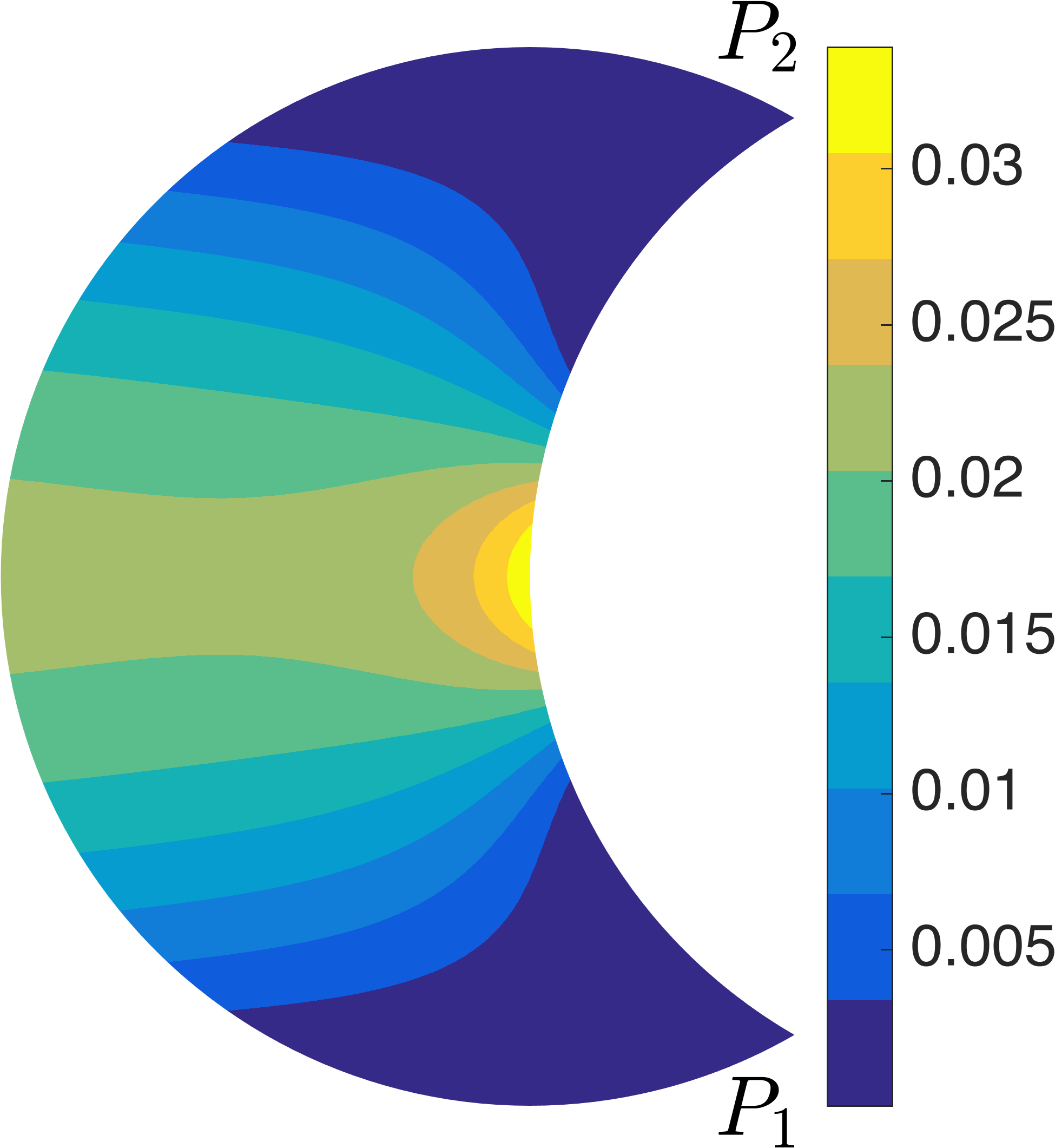}
    \caption{}
  \end{subfigure}
  \hfill
  \begin{subfigure}[b]{0.27\textwidth}
    \includegraphics[width=\textwidth]{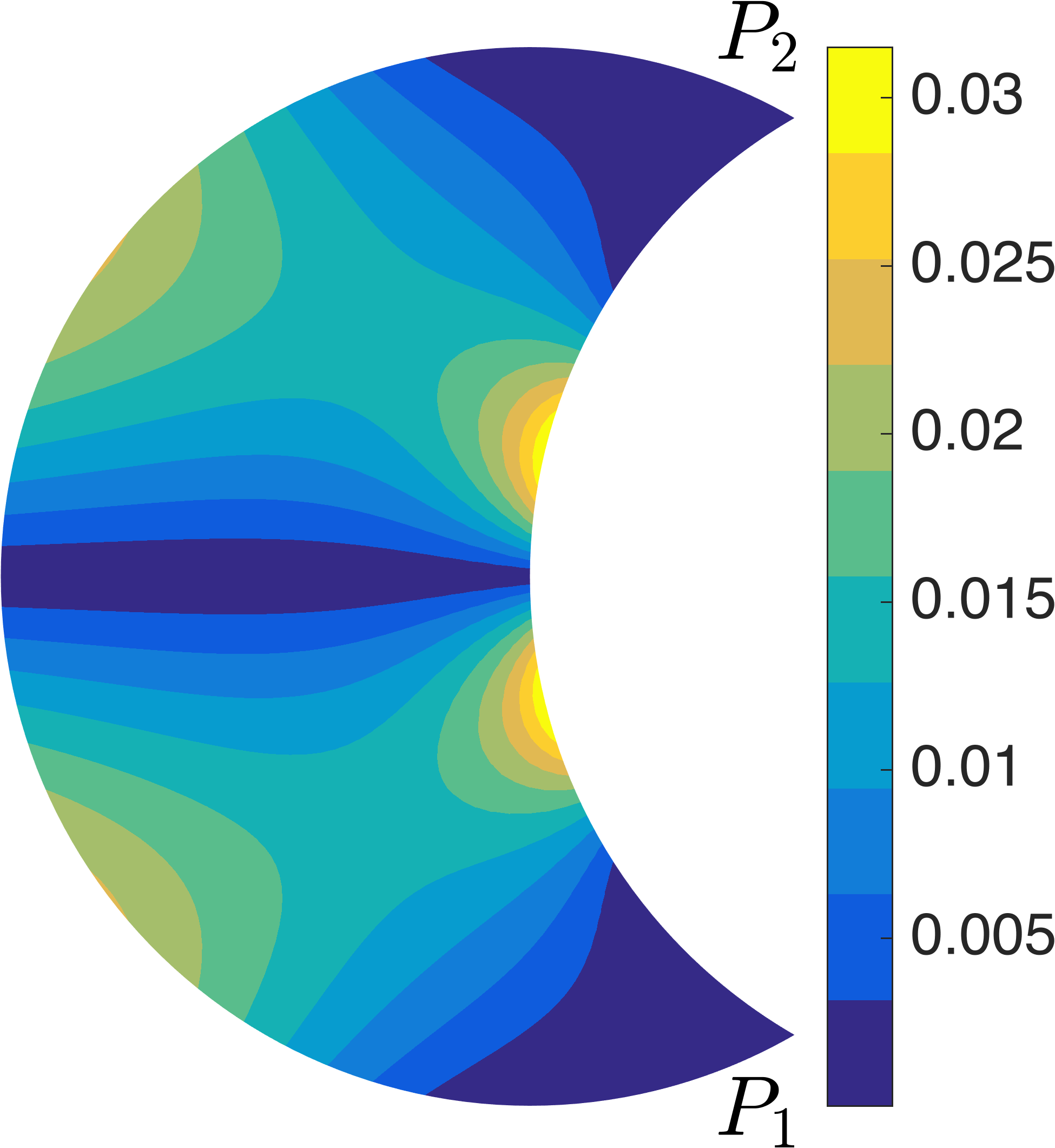}
    \caption{}
  \end{subfigure}
  \hfill
  \begin{subfigure}[b]{0.27\textwidth}
    \includegraphics[width=\textwidth]{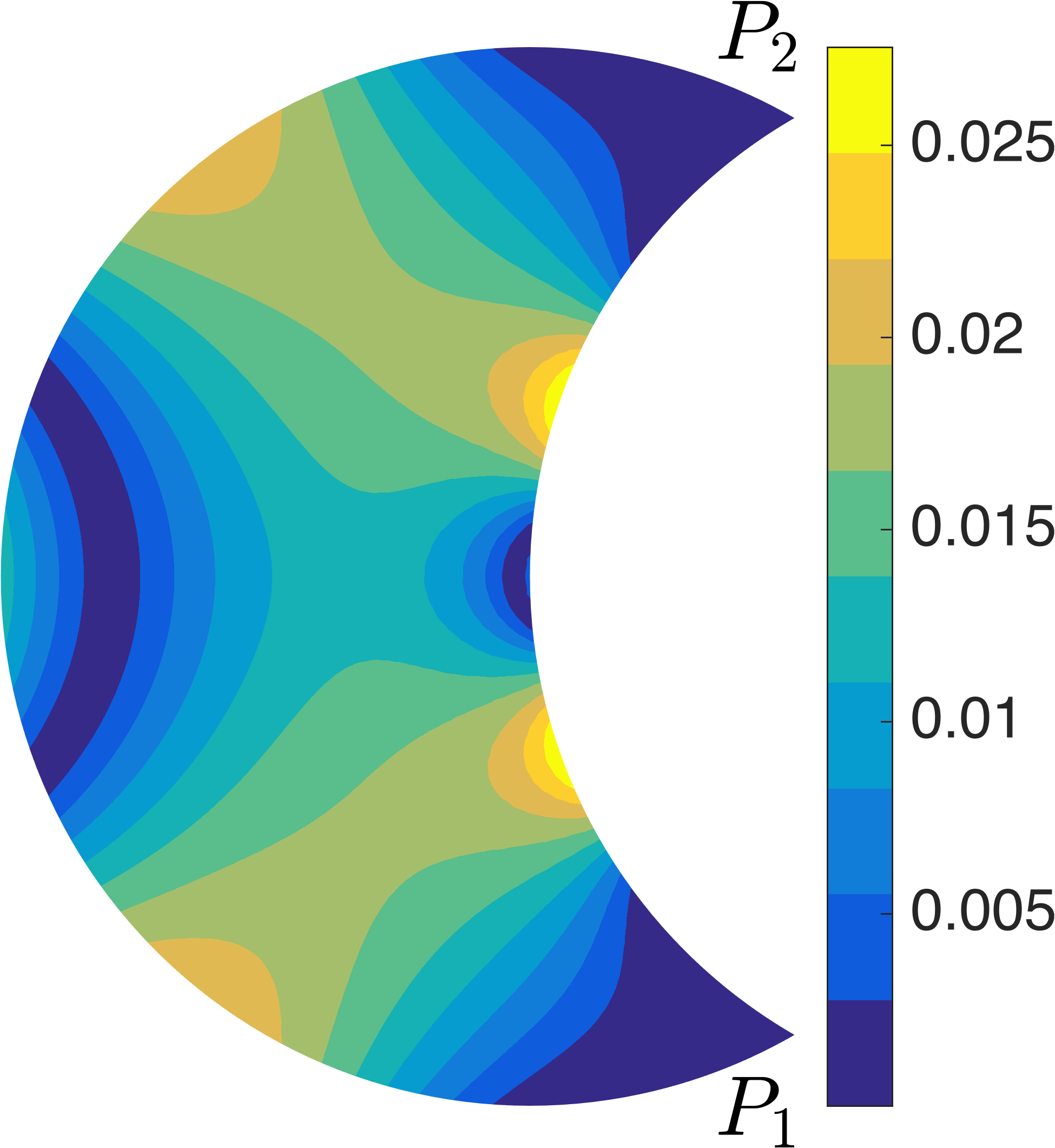}
    \caption{}
  \end{subfigure}
  \caption{The magnitude of transmission eigenfunctions for the moon shaped domain with refractive index $n=16$. (A) $|u_0^{(1)}|$; (B) $|u_0^{(2)}|$; (C) $|u_0^{(3)}|$.}
  \label{fig:moon}
\end{figure}

\begin{figure}[t]
  \centering
  \begin{subfigure}[b]{0.32\textwidth}
    \includegraphics[width=\textwidth]{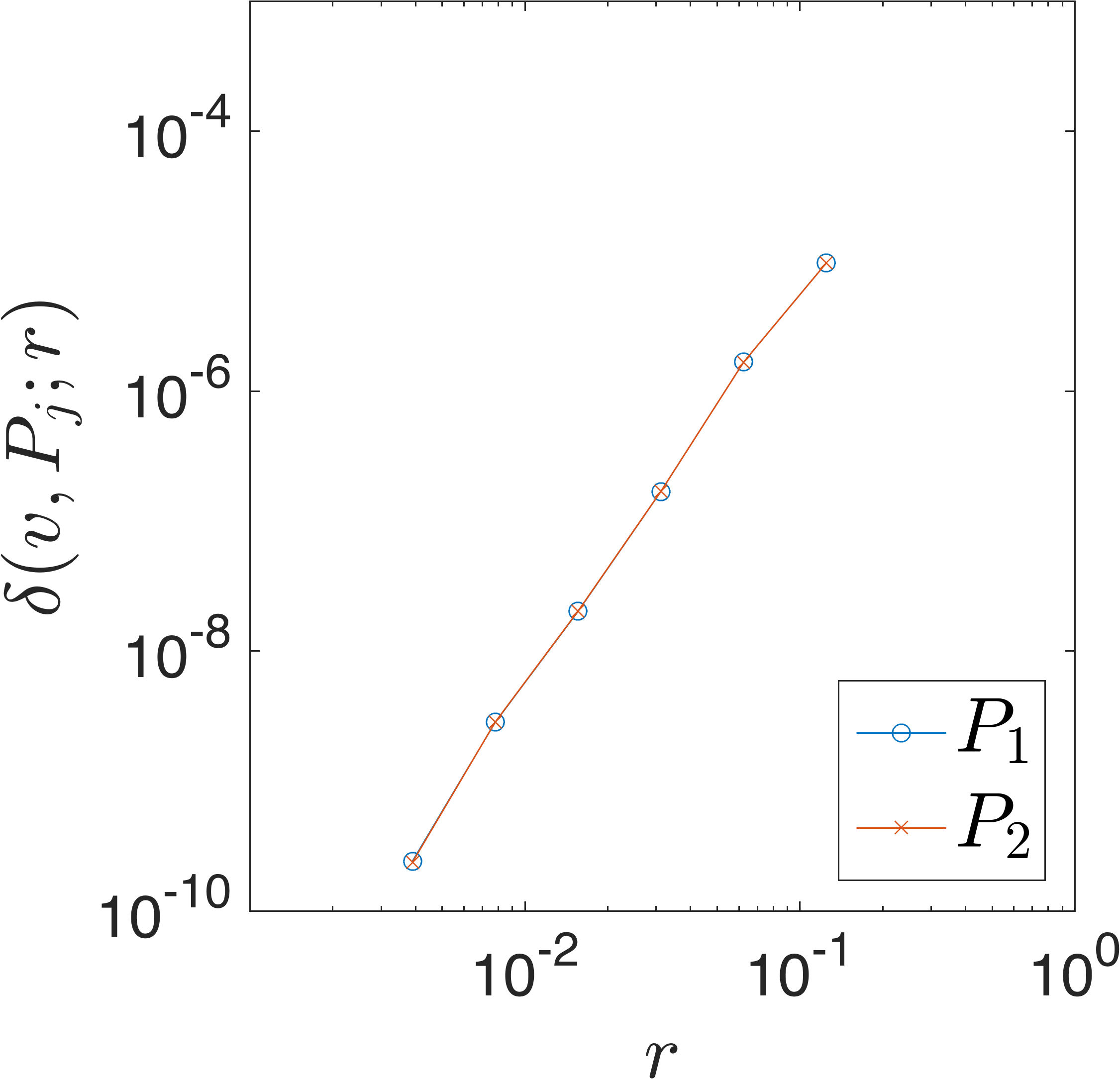}
    \caption{}
  \end{subfigure}
  \hfill
  \begin{subfigure}[b]{0.32\textwidth}
    \includegraphics[width=\textwidth]{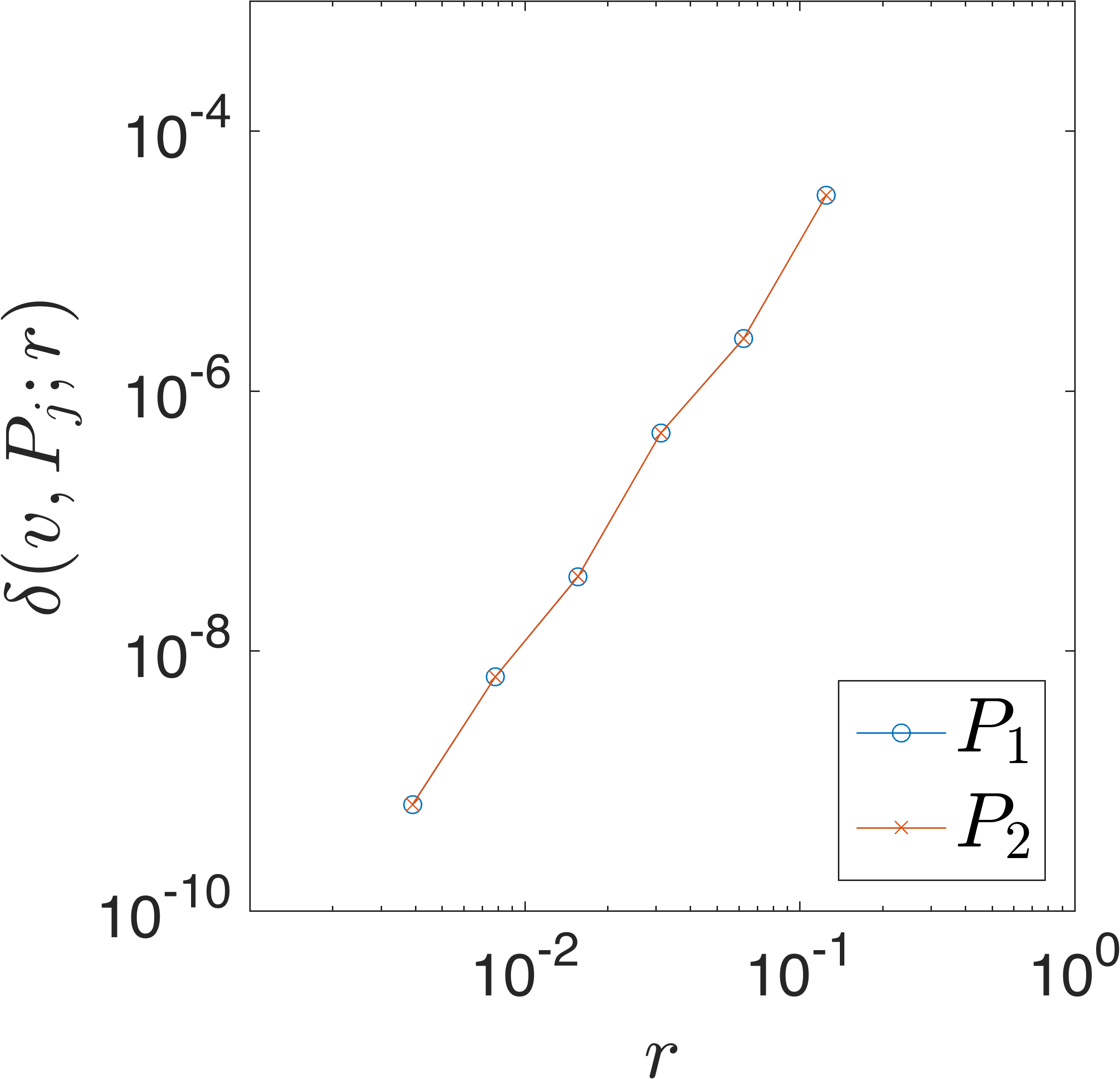}
    \caption{}
  \end{subfigure}
  \hfill
  \begin{subfigure}[b]{0.32\textwidth}
    \includegraphics[width=\textwidth]{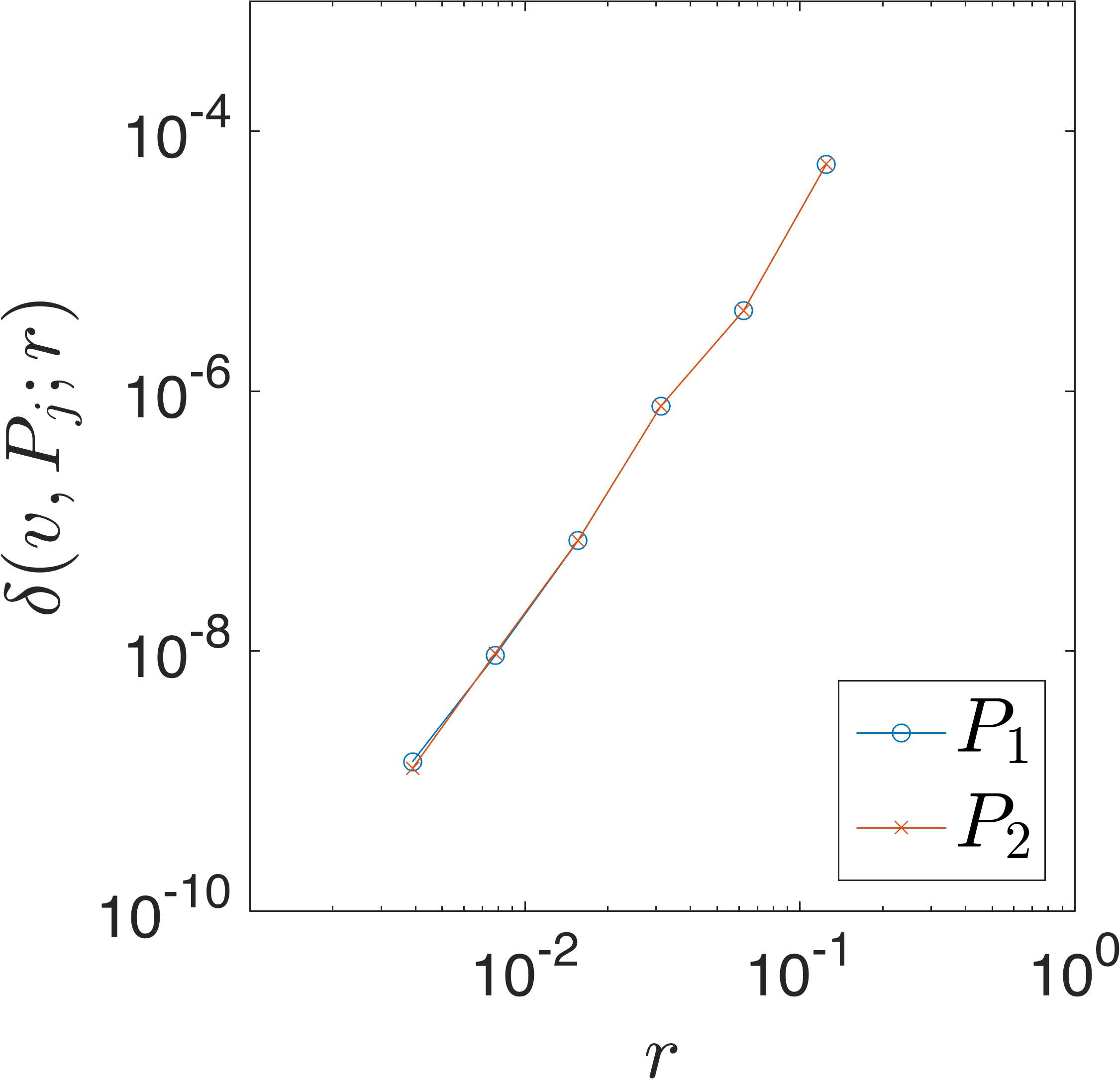}
    \caption{}
  \end{subfigure}
  \caption{The vanishing of the transmission eigenfunctions towards the corner points of the moon shaped domain with refractive index $n=16$. (A) $u_0^{(1)}$; (B) $u_0^{(2)}$; (C) $u_0^{(3)}$.}
  \label{fig:moon_delta}
\end{figure}

Let $D= D_1\setminus D_2$ be a moon shaped domain, where $D_1$ is the disk centered at $(-0.5,0)$ with radius $1$ and $D_2$ is the disk centered at $(0.5,0)$ with radius $1$. Clearly $P_1=(0,-\sqrt{3}/2)$ and $P_2=(0,\sqrt{3}/2)$ are two corner points of $D$ both with opening angle $\pi/3$. The first three real transmission eigenvalues are found to be $1.6832, 1.8153$ and $2.0688$ with corresponding eigenfunctions shown in Figure \ref{fig:moon}. The convergence of the average $L^2$ norm of the corresponding eigenfunctions towards the corner points is shown in Figure \ref{fig:moon_delta} and the order of convergence of shown in Table \ref{tab:moon_order}.
\begin{table}[t]
  \centering
  \begin{tabular}{ccc}
    \toprule
   & $P_1$ & $P_2$ \\
    \midrule
    $u_0^{(1)}$ & 3.06 & 3.06 \\[5pt]
    $u_0^{(2)}$ & 3.07 & 3.07 \\[5pt]
    $u_0^{(3)}$ & 3.03 & 3.05 \\
    \bottomrule
  \end{tabular}
  \caption{Order of convergence $\alpha(u_0^{(i)},P_j)$ of the transmission eigenfunctions for the moon shaped domain with refractive index $n=16$.}
  \label{tab:moon_order}
\end{table}
From Figure \ref{fig:moon} and Figure \ref{fig:moon_delta} we see the vanishing properties of the transmission eigenfunctions are still valid for corner points formed with curved boundary segments. Comparing the order of convergence of the corners $P_1,P_2,P_3$ in Table \ref{tab:equilateral_triangle_order}, corner $P_2$ in Table \ref{tab:right_triangle_order_u0}, corner $P_3$ in Table \ref{tab:arrow_order} and corners $P_1,P_2$ in Table \ref{tab:moon_order} we conjecture the order of convergence $\alpha(v,P)$ for a transmission eigenfunction $v$ at a corner $P$ does not depend on the shape of the corner but only on the opening angle $\omega(P)$.

\subsection{Relation between angle and convergence order}
From Example 4.2 and Example 4.3 we have found that the vanishing convergence order is related to the angle of the corner. In fact, the vanishing convergence rate increases as the sharpness of the corner increases. To find out more on the relation between angle and convergence rate, we perform a series of experiments both for interior angle less than $\pi$ and greater than $\pi$.

\begin{figure}[t]
  \centering
  \begin{subfigure}[b]{0.48\textwidth}
    \includegraphics[width=\textwidth]{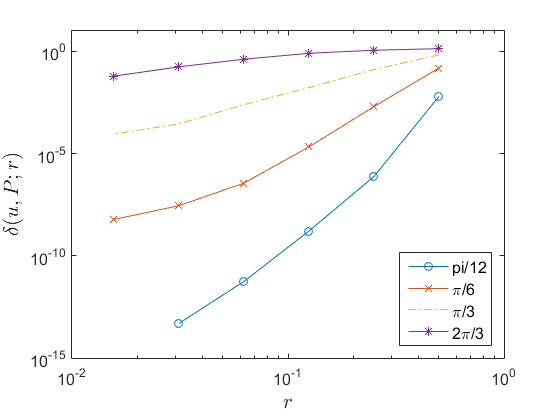}
    \caption{}
  \end{subfigure}
  \hfill
  \begin{subfigure}[b]{0.48\textwidth}
    \includegraphics[width=\textwidth]{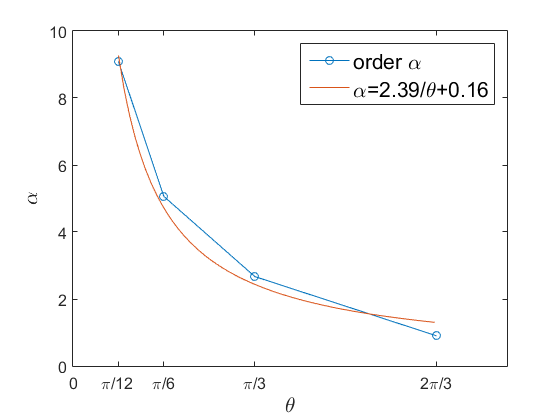}
    \caption{}
  \end{subfigure}
  \caption{(A) The vanishing of the transmission eigenfunctions towards corner of each isosceles triangle with different angles, refractive index $n=16$. (B) Inverse proportional fitting of the vanishing convergence rate of isosceles triangles with different angles.}
  \label{fig:delta_less_pi}
\end{figure}

We consider a series of isosceles triangles of the same height $1$. The angles of the top corners are $\frac{\pi}{12}$, $\frac{\pi}{6}$, $\frac{\pi}{3}$ and $\frac{2\pi}{3}$, respectively. All of the isosceles triangles have the same refractive index $n=16$. We compute the first transmission eigenvalue and corresponding eigenfunction for each triangle. The vanishing rate of the average $L^2$ norm of the eigenfunction near the top corner is shown for different angles in Figure \ref{fig:delta_less_pi}(A).
Fitting the vanishing rate, we have Figure \ref{fig:delta_less_pi}(B), from which we can see that the vanishing rate is inversely proportional to the interior angle when it is less than $\pi$.

\begin{figure}[t]
  \centering
  \begin{subfigure}[b]{0.48\textwidth}
    \includegraphics[width=\textwidth]{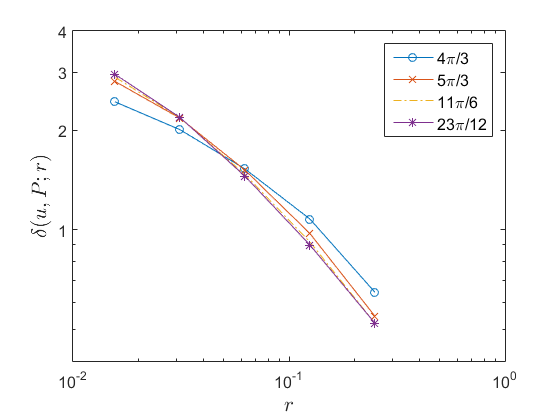}
    \caption{}
  \end{subfigure}
  \hfill
  \begin{subfigure}[b]{0.48\textwidth}
    \includegraphics[width=\textwidth]{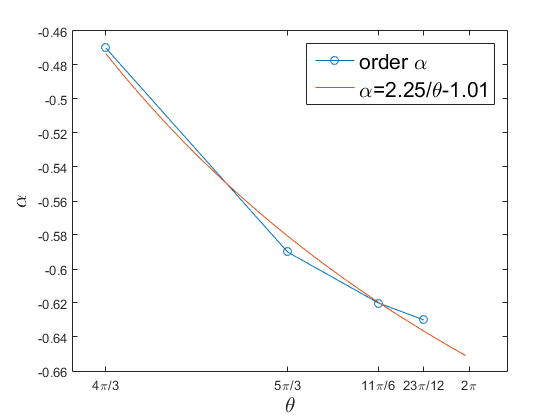}
    \caption{}
  \end{subfigure}
 \caption{(A) The localizing of the transmission eigenfunctions towards corner of each circular sector with different angles, refractive index $n=16$. (B) Inverse proportional fitting of the localizing convergence rate of circular sectors with different angles.}
  \label{fig:delta_gre_pi}
\end{figure}

Furthermore, we perform  experiments to find the relation between the localizing convergence rate and the interior angle when it is greater than $\pi$. Consider a series of circular sectors of radius $1$ and angles $\frac{4\pi}{3}$, $\frac{5\pi}{3}$, $\frac{11\pi}{6}$, $\frac{23\pi}{12}$, respectively. All of the circular sectors have the same refractive index $n=16$. We compute the first transmission eigenvalue and corresponding eigenfunction for each circular sector. The convergence rate of the average $L^2$ norm of the corresponding eigenfunctions towards different corners is shown in Figure \ref{fig:delta_gre_pi}(A).
Fitting the localizing convergence rate, we have Figure \ref{fig:delta_gre_pi}(B), from which we can see that the localizing convergence rate is also inversely proportional to the interior angle when it is greater than $\pi$. The localizing convergence rate increases as the angle of the corner increases. 

\begin{figure}[t]
  \centering
    \includegraphics[width=\textwidth]{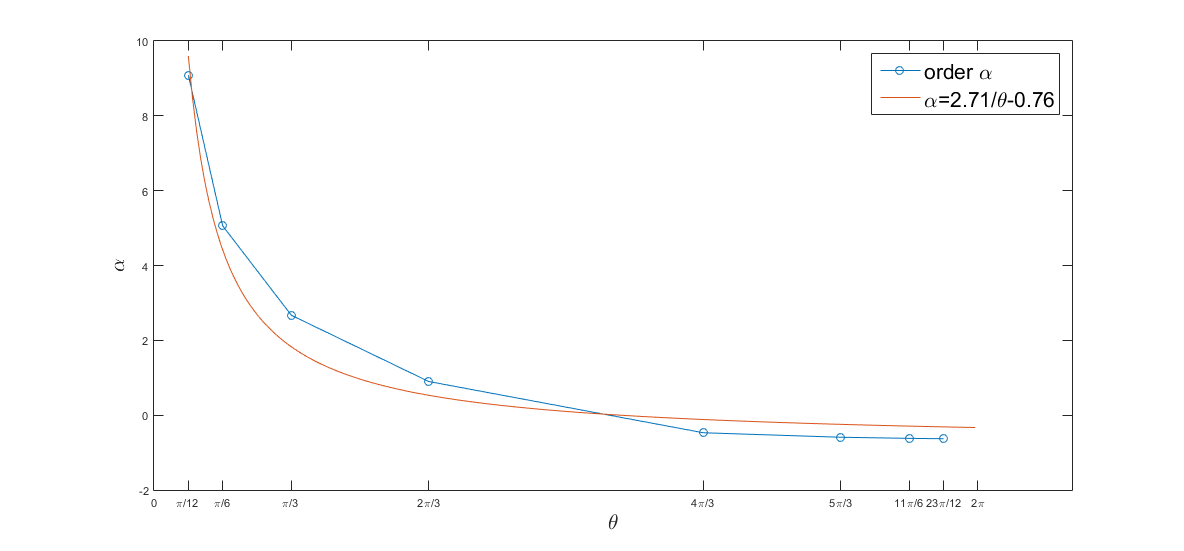}
  \caption{Inverse proportional fitting of the vanishing and localizing convergence rate of circular sectors with different angles.}
  \label{fig:whole}
\end{figure}

Collecting the vanishing convergence rate of angles less than $\pi$ and localization convergence rate of angles greater than $\pi$ together, i.e. fitting the convergence rate for angles range from $\frac{\pi}{12}$ to $\frac{23\pi}{12}$ we have Figure \ref{fig:whole}.
We find that the vanishing and localizing convergence rate is inversely proportional to the angle of corner. What a surprising discovery!

\subsection{On complex transmission eigenvalues and corresponding eigenfunctions}
All of the above examples considered the real transmission eigenvalues. In our numerical calculations, complex transmission eigenvalues also appear. Table \ref{tab:compl} list the first two complex transmission eigenvalues for several domains.
\begin{table}[t]
  \centering
  \begin{tabular}{ccc}
    \toprule
   & $k_1$ & $k_2$ \\
    \midrule
    Unit square & $4.5428\pm 0.5648i$ & $7.0416\pm 0.5847i$ \\[5pt]
    Equilateral triangle &  $3.6021\pm 0.4363i$&$5.5520\pm 0.4633i$ \\[5pt]
    Right triangle & $5.5472\pm 0.6826i$ & $5.5520\pm 0.4633i$ \\[5pt]
    Arrow shaped polygon & $5.1585\pm 0.7016i$ & $7.3080\pm 0.6746i$ \\[5pt]
    Moon shaped domain & $3.6210\pm 0.4946i$ & $4.8331\pm 0.4303$ \\[5pt]
    \bottomrule
  \end{tabular}
  \caption{Complex transmission eigenvalues.}
  \label{tab:compl}
\end{table}
The shapes of these domains are those considered in the previous examples. All of those domains have refractive index $n=16$. We consider both angle less than $\pi$ and angle greater than $\pi$. Our numerical experiments show that the complex transmission eigenvalues appear in terms of conjugate pair. The vanishing and localizing property of transmission eigenfunctions are still valid when the eigenvalues are complex. The numerical experiments also show that there exists no purely imaginary transmission eigenvalues when the refractive index is real-valued.

\section{Numerical experiments: three dimension}\label{sec:numerical-experiment-3D}
In this section, we are going to numerically investigate the vanishing and localizing property of eigenfunctions in the three dimensional case.

We shall first give the definition of vertex and edge points in $\mathbb{R}^3$, refering to \cite{MNP00} and the recent work on edge scattering \cite{EH1}.
\begin{defn}
Let $D\in \mathbb{R}^3$ be a bounded open set. A point $P\in\partial D$ is called a vertex if there exists a neighborhood $V$ of $P$, a diffeomorphism $\Psi$ of class $C^2$ and a polyhedral cone $\Pi$ with the vertex at $O$ such that
$$\nabla\Psi(P)=I_{3\times 3}\in\mathbb{R}^3,~\Psi(P)=O,$$
and $\Psi$ maps $V\cap\bar{D}$ onto a neighborhood of $O$ in $\bar{\Pi}$. $P$ is called an edge point of $D$ if
$$\Psi(V\cap D)=(\mathcal{K}_\omega\cap B_1)\times (-1,1)$$
for some $\omega(P)\in (0,2\pi)\backslash{\pi}$.
\end{defn}
We refer both vertex and edge points as the refractive index. In this section, we numerically show the vanishing and localizing property of the interior transmission eigenfunctions and investigate the convergence rates for both the vertices and edges. Similarly to the two-dimensional case, we define the convergence rate $\alpha$ of a vertex, convergence rate $\beta$ of an edge as follows.

Let $P$ be a vertex of $D$ and set
\begin{align}
\delta(v,P;r) = \frac{1}{\sqrt{|D_r(P)|}} \left\|v\right\|_{L^2 \left(D_P(r)\right)},
\label{deltaP-3D}
\end{align}
where $|D_r(P)|$ is the volume of $D_r(P)$ in three dimensions. Define the order of convergence rate $\alpha$ by the following asymptotic behavior of $\delta(v,P;r)$
  \begin{align*}
    \delta(v,P;r) \sim c(v,P) r^{-\alpha(v,P)} \quad {\rm as}~r \to 0
  \end{align*}
  for each $v=u_0^{(i)}$ or $v=u^{(i)}$, where $c>0$ and $\alpha>0$ are constants independent of $r$.

Let $E$ be an edge of $D$ and $B_r(E)$ be the (curved) cylinder with radius $r$ and with center line to be edge $E$. Let $S_r(E) = D \cap B_r(E)$. Define the average $L^2$ norm of a transmission eigenfunction $v$ in $S_r(E)$ as
\begin{align}
\delta(v,E;r) = \frac{1}{\sqrt{|S_r(E)|}} \left\|v\right\|_{L^2 \left(S_r(E)\right)},
\label{deltaE-3D}
\end{align}
where $|S_r(E)|$ is the volume of $S_r(E)$. Define the order of convergence rate $\beta$ by the following asymptotic behavior of $\delta(v,E;r)$
  \begin{align*}
    \delta(v,E;r) \sim c(v,E) r^{-\beta(v,E)} \quad {\rm as}~r \to 0
  \end{align*}
  for each $v=u_0^{(i)}$ or $v=u^{(i)}$, where $c>0$ and $\beta>0$ are constants independent of $r$. We would like to show that $\delta(v,P;r) \to 0$ and $\delta(v,E;r) \to 0$ as $r \to 0$ and estimate the rate of convergence.

\subsection{Example: unit cube}
\begin{figure}[t]
  \centering
   \begin{subfigure}[b]{0.32\textwidth}
      \includegraphics[width=\textwidth]{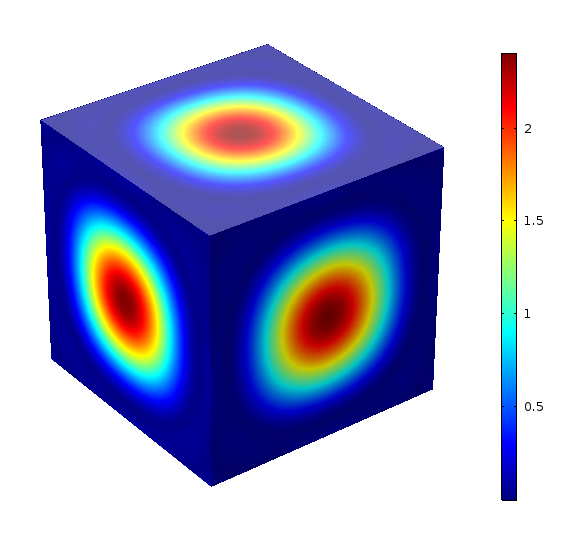}
      \caption{}
    \end{subfigure}
    \hfill
    \begin{subfigure}[b]{0.32\textwidth}
      \includegraphics[width=\textwidth]{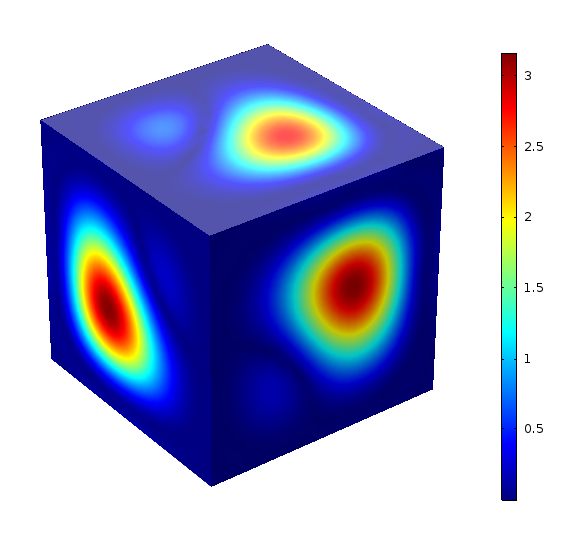}
      \caption{}
    \end{subfigure}
    \hfill
    \begin{subfigure}[b]{0.32\textwidth}
      \includegraphics[width=\textwidth]{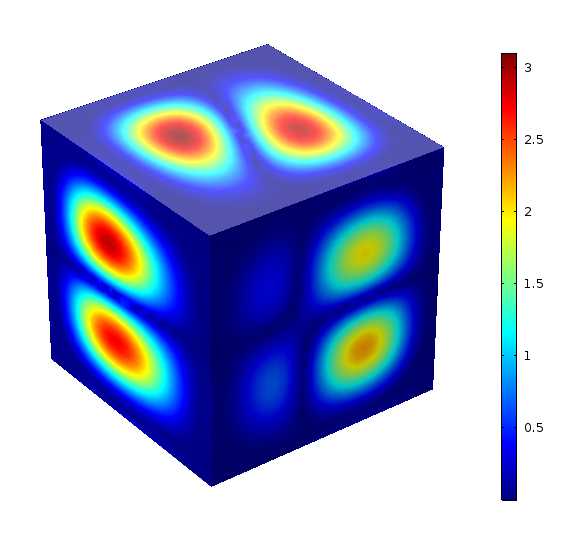}
      \caption{}
    \end{subfigure}
    \hfill
  \begin{subfigure}[b]{0.32\textwidth}
    \includegraphics[width=\textwidth]{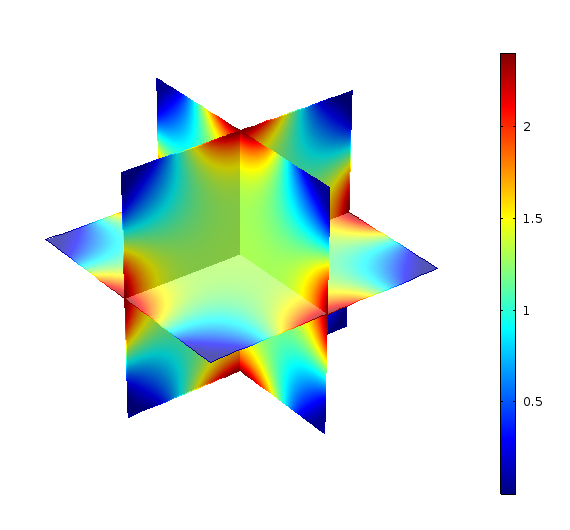}
    \caption{}
  \end{subfigure}
  \hfill
  \begin{subfigure}[b]{0.32\textwidth}
    \includegraphics[width=\textwidth]{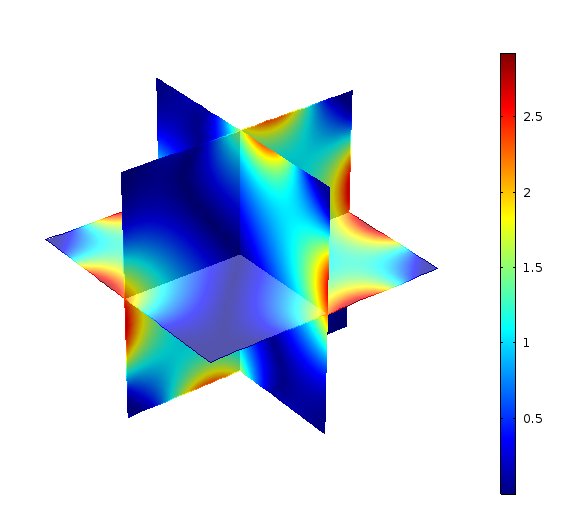}
    \caption{}
  \end{subfigure}
  \hfill
  \begin{subfigure}[b]{0.32\textwidth}
    \includegraphics[width=\textwidth]{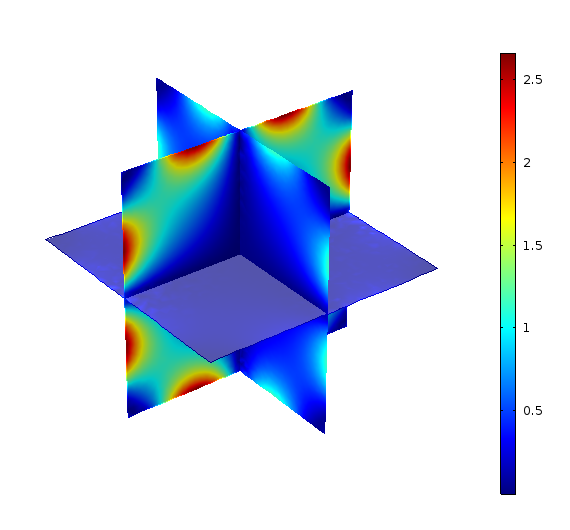}
    \caption{}
  \end{subfigure}
  \hfill
  \begin{subfigure}[b]{0.32\textwidth}
    \includegraphics[width=\textwidth]{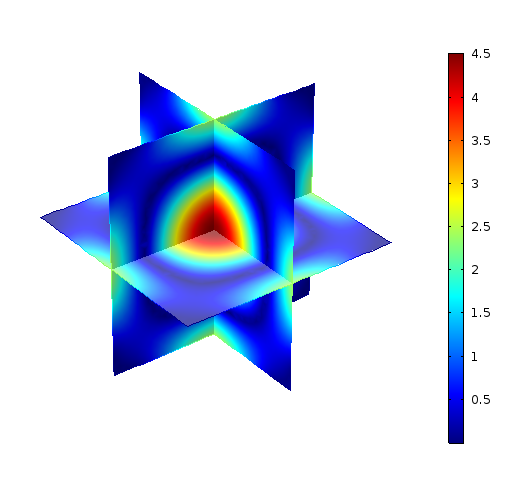}
    \caption{}
  \end{subfigure}
  \hfill
  \begin{subfigure}[b]{0.32\textwidth}
    \includegraphics[width=\textwidth]{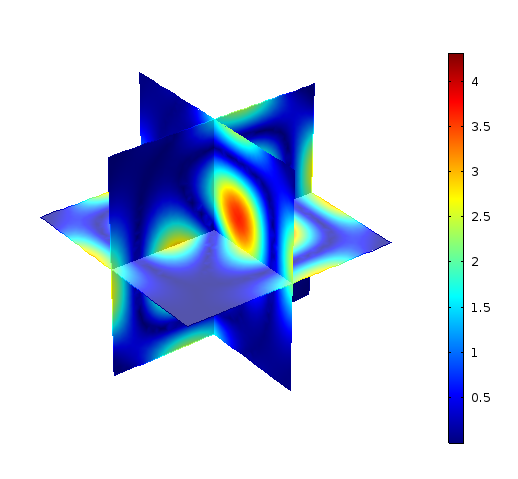}
    \caption{}
  \end{subfigure}
  \hfill
  \begin{subfigure}[b]{0.32\textwidth}
    \includegraphics[width=\textwidth]{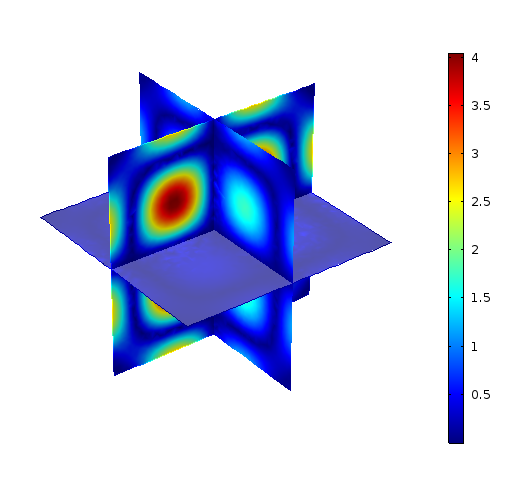}
    \caption{}
  \end{subfigure}
  \caption{The magnitude of transmission eigenfunctions for the unit cube with refractive index $n=16$. (A) surface of $|u_0^{(1)}|$ and $|u^{(1)}|$; (B) surface of $|u_0^{(2)}|$ and $|u^{(2)}|$; (C) surface of $|u_0^{(3)}|$ and $|u^{(3)}|$; (D) $|u_0^{(1)}|$; (E) $|u_0^{(2)}|$; (F) $|u_0^{(3)}|$; (G) $|u^{(1)}|$; (H) $|u^{(2)}|$; (I) $|u^{(3)}|$.}
  \label{fig:cube}
\end{figure}

Let $D$ be the unit cube of edge length $1$ centered at $0$. Let the refractive index $n=16$. The first three real transmission eigenvalues are numerically computed to be $k_1=2.0671$, $k_2=2.5860$ and $k_3=2.9884$. Denote by $u_0^{(i)}$ and $u^{(i)}$ the $i$-th transmission eigenfunctions corresponding to the $i$-th transmission eigenvalue. The magnitude of the transmission eigenfunctions for the first three eigenvalues are shown in Figure \ref{fig:cube}.
From Figure \ref{fig:cube} we see that the transmission eigenfunctions vanish both on the vertices and the edges.

\begin{figure}[t]
  \centering
  \begin{subfigure}[b]{0.32\textwidth}
    \includegraphics[width=\textwidth]{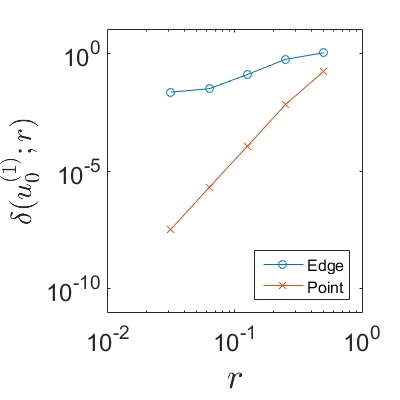}
    \caption{}
  \end{subfigure}
  \hfill
  \begin{subfigure}[b]{0.32\textwidth}
    \includegraphics[width=\textwidth]{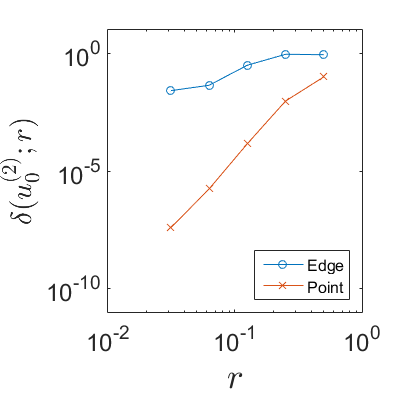}
    \caption{}
  \end{subfigure}
  \hfill
  \begin{subfigure}[b]{0.32\textwidth}
    \includegraphics[width=\textwidth]{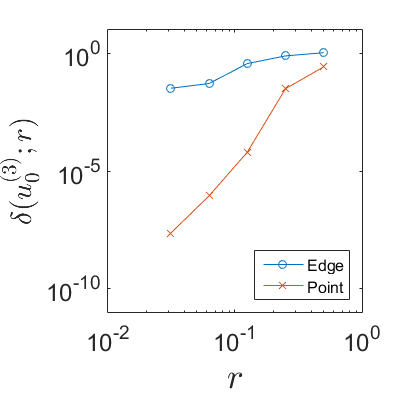}
    \caption{}
  \end{subfigure}
  \hfill
  \begin{subfigure}[b]{0.32\textwidth}
    \includegraphics[width=\textwidth]{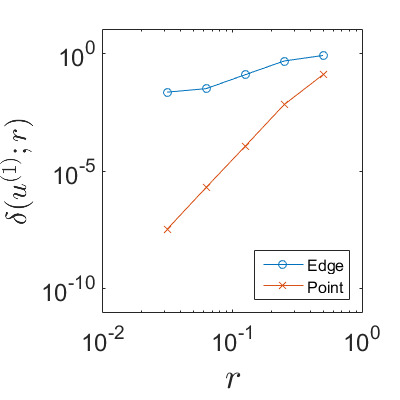}
    \caption{}
  \end{subfigure}
  \hfill
  \begin{subfigure}[b]{0.32\textwidth}
    \includegraphics[width=\textwidth]{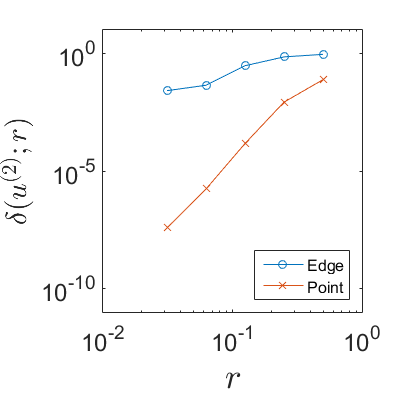}
    \caption{}
  \end{subfigure}
  \hfill
  \begin{subfigure}[b]{0.32\textwidth}
    \includegraphics[width=\textwidth]{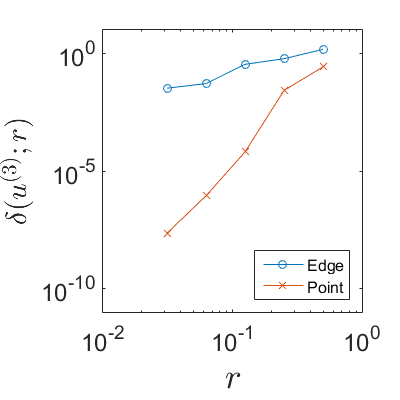}
    \caption{}
  \end{subfigure}
  \caption{The vanishing of the transmission eigenfunctions towards the vertex $P$ and the edge $E$ of the unit cube with refractive index $n=16$. (A) $u_0^{(1)}$; (B) $u_0^{(2)}$; (C) $u_0^{(3)}$; (D) $u^{(1)}$; (E) $u^{(2)}$; (F) $u^{(3)}$.}
  \label{fig:cube_delta}
\end{figure}

Let $P$ be the vertex $(-0.5,-0.5,-0.5)$ and $E$ be the edge connected by vertices $(-0.5,-0.5,-0.5)$ and $(-0.5,-0.5,0.5)$. We discretize $r$ in \eqref{deltaP-3D} and in \eqref{deltaE-3D} into $1/2$, $1/4$, $1/8$, $1/16$, $1/32$, simultaneously near point $P$ and edge $E$. In Figure \ref{fig:cube_delta} we plot $\delta(u_0^{(j)};r)$ and $\delta(u^{(j)}; r)$, $j=1,2,3$ versus $r$ for vertices and edges.
From Figure \ref{fig:cube_delta} we can see that each transmission eigenfunction vanishes near edges and vertices in the sense that $\delta(v,E;r) \to 0$ and $\delta(v,P;r) \to 0$ as $r \to 0$. We can also see that the convergence rate for vertices is faster than that for the edges. Table \ref{tab:cube_order} lists the order of convergence.
\begin{table}[t]
    \centering
    \begin{tabular}{ccccccc}
      \toprule
      &$ u_0^{(1)}$ & $u_0^{(2)}$ & $u_0^{(3)}$ & $u^{(1)}$& $u^{(2)}$& $u^{(3)}$ \\
      \midrule
   Vertex: $\alpha $&  5.64 & 5.51 & 6.21 & 5.54 & 5.41 & 6.19\\
      Edge: $\beta$ &1.53 & 1.45 & 0.40& 1.43 & 1.43 & 1.45\\
      \bottomrule
    \end{tabular}
  \caption{Order of convergence $\alpha(v,P)$ and $\beta(v,E)$ of the transmission eigenfunctions for the unit cube with refractive index $n=16$.}
  \label{tab:cube_order}
\end{table}
The order of convergence in Table \ref{tab:cube_order} is obtained by fitting the data in Figure \ref{fig:cube_delta} by linear polynomials. From this result we conjecture that the convergence rate of vertices is faster than that of edges.

\subsection{Example: Nut-shaped domain}
\begin{figure}[t]
  \centering
   \begin{subfigure}[b]{0.32\textwidth}
      \includegraphics[width=\textwidth]{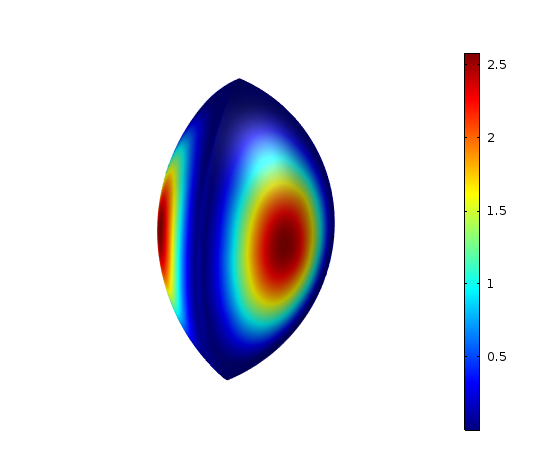}
      \caption{}
    \end{subfigure}
    \hfill
    \begin{subfigure}[b]{0.32\textwidth}
      \includegraphics[width=\textwidth]{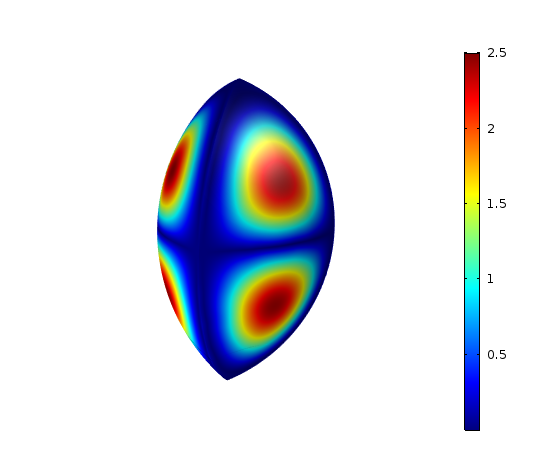}
      \caption{}
    \end{subfigure}
    \hfill
    \begin{subfigure}[b]{0.32\textwidth}
      \includegraphics[width=\textwidth]{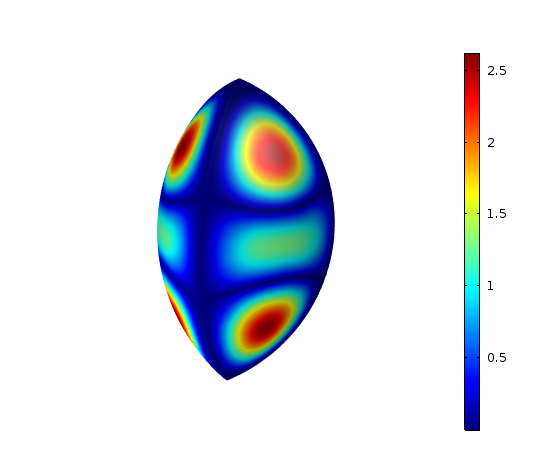}
      \caption{}
    \end{subfigure}
    \hfill
  \begin{subfigure}[b]{0.32\textwidth}
    \includegraphics[width=\textwidth]{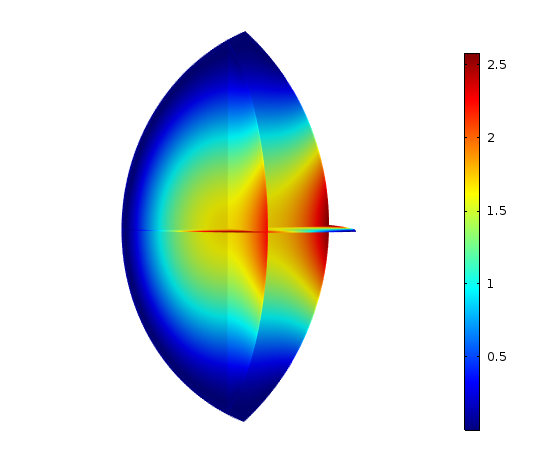}
    \caption{}
  \end{subfigure}
  \hfill
  \begin{subfigure}[b]{0.32\textwidth}
    \includegraphics[width=\textwidth]{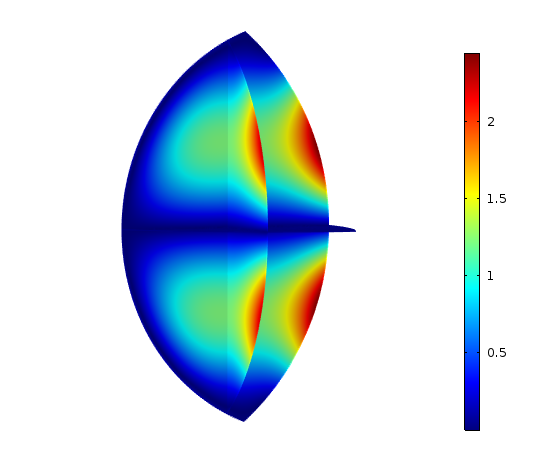}
    \caption{}
  \end{subfigure}
  \hfill
  \begin{subfigure}[b]{0.32\textwidth}
    \includegraphics[width=\textwidth]{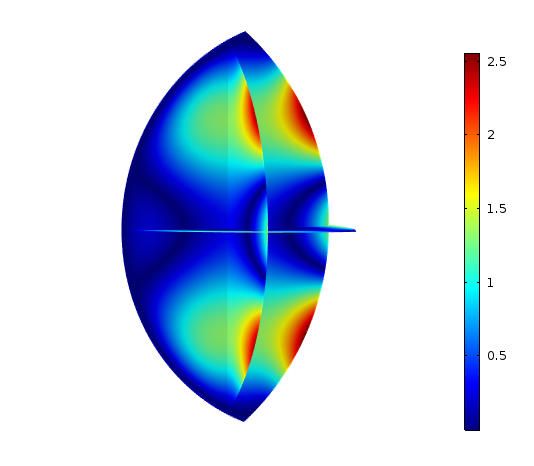}
    \caption{}
  \end{subfigure}
  \hfill
  \begin{subfigure}[b]{0.32\textwidth}
    \includegraphics[width=\textwidth]{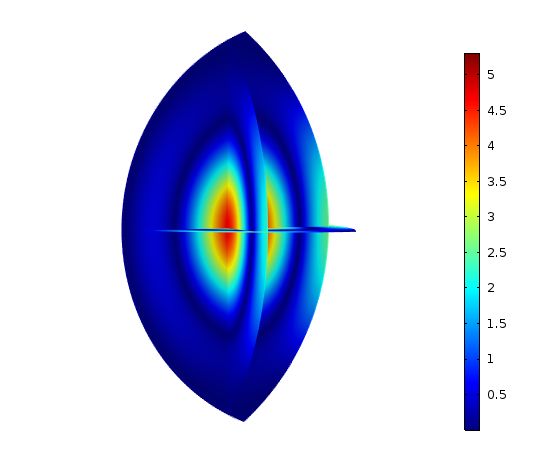}
    \caption{}
  \end{subfigure}
  \hfill
  \begin{subfigure}[b]{0.32\textwidth}
    \includegraphics[width=\textwidth]{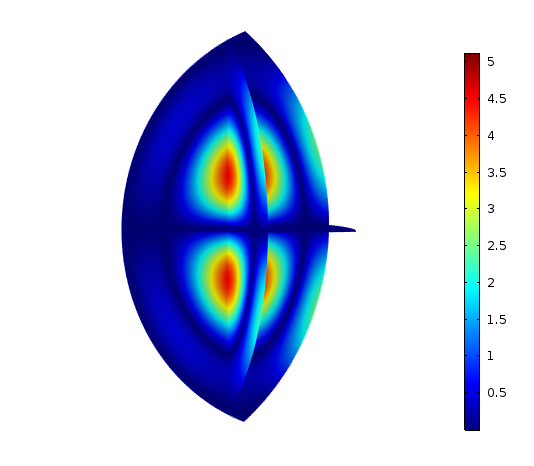}
    \caption{}
  \end{subfigure}
  \hfill
  \begin{subfigure}[b]{0.32\textwidth}
    \includegraphics[width=\textwidth]{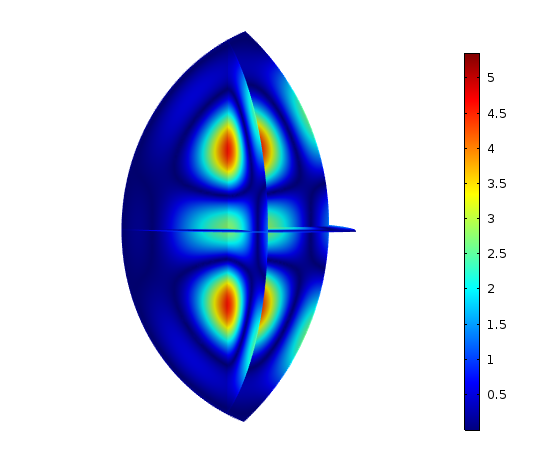}
    \caption{}
  \end{subfigure}
 \caption{The magnitude of transmission eigenfunctions for the nut shaped polyhedron with refractive index $n=16$. (A) surface of $|u_0^{(1)}|$ and $|u^{(1)}|$; (B) surface of $|u_0^{2)}|$ and $|u^{(2)}|$; (C) surface of $|u_0^{(3)}|$ and $|u^{(3)}|$; (D) $|u_0^{(1)}|$; (E) $|u_0^{(2)}|$; (F) $|u_0^{(3)}|$; (G) $|u^{(1)}|$; (H) $|u^{(2)}|$; (I) $|u^{(3)}|$.}
  \label{fig:spheres}
\end{figure}

In this example, we consider domain with curved angles. Let $D$ be the intersection of three balls. The three balls are all of radius $0.8$ and centred at $(0,1/2,0)$, $(-\sqrt{3}/4,-1/4,0)$ and $(\sqrt{3}/4,-1/4,0)$, respectively.  Let the refractive index of $D$ be $n=16$. The first three real transmission eigenvalues are computed to be $3.2592, 3.6796, 4.1936$ with corresponding eigenfunctions shown in Figure \ref{fig:spheres}. It is easy to see that the transmission eigenfunctions vanish on all the vertexes and all the edges.

\begin{figure}[t]
  \centering
  \begin{subfigure}[b]{0.32\textwidth}
    \includegraphics[width=\textwidth]{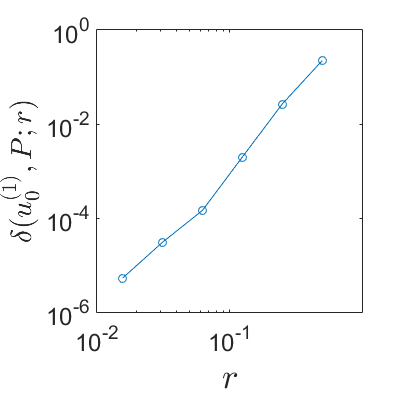}
    \caption{}
  \end{subfigure}
  \hfill
  \begin{subfigure}[b]{0.32\textwidth}
    \includegraphics[width=\textwidth]{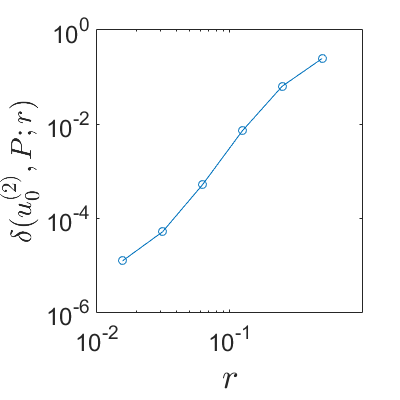}
    \caption{}
  \end{subfigure}
  \hfill
  \begin{subfigure}[b]{0.32\textwidth}
    \includegraphics[width=\textwidth]{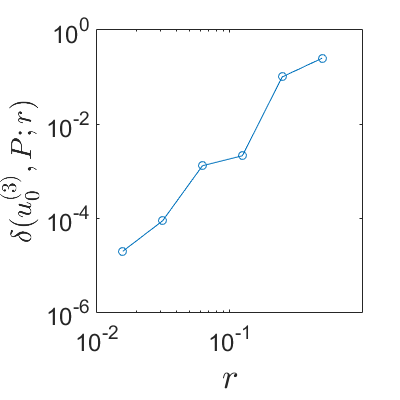}
    \caption{}
  \end{subfigure}
  \hfill
  \begin{subfigure}[b]{0.32\textwidth}
    \includegraphics[width=\textwidth]{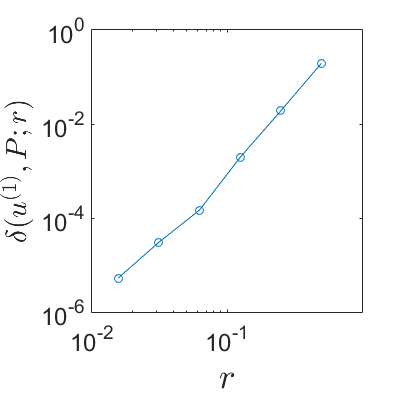}
    \caption{}
  \end{subfigure}
  \hfill
  \begin{subfigure}[b]{0.32\textwidth}
    \includegraphics[width=\textwidth]{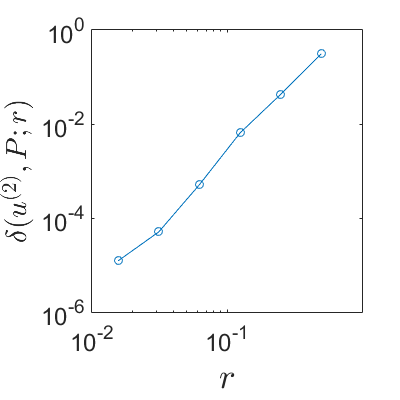}
    \caption{}
  \end{subfigure}
  \hfill
  \begin{subfigure}[b]{0.32\textwidth}
    \includegraphics[width=\textwidth]{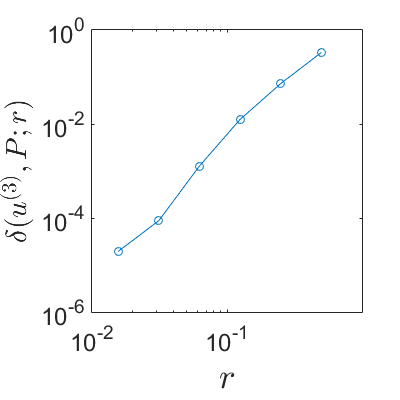}
    \caption{}
  \end{subfigure}
  \caption{The vanishing of the transmission eigenfunctions towards the corner point $P$ of the nut shaped polyhedron with refractive index $n=16$. (A) $u_0^{(1)}$; (B) $u_0^{(2)}$; (C) $u_0^{(3)}$; (D) $u^{(1)}$; (E) $u^{(2)}$; (F) $u^{(3)}$.}
  \label{fig:spheres_delta}
\end{figure}

Let $P$ be the top vertex $(0,0,0.6245)$. In Figure \ref{fig:spheres_delta} we plot $\delta(u_0^{j},P;r)$ and $\delta(u^{j}, P; r)$ versus $r$ for $j=1,2,3$.
From Figure \ref{fig:spheres_delta} we can see that each transmission eigenfunction vanishes towards vertex $P$ in the sense that $\delta(v,P;r) \to 0$ as $r \to 0$. Fitting the data in Figure \ref{fig:spheres_delta} by linear polynomials we obtain the estimates of the convergence order shown in Table \ref{tab:spheres_order}.
\begin{table}[t]
    \centering
    \begin{tabular}{ccccccc}
      \toprule
      &$ u_0^{(1)}$ & $u_0^{(2)}$ & $u_0^{(3)}$ & $u^{(1)}$& $u^{(2)}$& $u^{(3)}$ \\
      \midrule
       Vertex & 4.57 & 4.46 &4.27 &4.5 & 4.46 & 4.36\\
      \bottomrule
    \end{tabular}
  \caption{Order of convergence $\alpha(v,P)$ of the transmission eigenfunctions for the nut shaped polyhedron with refractive index $n=16$.}
  \label{tab:spheres_order}
\end{table}
From this result we see that the vanishing properties of the transmission eigenfunctions are still valid for vertex formed with curved boundary segments in three dimension.


\subsection{Example: Cone}
\begin{figure}[t]
  \centering
   \begin{subfigure}[b]{0.32\textwidth}
      \includegraphics[width=\textwidth]{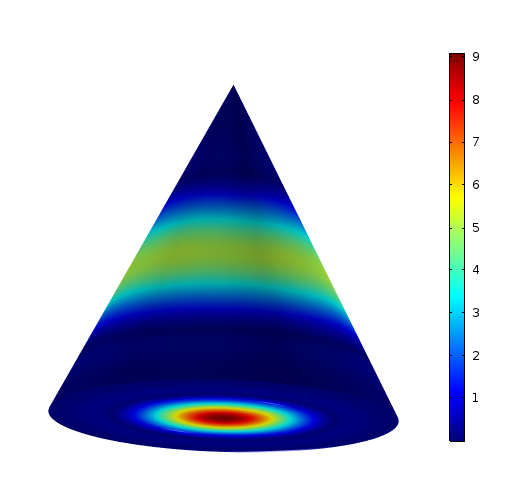}
      \caption{}
    \end{subfigure}
    \hfill
    \begin{subfigure}[b]{0.32\textwidth}
      \includegraphics[width=\textwidth]{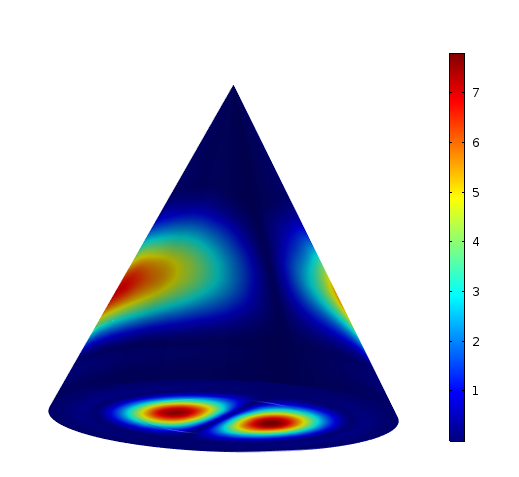}
      \caption{}
    \end{subfigure}
    \hfill
    \begin{subfigure}[b]{0.32\textwidth}
      \includegraphics[width=\textwidth]{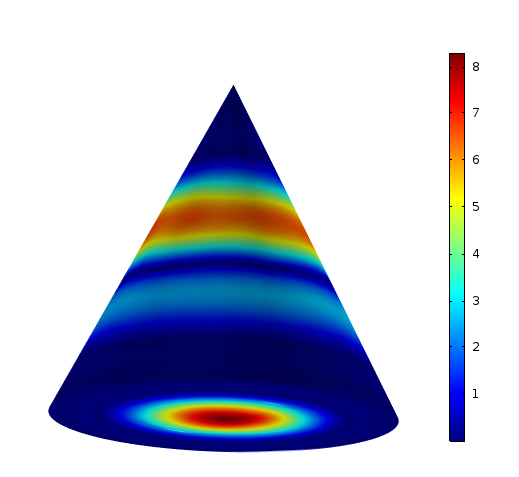}
      \caption{}
    \end{subfigure}
    \hfill
  \begin{subfigure}[b]{0.32\textwidth}
    \includegraphics[width=\textwidth]{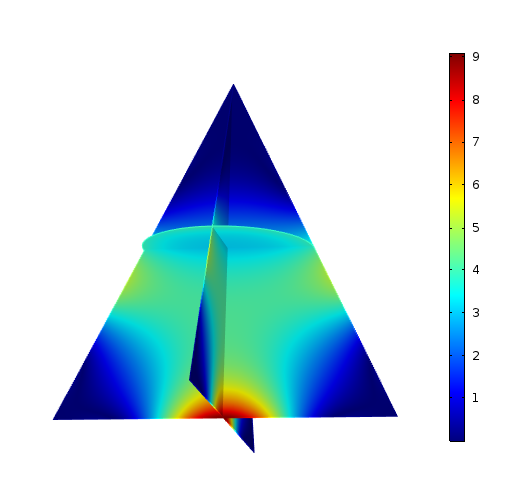}
    \caption{}
  \end{subfigure}
  \hfill
  \begin{subfigure}[b]{0.32\textwidth}
    \includegraphics[width=\textwidth]{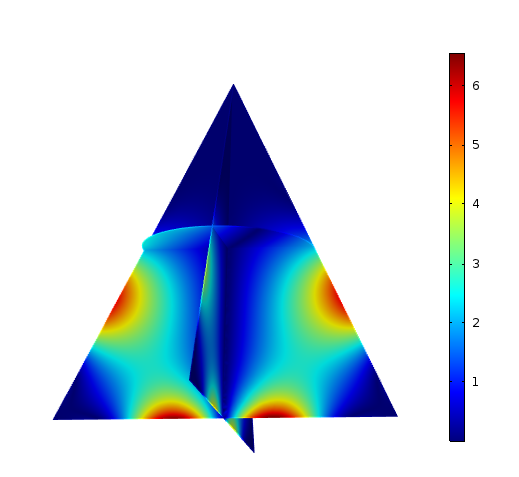}
    \caption{}
  \end{subfigure}
  \hfill
  \begin{subfigure}[b]{0.32\textwidth}
    \includegraphics[width=\textwidth]{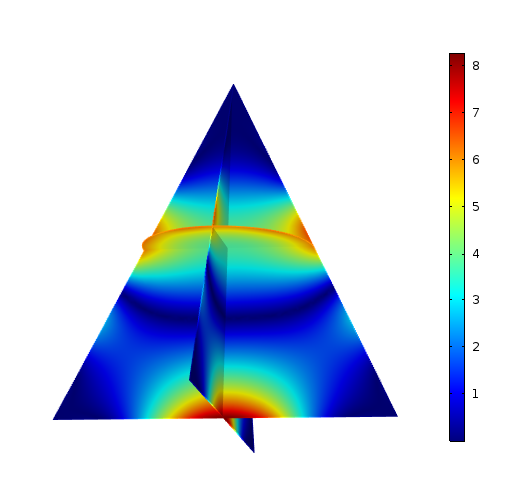}
    \caption{}
  \end{subfigure}
  \hfill
  \begin{subfigure}[b]{0.32\textwidth}
    \includegraphics[width=\textwidth]{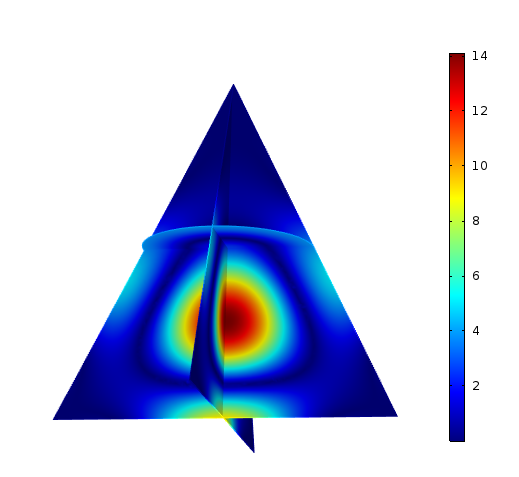}
    \caption{}
  \end{subfigure}
  \hfill
  \begin{subfigure}[b]{0.32\textwidth}
    \includegraphics[width=\textwidth]{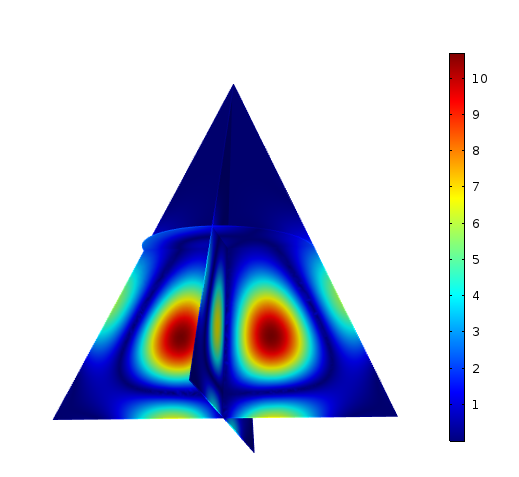}
    \caption{}
  \end{subfigure}
  \hfill
  \begin{subfigure}[b]{0.32\textwidth}
    \includegraphics[width=\textwidth]{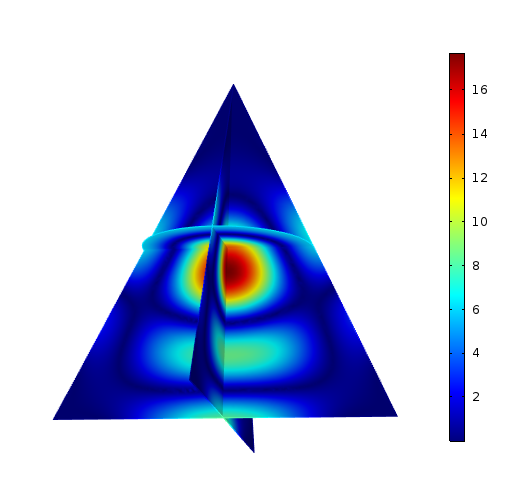}
    \caption{}
  \end{subfigure}
 \caption{The magnitude of transmission eigenfunctions for the cone with refractive index $n=16$. (A) surface of $|u_0^{(1)}|$ and $|u^{(1)}|$; (B) surface of $|u_0^{2)}|$ and $|u^{(2)}|$; (C) surface of $|u_0^{(3)}|$ and $|u^{(3)}|$; (D) $|u_0^{(1)}|$; (E) $|u_0^{(2)}|$; (F) $|u_0^{(3)}|$; (G) $|u^{(1)}|$; (H) $|u^{(2)}|$; (I) $|u^{(3)}|$.}
  \label{fig:cone}
\end{figure}

In this example, we consider a domain with curved edges. Let $D$ be a cone of top point $(0,0,1)$ and bottom radius $1/2$. Let the refractive index of $D$ be $16$. The first three real transmission eigenvalues are computed to be $k_1=3.3668, k_2=4.1259,k_3=4.1710$, where $k_2$ has multiplicity $2$. The corresponding eigenfunctions are shown in Figure \ref{fig:cone}.

\begin{figure}[t]
  \centering
  \begin{subfigure}[b]{0.32\textwidth}
    \includegraphics[width=\textwidth]{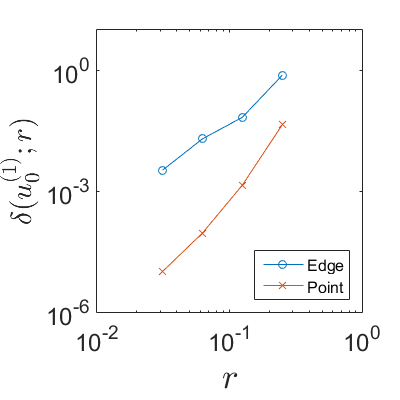}
    \caption{}
  \end{subfigure}
  \hfill
  \begin{subfigure}[b]{0.32\textwidth}
    \includegraphics[width=\textwidth]{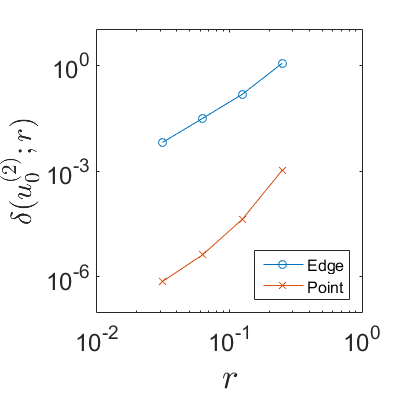}
    \caption{}
  \end{subfigure}
  \hfill
  \begin{subfigure}[b]{0.32\textwidth}
    \includegraphics[width=\textwidth]{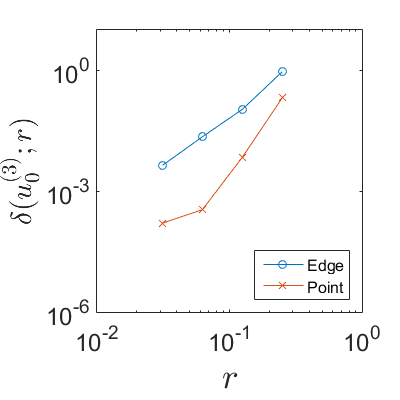}
    \caption{}
  \end{subfigure}
  \hfill
  \begin{subfigure}[b]{0.32\textwidth}
    \includegraphics[width=\textwidth]{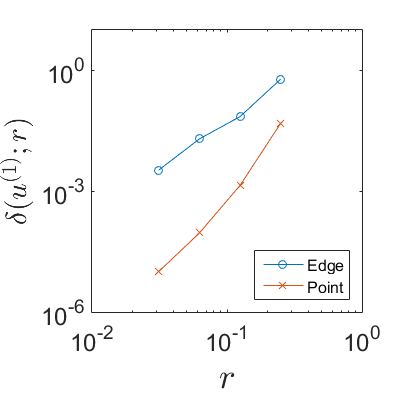}
    \caption{}
  \end{subfigure}
  \hfill
  \begin{subfigure}[b]{0.32\textwidth}
    \includegraphics[width=\textwidth]{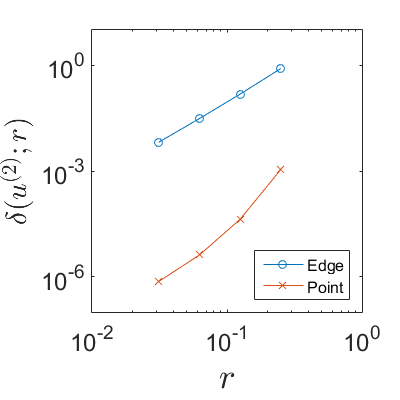}
    \caption{}
  \end{subfigure}
  \hfill
  \begin{subfigure}[b]{0.32\textwidth}
    \includegraphics[width=\textwidth]{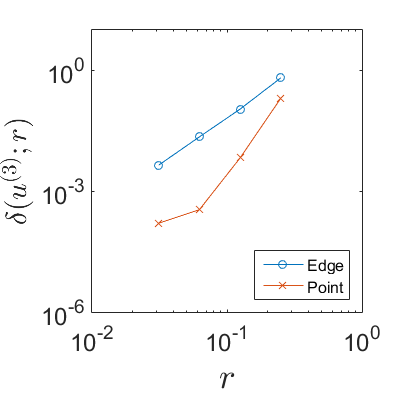}
    \caption{}
  \end{subfigure}
  \caption{The vanishing of the transmission eigenfunctions towards the edge $E$ and vertex $P$ of the cone with refractive index $n=16$. (A) $u_0^{(1)}$; (B) $u_0^{(2)}$; (C) $u_0^{(3)}$; (D) $u^{(1)}$; (E) $u^{(2)}$; (F) $u^{(3)}$.}
  \label{fig:cone_delta}
\end{figure}

Let $P$ be the vertex $(0,0,1)$ and $E$ be edge of bottom circle with radius $1/2$ and centering at $(0,0,0)$. We discretize $r$ in \eqref{deltaP-3D} and in \eqref{deltaE-3D} into $1/2$, $1/4$, $1/8$, $1/16$, $1/32$, simultaneously near edge and vertex. In Figure \ref{fig:cone_delta} we plot $\delta(u_0^{j},E;r)$, $\delta(u^{j}, E; r)$ and $\delta(u_0^{j},P;r)$, $\delta(u^{j}, P; r)$ versus $r$ for $j=1,2,3$.
From Figure \ref{fig:cone_delta} we can see that each transmission eigenfunction vanishes towards vertex $P$ and edge $E$  in the sense that $\delta(v,P;r) \to 0$, $\delta(v,E;r) \to 0$ as $r \to 0$. Fitting the data in Figure \ref{fig:cone_delta} by linear polynomials we obtain the estimates of the convergence order of vertex and edge shown in Table \ref{tab:cone_order}.
\begin{table}[t]
    \centering
        \begin{tabular}{ccccccc}
          \toprule
        &  $ u_0^{(1)}$ & $u_0^{(2)}$ & $u_0^{(3)}$ & $u^{(1)}$& $u^{(2)}$& $u^{(3)}$ \\
          \midrule
         Vertex: $\alpha$ & $$  4.02 & 3.48 &3.53 &4.06 & 3.50&3.51\\
          Edge: $\beta$ & $$  2.53 & 2.47 & 2.54 & 2.44 & 2.31&2.38\\
          \bottomrule
        \end{tabular}
  \caption{Order of convergence $\alpha(v,P)$, $\beta(v,E)$ of the transmission eigenfunctions for the cone with refractive index $n=16$.}
  \label{tab:cone_order}
\end{table}
From this result we see again that the convergence rate of vertex is faster than that of edge.

\subsection{Example: Spherical hat}
\begin{figure}[t]
  \centering
   \begin{subfigure}[b]{0.32\textwidth}
      \includegraphics[width=\textwidth]{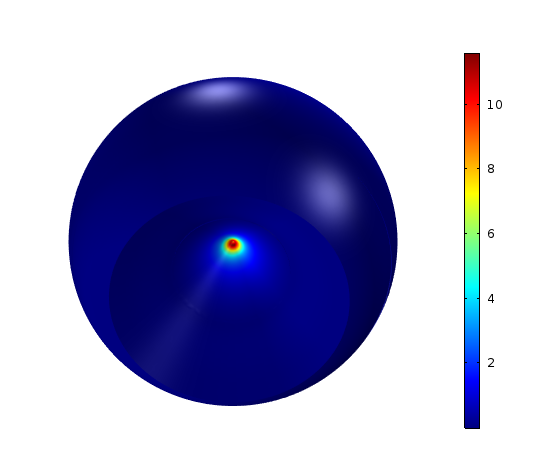}
      \caption{}
    \end{subfigure}
    \hfill
    \begin{subfigure}[b]{0.32\textwidth}
      \includegraphics[width=\textwidth]{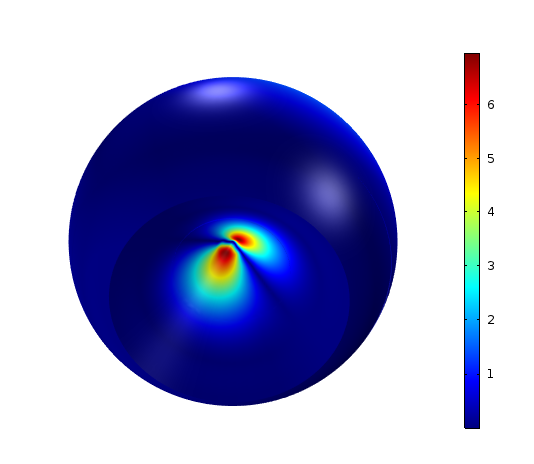}
      \caption{}
    \end{subfigure}
    \hfill
    \begin{subfigure}[b]{0.32\textwidth}
      \includegraphics[width=\textwidth]{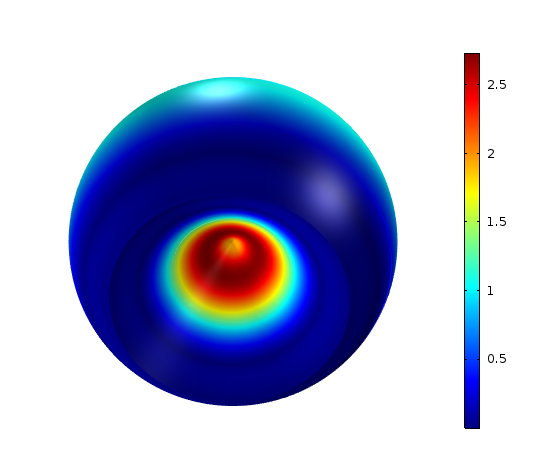}
      \caption{}
    \end{subfigure}
    \hfill
  \begin{subfigure}[b]{0.32\textwidth}
    \includegraphics[width=\textwidth]{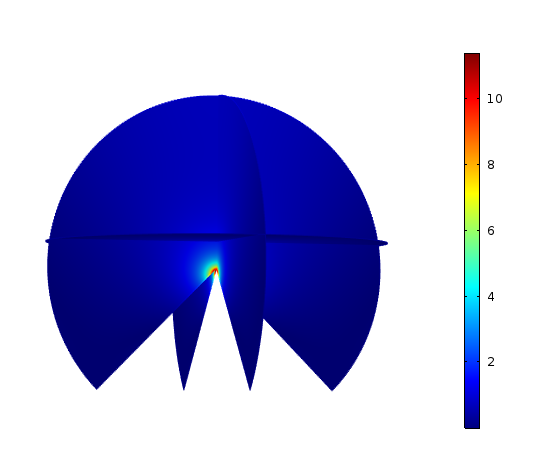}
    \caption{}
  \end{subfigure}
  \hfill
  \begin{subfigure}[b]{0.32\textwidth}
    \includegraphics[width=\textwidth]{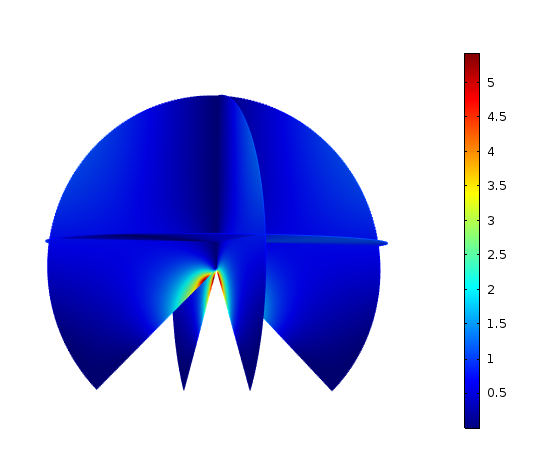}
    \caption{}
  \end{subfigure}
  \hfill
  \begin{subfigure}[b]{0.32\textwidth}
    \includegraphics[width=\textwidth]{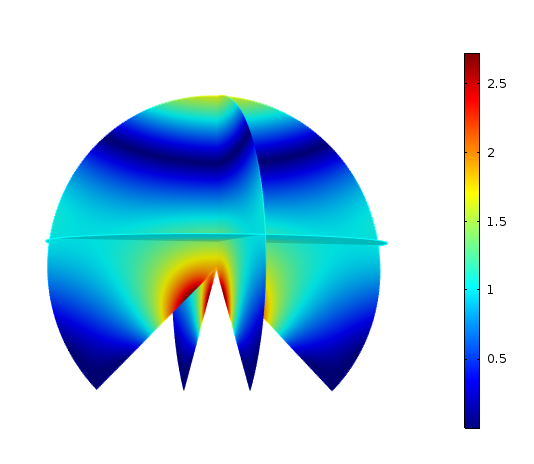}
    \caption{}
  \end{subfigure}
  \hfill
  \begin{subfigure}[b]{0.32\textwidth}
    \includegraphics[width=\textwidth]{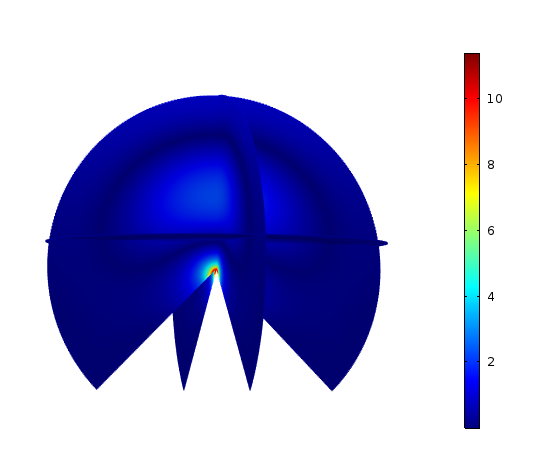}
    \caption{}
  \end{subfigure}
  \hfill
  \begin{subfigure}[b]{0.32\textwidth}
    \includegraphics[width=\textwidth]{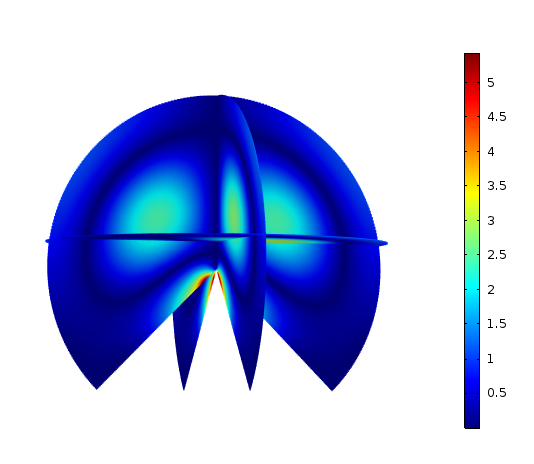}
    \caption{}
  \end{subfigure}
  \hfill
  \begin{subfigure}[b]{0.32\textwidth}
    \includegraphics[width=\textwidth]{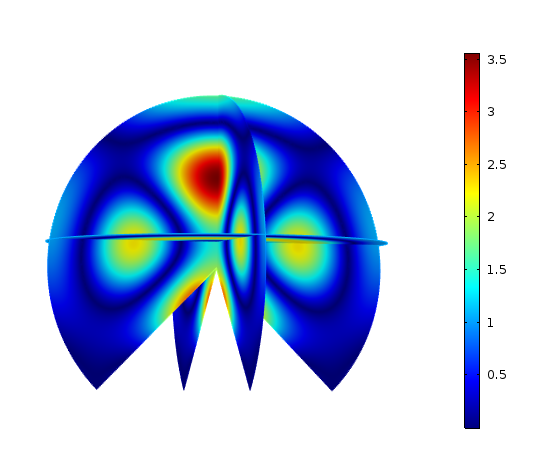}
    \caption{}
  \end{subfigure}
 \caption{The magnitude of transmission eigenfunctions for the spherical hat with refractive index $n=16$. (A) surface of $|u_0^{(1)}|$ and $|u^{(1)}|$; (B) surface of $|u_0^{2)}|$ and $|u^{(2)}|$; (C) surface of $|u_0^{(3)}|$ and $|u^{(3)}|$; (D) $|u_0^{(1)}|$; (E) $|u_0^{(2)}|$; (F) $|u_0^{(3)}|$; (G) $|u^{(1)}|$; (H) $|u^{(2)}|$; (I) $|u^{(3)}|$.}
  \label{fig:sherical_hat}
\end{figure}

Similar to the two dimensional case, we consider the vertex with angles greater than $\pi$. Let $B$ be the unit ball and $C$ be a cone with top point $(0,0,0)$ and bottom radius $1$. We consider the domain $D:=B\backslash C$. Let the refractive index of $D$ be $16$. The first four real transmission eigenvalues are computed to be $k_1=1.5995, k_2=1.6990,k_3=1.8956,k_4=1.9220$, where $k_2,k_3$ has multiplicity $2$. The corresponding eigenfunctions for $k_1,k_2,k_4$ are shown in Figure \ref{fig:sherical_hat}.

\begin{figure}[t]
  \centering
  \begin{subfigure}[b]{0.32\textwidth}
    \includegraphics[width=\textwidth]{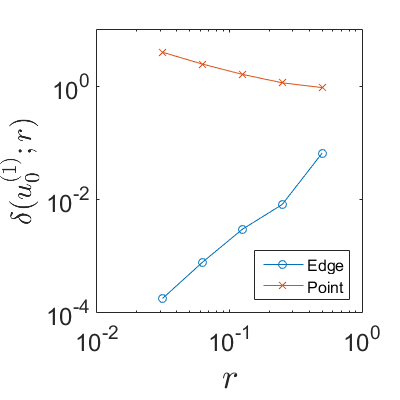}
    \caption{}
  \end{subfigure}
  \hfill
  \begin{subfigure}[b]{0.32\textwidth}
    \includegraphics[width=\textwidth]{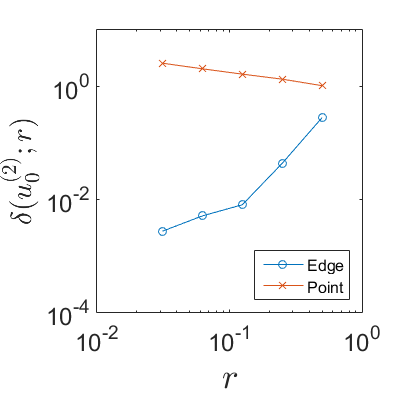}
    \caption{}
  \end{subfigure}
  \hfill
  \begin{subfigure}[b]{0.32\textwidth}
    \includegraphics[width=\textwidth]{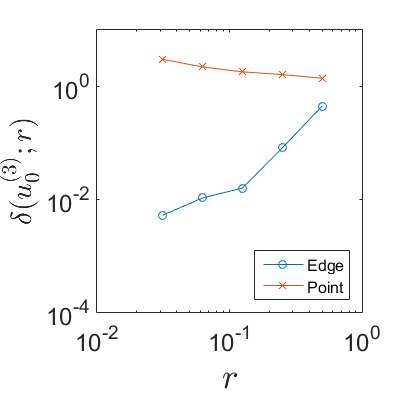}
    \caption{}
  \end{subfigure}
  \hfill
  \begin{subfigure}[b]{0.32\textwidth}
    \includegraphics[width=\textwidth]{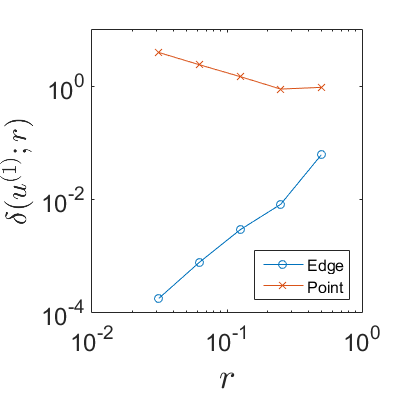}
    \caption{}
  \end{subfigure}
  \hfill
  \begin{subfigure}[b]{0.32\textwidth}
    \includegraphics[width=\textwidth]{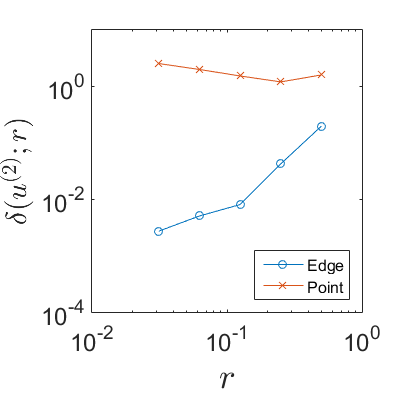}
    \caption{}
  \end{subfigure}
  \hfill
  \begin{subfigure}[b]{0.32\textwidth}
    \includegraphics[width=\textwidth]{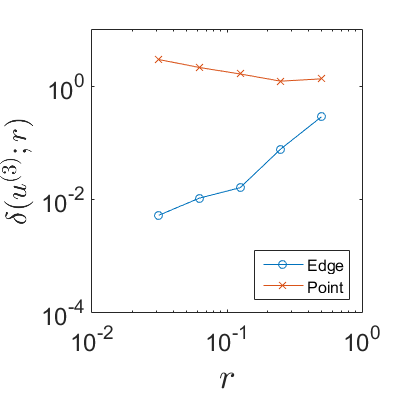}
    \caption{}
  \end{subfigure}
  \caption{The vanishing of the transmission eigenfunctions towards the edge $E$ and vertex $P$ of the spherical hat with refractive index $n=16$. (A) $u_0^{(1)}$; (B) $u_0^{(2)}$; (C) $u_0^{(3)}$; (D) $u^{(1)}$; (E) $u^{(2)}$; (F) $u^{(3)}$.}
  \label{fig:hat_delta}
\end{figure}

Let $P$ be the vertex $(0,0,0)$ and $E$ be the circular bottom edge. We discretize $r$ in \eqref{deltaP-3D} and in \eqref{deltaE-3D} into $1/2$, $1/4$, $1/8$, $1/16$, $1/32$, simultaneously near edge and vertex. In Figure \ref{fig:hat_delta} we plot $\delta(u_0^{j},E;r)$, $\delta(u^{j}, E; r)$ and $\delta(u_0^{j},P;r)$, $\delta(u^{j}, P; r)$ versus $r$ for $j=1,2,3$.
From Figure \ref{fig:hat_delta} we can see that each transmission eigenfunction vanishes towards the edge $E$ in the sense that $\delta(v,E;r) \to 0$ as $r \to 0$, but not for vertex $P$, where $P$ posses angle larger than $\pi$. The order of convergence rate of edge $E$ is listed in Table \ref{tab:hat_order}.
\begin{table}[t]
    \centering
    \begin{tabular}{ccccccc}
      \toprule
    &  $ u_0^{(1)}$ & $u_0^{(2)}$ & $u_0^{(3)}$ & $u^{(1)}$& $u^{(2)}$& $u^{(3)}$ \\
      \midrule
     Edge: $\beta$& $$  2.05 & 1.65 & 1.58 & 2.03 & 1.54&1.45\\
      \bottomrule
    \end{tabular}
  \caption{Order of convergence $\beta(v,E)$ of the transmission eigenfunctions for the spherical hat with refractive index $n=16$.}
  \label{tab:hat_order}
\end{table}

Figure \ref{fig:sherical_hat} and \ref{fig:hat_delta} clearly show the localizing property of transmission eigenfunctions at vertex whose angle is larger than $\pi$. From the two dimensional Figure \ref{fig:arrow}(A-C) and three dimensional Figure \ref{fig:sherical_hat}(A-C), we can see that the localizing behavior depends on eigenfunctions and eigenvalue multiplicity. Furthermore, for the first transmission eigenvalue, the corresponding eigenfunction blows up towards exactly the vertex, while for the eigenvalues posses multiplicity, the corresponding eigenfunctions blow up near the vertex from different directions which is rather complicate.

\section{Concluding remarks}\label{sect:con}

In this work, we numerically invesigate the vanishing and localizing properties of the interior transmission eigenfunctions at singular points on the support of the underlying refractive index, i.e. near points where the boundary tangent is not defined. We numerically show that if the interior angle of a corner is less than $\pi$, then the transmission eigenfunctions vanish near the corner, whereas if the interior angle is bigger than $\pi$, the transmission eigenfunctions localize near the corner. Furthermore, we estimate the order of convergence rate and find that it is related to the angle of the corner. In the three dimensional case, we also present the vanishing property of transmission eigenfunctions on edges. It turns out that edges also posses the vanishing phenomena. In the examples of cube and cone, the results show that the convergence rate of vertex is faster than that of the edge. 
On the one hand, our numerical results clearly verify the theoretical study in \cite{EL} on the vanishing property of transmission eigenfunctions near corners. On the other hand, the numerical results indicate that the geometric properties of the transmission eigenfunctions can be much more delicate and intriguing than the one theoretically justified in \cite{EL}. Our numerical study is by no means exclusive and complete, and it opens up a new research direction for many further developments. 
As one possible application, the vanishing and localizing behaviours of transmission eigenfunctions clearly carry the geometric information of the underlying refractive index $n$, and hence they can be used in inverse scattering problems of recovering the refractive index from exterior measurements. We shall report this finding in a forthcoming paper.

\end{document}